\documentclass[10pt, reqno]{amsart}
\usepackage{graphicx, amssymb, amsmath, amsthm}
\numberwithin{equation}{section}

\usepackage{color}
\usepackage{epsfig}
\usepackage[backref=page]{hyperref}

\hypersetup{urlcolor=blue, citecolor=red}

\usepackage{hyperref}

\let\Re=\undefined\DeclareMathOperator*{\Re}{Re}
\let\Im=\undefined\DeclareMathOperator*{\Im}{Im}

\newcommand{\Z}{\mathbb{Z}}

\newcommand{\HH}{\mathcal{H}}

\newtheorem{theorem}{Theorem}[section]

\newtheorem{lemma}[theorem]{Lemma}

\newtheorem{corollary}[theorem]{Corollary}

\newtheorem{proposition}[theorem]{Proposition}

\theoremstyle{definition}
\newtheorem{definition}[theorem]{Definition}
\newtheorem{remark}[theorem]{Remark}

\makeatletter
\newcommand{\Extend}[5]{\ext@arrow0099{\arrowfill@#1#2#3}{#4}{#5}}
\makeatother

\begin{document}
\title[dispersive estimates on product cones]{Decay estimates for one Aharonov-Bohm solenoid in a uniform magnetic field III: Product cones}

\author{Haoran Wang}
\address{Institute of Mathematics, Henan Academy of Sciences, Zhengzhou 450046, China}
\email{wanghaoran@hnas.ac.cn}

\begin{abstract}
The goal of a recently launched project is to extend the Euclidean models in \cite{Wang24,WZZ25-AHP,WZZ25-JDE} to a more general setting of conically singular spaces. In this paper, the main results include a weighted dispersive inequality for the Schr\"odinger equation and a dispersive estimate for the wave equation both with one Aharonov-Bohm solenoid in a uniform magnetic field on the product cone $X=\mathcal{C}(\mathbb{S}_\sigma^1)=(0,+\infty)_r\times\mathbb{S}_\sigma^1$ endowed with the flat metric $g=dr^2+r^2d\theta^2$, where $\mathbb{S}_\sigma^1\simeq\mathbb{R}/2\pi\sigma\mathbb{Z}$ denotes the circle of radius $\sigma\geq1$ in the Euclidean plane $\mathbb{R}^2$. As a byproduct, we also give the corresponding Strichartz estimates for these equations via the abstract argument of Keel-Tao.
\end{abstract}

\maketitle

\begin{center}
 \begin{minipage}{120mm}
 \end{minipage}
 \end{center}

\section{Introduction}

In the previous papers \cite{Wang24,WZZ25-AHP,WZZ25-JDE}, decay estimates for one Aharonov-Bohm solenoid in a uniform magnetic field have been obtained for solutions to the Schr\"odinger equation on $\mathbb{R}^2$
\begin{equation}\label{eq:S}
\begin{cases}
i\partial_tv(t,x)+H_{\alpha,B_0}v(t,x)=0,\\
v(0,x)=v_0(x)
\end{cases}
\end{equation}
and the wave equation on $\mathbb{R}^2$
\begin{equation}\label{eq:W}
\begin{cases}
\partial_{tt}u(t,x)+H_{\alpha,B_0}u(t,x)=0,\\
u(0,x)=u_0(x),\qquad \partial_tu(0,x)=u_1(x),
\end{cases}
\end{equation}
where $H_{\alpha,B_0}$ denotes the associated magnetic Schr\"odinger operator
\begin{equation}\label{H-A}
H_{\alpha,B_0}=\left(-i\nabla+\left(\frac{\alpha}{|x|^2}+\frac{B_0}{2}\right)(-x_2,x_1)\right)^2,\quad x=(x_1,x_2)\in\mathbb{R}^2\setminus\{0\}
\end{equation}
with $\alpha\in\mathbb{R}$ being the total flux of the Aharonov-Bohm solenoid around the unit circle $\mathbb{S}^1$ and $B_0>0$ being the strength of the up-oriented uniform magnetic field. Since the case of integer-valued total flux (i.e. $\alpha\in\mathbb{Z}$) can be reduced to the vanishing case (i.e. $\alpha=0$) via a suitable gauge transformation (see e.g. \cite{ESV02}), it is a convention to consider the reduced range $\alpha\in(0,1)$ instead of $\alpha\in\mathbb{R}$ without loss of generality. We stress that the Hamiltonian \eqref{H-A} is a two-dimensional model and since the magnetic field $B(x)$ generated by a general vector-valued potential $A(x)=(A_1(x),A_2(x))$ in dimension two is given by
\begin{equation}\label{B-n}
B(x):=\frac{\partial A_2}{\partial x_1}-\frac{\partial A_1}{\partial x_2},
\end{equation}
the resulting magnetic field of the model \eqref{H-A} is actually a superposition of the uniform magnetic field and the singular Aharonov-Bohm magnetic field
\begin{equation*}
B(x)=B_0+\alpha \delta(x),\quad B_0>0,\quad \alpha\in(0,1),
\end{equation*}
where $\delta$ denotes Dirac's delta distribution. Due to the presence of the growing potential $\frac{B_0}{2}(-x_2,x_1)$, the spectrum of the Schr\"odinger operator \eqref{H-A} consists of purely discrete eigenvalues so that the dispersive behaviors of the Schr\"odinger and wave equations related to $H_{\alpha,B_0}$ are different from the scaling-critical models in \cite{FZZ22,GYZZ22}.

\subsection{Motivations}

If we use the polar coordinates $(r,\theta)$ of the Euclidean plane $\mathbb{R}^2$, the operator \eqref{H-A} has the following expression
\begin{equation}\label{polar}
H_{\alpha,B_0}=-\partial_{rr}-\frac{1}{r}\partial_r+\frac{1}{r^2}\left(-i\nabla_{\mathbb{S}^1}+\alpha+\frac{B_0}{2}r^2\right)^2,\quad \alpha\in(0,1),B_0>0,
\end{equation}
where $\mathbb{S}^1$ denotes the unit circle in the Euclidean plane $\mathbb{R}^2$. A natural way to generalize the operator \eqref{polar} is to consider it on the product cone $X:=X(\mathbb{S}_\sigma^1)=(0,+\infty)_r\times\mathbb{S}_\sigma^1$ equipped with the flat metric $g=dr^2+r^2d\theta^2$
\begin{equation}\label{operator}
H_{\alpha,B_0,\sigma}:=-\partial_{rr}-\frac{1}{r}\partial_r+\frac{1}{r^2}\left(-i\nabla_{\mathbb{S}_\sigma^1}+\alpha+\frac{B_0}{2}r^2\right)^2,\quad \alpha\in(0,1/\sigma),B_0>0,
\end{equation}
where $\mathbb{S}_\sigma^1$ stands for the circle of radius $\sigma\geq1$ in the Euclidean plane $\mathbb{R}^2$. Here, the requirement $\alpha\in(0,1/\sigma)\subseteq(0,1)$ comes from the reduction of the total flux $\alpha\in\mathbb{R}\setminus\sigma^{-1}\mathbb{Z}$ of the Aharonov-Bohm solenoid around the dilated circle $\mathbb{S}_\sigma^1$.
In this paper, we aim to explore the decay properties for solutions to the Schr\"odinger equation \eqref{eq:S} and the wave equation \eqref{eq:W} associated to the operator $H_{\alpha,B_0,\sigma}$ on the product cone $X=(0,+\infty)_r\times\mathbb{S}_\sigma^1$. To avoid confusion, we will sometimes write $(r,\theta)\in(0,+\infty)_r\times\mathbb{S}_\sigma^1$ instead of $x\in X$. As the self-adjoint extension theory of the Euclidean operator $H_{\alpha,B_0}$ on $L^2(\mathbb{R}^2)$ (see e.g. \cite{ESV02}), Friedrichs extension of the conic operator \eqref{operator} yields a self-adjoint operator on $L^2(X)$, which we still denote by $H_{\alpha,B_0,\sigma}$ throughout the whole paper. For the self-adjoint operator $H_{\alpha,B_0,\sigma}$, time-evolutionary propagators like $e^{itH_{\alpha,B_0,\sigma}},e^{-tH_{\alpha,B_0,\sigma}},e^{it\sqrt{H_{\alpha,B_0,\sigma}}}$ will be meaningful and thus we are allowed to study a large class of dispersive estimates, such as time-decay (perhaps local in time), Strichartz and local smoothing, for dispersive equations.

The general form of the Schr\"odinger operator with electromagnetic potentials on the product cone $C(Y):=(0,+\infty)_r\times Y$ is given by
\begin{equation}\label{op:general}
\mathcal{L}_{A,V}=\frac{1}{\sqrt{|g|}}\sum_{j,k=1}^n\left(i\frac{\partial}{\partial x_j}+A_j\right)\sqrt{|g|}g^{jk}\left(i\frac{\partial}{\partial x_k}+A_k\right)+V(x), x\in C(Y),
\end{equation}
where $V:C(Y)\rightarrow\mathbb{R}$ is a real-valued electric potential and the 1-form $A=A_jdx^j$ represents the vector-valued magnetic potential. Here $C(Y)$ stands for the general product cone endowed with the metric $g=dr^2+r^2h$, where $(Y,h)$ denotes a compact Riemannian manifold representing the cross-section of the cone $C(Y)$. It should be pointed out that the general product cone $C(Y)$ preserves the dilation symmetry of the usual Euclidean spaces (i.e. the special case $Y=\mathbb{S}^{n-1}$), but no other symmetries in general. In the Euclidean case $Y=\mathbb{S}^{n-1}$, decay estimates for the Schr\"odinger and wave equations have been extensively studied with various kinds of asymptotic decaying and regularity conditions imposed on the potential $V$. In their seminal work, Journ\'{e}-Soffer-Sogge \cite{JSS91} obtained the dispersive estimate for the Schr\"odinger operator with strongly decaying and regular potentials in $\mathbb{R}^n(n\geq3)$. Since then, dispersive estimates have been explored in many different situations, among which Goldberg, together with his collaborators \cite{BG12,Gol06,Gold06,GS04,GV06} systematically explored the regularity conditions required for the potential in higher dimensions and the decay rates in dimension three in the absence of embedded resonances and eigenvalues. Subsequent results on Schr\"odinger operators with potentials in $\mathbb{R}^n$ have extended dispersive estimate to the setting of $-\Delta+V$ admitting resonances at zero energy, thereby resulting in a slower decay rate of the dispersive estimate (see e.g. Erdo\u{g}an-Schlag \cite{ES06} for $n=3$, Erdo\u{g}an-Green \cite{EG13,EG13-CMP} for $n=2$, Green \cite{Green12} for $n=5$, as well as Goldberg-Green \cite{GG15,GG16,GG17} for $n\geq4$). For the critically decaying electric potential $V$ (including the inverse-square potential), Strichartz estimates for the Schr\"odinger and wave equations were established in \cite{BPST03,BPST04,PST03}. Recently, Blair-Sire-Sogge \cite{BSS21} pushed the construction of the spectral measure for $-\Delta_g+V$ on compact manifolds to the situations where the regularity of the potential $V$ is at very critical levels, although they did not establish the dispersive estimate directly. This is of course not an exhaustive list of references in the purely electric potential case, but these results are representative for the techniques involved in the Euclidean setting (based on Stone's formula). The survey paper of Schlag \cite{Sch07} provides an excellent overview of the key ideas involved. Another well-known method to obtain the dispersive estimate for the Schr\"odinger equation in the purely electric potential case is the so-called wave operator method developed by Yajima \cite{Yaj4,Yaj5,Yaj6,Yaj7}.

In the presence of a magnetic potential, the picture becomes more subtle due to the failure of the above mentioned methods in general. For some earlier results, we refer to e.g. \cite{CS01,DP05,DF07,DF08,DFVV10,EGS08,EGS09,GST07}. The critical magnetic potentials occupy a special position, for example, Fanelli-Felli-Fontelos-Primo \cite{FFFP13,FFFP15} obtained a general representation formula for solutions to the Schr\"odinger equation via the pseudo-conformal transformation and then proved the dispersive estimates for the most general electromagnetic potential case in dimension two. As for the wave equation and Klein-Gordon equation with critical magnetic potentials, dispersive estimates were established recently by Fanelli-Zhang-Zheng \cite{FZZ22} via Hankel transform and by Gao-Yin-Zhang-Zheng \cite{GYZZ22} via the spectral measure argument, respectively. A typical example of the critical magnetic potential is given by the two-dimensional model $\alpha|x|^{-2}(-x_2,x_1)$, which was originally introduced by Aharonov-Bohm \cite{AB59} to model a well-known physical phenomenon (i.e. Aharonov-Bohm effect) describing the dynamics of a charged particle when passing by an infinitely thin solenoid (i.e. the Aharonov-Bohm solenoid). The Schr\"odinger operator with such a magnetic potential is often referred to as the Aharonov-Bohm Hamiltonian in the physical literature. One of the most famous growing potential model in quantum mechanics is given by the uniform magnetic potential $\frac{B_0}{2}(-x_2,x_1)$ and Schr\"odinger operator with this potential is often referred to as the Landau Hamiltonian, which was originally introduced by Landau \cite{Lan30} and systematically studied by Avron-Herbst-Simon \cite{AHS1,AHS2,AHS3} from the perspective of mathematical physics. Moreover, it is not hard to realize that the Landau Hamiltonian actually coincides with the one-dimensional twisted Laplacian (also called the special Hermite operator, see e.g. \cite{DPR10}). The operator $H_{\alpha,B_0,\sigma}$ concerned in the present paper actually models the superposition effect of one Aharonov-Bohm solenoid in the uniform magnetic field (formally, the coupling of the Aharonov-Bohm Hamiltonian and the Landau Hamiltonian).

The consideration of conic manifolds arises naturally in the recent work of Hintz-Vasy \cite{HV15,HV18} and Hafner-Hintz-Vasy \cite{HHV20} on the black holes in general relativity since the cosmic space is not always flat (see also \cite{Hin22}). When the potentials $A,V$ in \eqref{op:general} vanish, the operator $\mathcal{L}_{A,V}$ is reduced to the Laplacian defined on product cones. In the seminal papers of Cheeger-Taylor \cite{CT82-I,CT82-II}, the spectral measure for the Laplacian on product cones was firstly constructed by using the Lipschitz-Hankel integral formula and the Strichartz estimates for the corresponding wave equation were established as an application. Since then, there have been many results concerning evolutionary equations on product cones, especially the decay estimates for dispersive equations (see e.g. \cite{BFHM12,BFM13,For10,FHH18,MW04,ZZ20,Zhang22}). In the case of critically decaying potential $V(x)=\frac{a(\theta)}{|x|^2}$, the Schr\"odinger operator $\mathcal{L}_{0,V}$ has attracted several attentions from other topics. For instance, asymptotical behavior of the Schr\"odinger propagator was studied in \cite{Wang06}, Riesz transform was investigated in \cite{HL14,Li99} and the Strichartz and restriction estimates were obtained in \cite{Zhang15,ZZ20}. Dispersive estimate for the Schr\"odinger equation on product cones has been established recently by Zhang-Zheng \cite{ZZ17-CPDE,ZZ20} via the spectral measure technique, which is more matchable for the cases of decaying potential $V$ (see also \cite{KM22}). The Schwartz kernel of the spectral measure related to Schr\"odinger operators in conically singular spaces has been systematically studied by Hassell-Vasy \cite{HV01} and Guillarmou-Hassell-Sikora \cite{GHS13}. It has a lot of applications to the harmonic analysis problems (especially to those under the $L^p$-frame) in conically singular spaces. For more examples, we refer to Hassell-Zhang \cite{HZ16} for the Strichartz estimates of Schr\"odinger equation, of which the arguments heavily rely on the microlocal techniques developed by Melrose \cite{Mel94} due to the generality and complexity of the geometry therein, and to Zhang \cite{Zhang22} for the construction of Schwartz kernels of the resolvent and spectral measure associated to the Laplacian $\Delta_g$ on the special product cone $X:=\mathcal{C}(\mathbb{S}_\sigma^1)=(0,+\infty)_r\times\mathbb{S}^1_\sigma$ equipped with the metric $g=dr^2+r^2d\theta^2$, where $\mathbb{S}_\sigma^1\simeq\mathbb{R}/2\pi\sigma\mathbb{Z}$ stands for the planar circle of radius $\sigma>0$ in the Euclidean space $\mathbb{R}^2$. Dispersive estimates for the Schr\"odinger and wave equations have also been explored in several other geometries, for example, the hyperbolic space by Borthwick-Marzuola \cite{BM15}, asymptotically Euclidean space by Bouclet-Burq \cite{BB21}, asymptotically conic manifolds by Hassell-Zhang \cite{HZ16} and manifolds with conic ends by Schlag-Soffer-Staubach \cite{SSS10-I,SSS10-II}. The up-to-date progress on decay and Strichartz estimates for higher dimensional Schr\"odinger and wave equations on conically singular spaces was made by Jia-Zhang \cite{JZ1,JZ2} and Fanelli-Su-Wang-Zhang-Zheng \cite{FSWZZ24}.

\subsection{Main results}


Now we state the first result concerning the weighted dispersive inequality for the Schr\"odinger equation \eqref{eq:S} associated to the operator \eqref{operator} on the product cone $(X,g)$.
\begin{theorem}[Weighted dispersive estimate for Schr\"odinger]\label{thm:S}
Let $X:=\mathcal{C}(\mathbb{S}_\sigma^1)=(0,+\infty)_r\times\mathbb{S}^1_\sigma$ be the product cone endowed with a flat metric $g=dr^2+r^2d\theta^2$ and $H_{\alpha,B_0,\sigma}$ the Friedrichs extension of the operator \eqref{operator} with $\alpha\in\mathbb{R}/\sigma^{-1}\mathbb{Z}$ and $B_0>0,\sigma\geq1$. If we set $\kappa_\sigma=\text{dist}(\alpha,\sigma^{-1}\mathbb{Z})$ for a fixed $\sigma\geq1$, then there exists a constant $C=C_{B_0,\sigma}>0$ such that, for all $\gamma\in[0,\kappa_\sigma]$ and $t\neq\frac{\pi}{B_0}\mathbb{Z}$, it holds the weighted dispersive estimate
\begin{equation}\label{decay:S}
\||\cdot|^{-\gamma}e^{itH_{\alpha,B_0,\sigma}}|\cdot|^{-\gamma}\|_{L^1(X)\rightarrow L^\infty(X)}\leq C|\sin(tB_0)|^{-1-\gamma}.
\end{equation}
\end{theorem}
\begin{remark}
As mentioned in \eqref{operator}, we shall consider the reduction range $\alpha\in(0,1/\sigma)$ throughout the whole paper. It is not difficult to see that the maximal possible value of $\kappa_\sigma$ relies on the radius $\sigma$ of the sectional circle $\mathbb{S}_\sigma^1$ of the product cone $X$ and is given by $\frac{1}{2\sigma}$. Correspondingly, we have
\begin{equation*}
\||\cdot|^{-\frac{1}{2\sigma}}e^{itH_{\alpha,B_0,\sigma}}|\cdot|^{-\frac{1}{2\sigma}}\|_{L^1(X)\rightarrow L^\infty(X)}\leq C|\sin(tB_0)|^{-1-\frac{1}{2\sigma}},\quad \alpha\in(2\sigma)^{-1}\mathbb{Z}.
\end{equation*}
In the case of vanishing weight (i.e. $\gamma=0$), \eqref{decay:S} is reduced to the usual dispersive estimate for the Schr\"odinger equation \eqref{eq:S}
\begin{equation}\label{dis:S}
\|e^{itH_{\alpha,B_0,\sigma}}f\|_{L^\infty(X)}\lesssim_\sigma\frac{B_0}{|\sin(tB_0)|}\|f\|_{L^1(X)},\quad t\neq\frac{\pi}{B_0}\mathbb{Z},
\end{equation}
which further covers the Euclidean result of \cite[Theorem 1.1]{WZZ25-AHP} (i.e. the case $\sigma=1$). The optimal decay rate in \eqref{decay:S} depends on the ground level $\kappa_\sigma^2$ of the angular part of the Aharonov-Bohm Hamiltonian or the total flux $\alpha$ of the Aharonov-Bohm solenoid around the circle $\mathbb{S}_\sigma^1$. For this reason, the weighted estimate \eqref{decay:S} is essentially a corollary of the dispersive estimate \eqref{dis:S}.
The appearance of the sine factor $\sin(tB_0)$ is natural due to the confinement effect of the uniform magnetic field $B_0$. This phenomenon is closely related to the Euclidean harmonic oscillator $H=-\Delta+|x|^2$, for which Koch and Tataru proved in \cite{KT05}
\begin{equation*}
  \|e^{itH}\|_{L^1(\mathbb{R}^n)\rightarrow L^\infty(\mathbb{R}^n)}\lesssim_n|\sin(2t)|^{-\frac{n}{2}}.
\end{equation*}
\end{remark}
\begin{remark}
The proof of the Schr\"odinger dispersive estimate \eqref{dis:S} reduces to the establishment of a uniform bound (by $\frac{B_0}{|\sin(tB_0)|}$) for the kernel of the Schr\"odinger propagator $e^{itH_{\alpha,B_0,\sigma}}$ and we would like to achieve this goal by three different approaches, which are all based on a series representation of the Schr\"odinger kernel obtained by the spectral theorem. Due to the appearance of the uniform magnetic field $B_0$, the spectrum of the Schr\"odinger operator $H_{\alpha,B_0,\sigma}$ consists of purely discrete eigenvalues and it is thus possible to obtain a series representation for the corresponding Schr\"odinger kernel by the spectral theorem. The first approach is to directly construct a closed expression for the Schr\"odinger kernel by the Poisson summation formula from the series expression of the Schr\"odinger kernel and the second approach relies on applying the universal covering space technique, which was originally developed by Schulman \cite{Sch68} in connection with the Feynman path integral on multiply connected spaces, to obtain the same closed expression for the Schr\"odinger kernel, thereby verifying a connection formula for the Schr\"odinger evolutions on two different product spaces. The expected Schr\"odinger dispersive estimate \eqref{dis:S} then follows by bounding the explicit representation formula of the Schr\"odinger kernel via elementary calculus. The third approach is more straightforward and based on a uniform boundedness result summarized from the recent preprint \cite[Section 2]{GWXZ25}
\begin{equation*}
\sup_{r\geq0\atop |\delta|<R}\left|\sum_{k\in\mathbb{Z}}e^{i\frac{k}{\sigma}\delta}I_{|k/\sigma+\alpha|}(ir)\right|<+\infty,\quad R<+\infty.
\end{equation*}
\end{remark}
Before stating the second result, let us introduce some preliminary notations.
Let $\varphi\in C_c^\infty(\mathbb{R}_+)$ satisfy $0\leq\varphi\leq1,\text{supp}\varphi\subset[1/2,1]$ and
\begin{equation}\label{LP-dp}
\sum_{j\in\mathbb{Z}}\varphi_j(\lambda)=1,\quad \forall\lambda>0,\qquad \varphi_j(\lambda):=\varphi(2^{-j}\lambda),\quad j\in\mathbb{Z}.
\end{equation}
The Besov space related to the operator $H_{\alpha,B_0,\sigma}$ is defined as follows:
\begin{definition}[Besov spaces related to $H_{\alpha,B_0,\sigma}$]\label{def:besov}
For $s\in\mathbb{R}$ and $1\leq p,q<\infty$, the homogeneous Besov norm $\|\cdot\|_{\dot{\mathcal{B}}^s_{p,q,\sigma}(X)}$ is defined by
\begin{align}
\|f\|_{\dot{\mathcal{B}}^s_{p,q,\sigma}(X)}:&=\left(\sum_{j\in\mathbb{Z}}2^{jsq}\|\varphi_j(\sqrt{H_{\alpha,B_0,\sigma}})f\|_{L^p(X)}^q\right)^{1/q}\label{Besov}\\
&=\left(\sum_{j\in\mathbb{Z}}2^{jsq}\left(\int_0^\infty\int_{\mathbb{S}_\sigma^1}
\left|\left(\varphi_j(\sqrt{H_{\alpha,B_0,\sigma}})f\right)(r,\theta)\right|^p\mathrm{d}S(\theta)r\mathrm{d}r\right)^{q/p}\right)^{1/q}.\nonumber
\end{align}
In particular, for $p=q=2$, we have the Sobolev norm
\begin{equation}\label{Sobolev1}
\|f\|_{\dot{\mathcal{H}}^s_{\alpha,B_0,\sigma}(X)}:=\|f\|_{\dot{\mathcal{B}}^s_{2,2,\sigma}(X)}.
\end{equation}
\end{definition}
\begin{remark}
The Sobolev space can be defined alternatively by
\begin{equation*}
\dot{\mathcal{H}}^s_{\alpha,B_0,\sigma}(X):=H_{\alpha,B_0,\sigma}^{-s/2}L^2(X)
\end{equation*}
with the following norm
\begin{equation}\label{Sobolev2}
\begin{split}
\|f\|_{\dot{\mathcal{H}}^s_{\alpha,B_0,\sigma}(X)}:&=\left\|H_{\alpha,B_0,\sigma}^{s/2}f\right\|_{L^2(X)}\\
&=\left(\int_0^\infty\int_{\mathbb{S}_\sigma^1}\left|\left(\varphi_j(\sqrt{H_{\alpha,B_0,\sigma}})f\right)(r,\theta)\right|^2\mathrm{d}S(\theta)r\mathrm{d}r\right)^{1/2}.
\end{split}
\end{equation}
In view of the spectral theory for operators on $L^2$, the norms in \eqref{Sobolev1} and \eqref{Sobolev2} are equivalent (see Proposition \ref{prop:E-S-B} below for a proof).
\end{remark}
Now we state the second result concerning the decay estimate for the wave equation \eqref{eq:W} associated to the operator \eqref{operator} on the product cone $(X,g)$.
\begin{theorem}[Decay estimate for wave]\label{thm:W}
Let $X:=\mathcal{C}(\mathbb{S}_\sigma^1)=(0,+\infty)_r\times\mathbb{S}^1_\sigma$ be the product cone endowed with a flat metric $g=dr^2+r^2d\theta^2$ and $H_{\alpha,B_0,\sigma}$ the Friedrichs extension of the operator \eqref{operator} with $\alpha\in(0,1/\sigma)$ and $B_0>0,\sigma\geq1$, then there exists a constant $C=C_{B_0,\sigma}>0$ such that, for any finite $T>0$, it holds the dispersive estimate for the wave equation \eqref{eq:W}
\begin{equation}\label{dis:W}
\left\|\frac{\sin(t\sqrt{H_{\alpha,B_0,\sigma}})}{\sqrt{H_{\alpha,B_0,\sigma}}}f\right\|_{L^\infty(X)}\leq C|t|^{-1/2}\|f\|_{\dot{\mathcal{B}}^{1/2}_{1,1,\sigma}(X)},\quad\forall t\in I:=[0,T],
\end{equation}
where $\|\cdot\|_{\dot{\mathcal{B}}^s_{p,q,\sigma}(X)}$ denotes the homogeneous Besov norm in \eqref{Besov}.
\end{theorem}
\begin{remark}
The result \eqref{dis:W} is a naturally generalization of the Euclidean counterpart (corresponding to $\sigma=1$) in \cite[Theorem 1.4]{WZZ25-JDE} provided that the radius of the sectional circle $\mathbb{S}_\sigma^1$ of the product cone $X$ is fixed. In the Euclidean plane $\mathbb{R}^2$ (i.e. $\sigma=1$), if the Aharonov-Bohm solenoid is further removed (i.e. $\alpha=0$), the inequality \eqref{dis:W} will be reduced to the dispersive estimate for the wave equation associated to the Landau Hamiltonian $\mathcal{L}:=(-i\nabla+\frac{B_0}{2}(-x_2,x_1))^2$ . In fact, dispersive estimate for the wave equation associated to the Landau Hamiltonian $\mathcal{L}$ has been obtained by D'Ancona-Pierfelice-Ricci \cite[Theorem 1.1]{DPR10}
\begin{equation*}
\|e^{it\sqrt{H_{\alpha,B_0,\sigma}}}f\|_{L^\infty(\mathbb{R}^2)}\lesssim|t|^{-1/2}\|f\|_{\dot{\mathcal{B}}^{3/2}_{1,1}(\mathbb{R}^2)},\quad\forall |t|\leq T<\frac{\pi}{2B_0},
\end{equation*}
compared to which, we are able to relax their restriction $T<\frac{\pi}{2B_0}$ to any finite $T<+\infty$ here by following their arguments with some modifications.
\end{remark}
\begin{remark}
By the definition of the homogeneous Besov norm \eqref{Besov}, the desired estimate \eqref{dis:W} is reduced to the following $L^1\rightarrow L^\infty$ estimates
\begin{equation}\label{dis:W-1}
\begin{split}
&\left\|\varphi(2^{-j}\sqrt{H_{\alpha,B_0,\sigma}})e^{it\sqrt{H_{\alpha,B_0,\sigma}}}f\right\|_{L^\infty(X)}\\
&\qquad\qquad \lesssim2^{2j}\big(1+2^jt\big)^{-N}\|\varphi(2^{-j}\sqrt{H_{\alpha,B_0,\sigma}})f\|_{L^1(X)} \quad\text{for}\quad 2^jt\lesssim1
\end{split}
\end{equation}
and
\begin{align}\label{dis:W-2}
&\left\|\varphi(2^{-j}\sqrt{H_{\alpha,B_0,\sigma}})e^{it\sqrt{H_{\alpha,B_0,\sigma}}}f\right\|_{L^\infty(X)}\\
&\lesssim2^{2j}\big(1+2^jt\big)^{-1/2}\|\varphi(2^{-j}\sqrt{H_{\alpha,B_0,\sigma}})f\|_{L^1(X)} \quad\text{for}\quad 2^jt\lesssim1, 2^{-j}|t|\leq\frac{\pi}{2B_0}.\nonumber
\end{align}
The possibility of obtaining an explicit representation formula for the Schr\"odinger kernel is mainly due to the separation feature of the explicit eigenvalues of the operator $H_{\alpha,B_0,\sigma}$. However, this significant feature fails in obtaining a proper representation for the kernel of the half-wave propagator $e^{it\sqrt{H_{\alpha,B_0,\sigma}}}$ due to the square root of the eigenvalue, which prevents us from directly proving the estimates \eqref{dis:W-1} and \eqref{dis:W-2}. In addition, the appearance of the uniform magnetic field $B_0$ of the operator $H_{\alpha,B_0,\sigma}$ results in the failure of the spectral measure argument via Stone's formula in \cite{GYZZ22,Zhang22}. To avoid constructing the spectral measure, we turn to prove the associated Bernstein inequality to address the low frequency. As for the high frequency, thanks to the classical subordination formula
\begin{equation}\label{sub:form}
e^{-z\sqrt{y}}=\frac{z}{2\sqrt{\pi}}\int_0^\infty e^{-sy-\frac{z^2}{4s}}s^{-\frac{3}{2}}\mathrm{d}s,\quad z>0,
\end{equation}
we are allowed to obtain the expected decay estimate for the half-wave propagator by using the dispersive estimate \eqref{dis:S} for the Schr\"odinger propagator. In order to establish the Littlewood-Paley theory (including the Bernstein inequalities and the square function estimates) associated to the operator $H_{\alpha,B_0,\sigma}$, a standard mechanism is to obtain a Gaussian-type upper bound for the related heat kernel. Gaussian boundedness of heat kernels associated to operators on various manifolds is of independent interest. Our strategy is again to apply Schulman's abstract ansatz (i.e. the universal covering space technique) so that it is possible to obtain an explicit representation formula for the heat kernel associated to the operator $H_{\alpha,B_0,\sigma}$ and then the desired Gaussian-type bound follows by suitably analysing this explicit heat kernel. Hence, the main efforts to obtain \eqref{dis:W} will be devoted to prove the Gaussian heat kernel estimate.
\end{remark}
The idea of deducing a dispersive estimate for the half-wave flow $e^{it\sqrt{H_{\alpha,B_0,\sigma}}}$ from the dispersive estimate for the corresponding Schr\"odinger flow $e^{itH_{\alpha,B_0,\sigma}}$ was originally suggested by M\"uller-Seeger \cite{MS15} (via the above subordination formula \eqref{sub:form} bridging such two flows) and applied later by D'Ancona-Pierfelice-Ricci \cite{DPR10} to deal with the Hermite operator and the twisted Laplacian and recently by Bui-D'Ancona-Duong-M\"uller \cite{BDDM19} to address general nonnegative self-adjoint operators $L$ with power-like behaviors near $0$ and $\infty$ on metric measure spaces provided that both the Schr\"odinger dispersive estimate of the form $\|e^{itL}\|_{L^1\rightarrow L^\infty}\lesssim|t|^{-n/2}$ and the Gaussian heat kernel estimate of the type $|e^{-tL}(x,y)|\lesssim\frac{1}{\mu(B(x,\sqrt{t}))}e^{-\frac{d(x,y)^2}{ct}}$ are valid.

Therefore, in view of \cite[Theorem 1.1]{BDDM19}, we in fact obtain more general results for operators $H_{\alpha,B_0,\sigma}^\nu$ of fractional order $\nu\in(0,1)$.
\begin{theorem}\label{thm:F}
Assume $\nu\in(0,1)$ and $s>2-\nu$. Then, under the conditions of Theorem \ref{thm:S}, we have
\begin{equation}\label{dis:Sobolev}
\|e^{itH_{\alpha,B_0,\sigma}^\nu}f\|_{L^\infty(X)}\lesssim|t|^{-1/2}\|(I+H_{\alpha,B_0,\sigma})^{s/2}f\|_{L^1(X)},\quad\forall t\in I:=[0,T]
\end{equation}
and
\begin{equation}\label{dis:Besov}
\|e^{itH_{\alpha,B_0,\sigma}^\nu}f\|_{L^\infty(X)}\lesssim|t|^{-1/2}\|f\|_{\dot{\mathcal{B}}^{2-\nu}_{1,1,\sigma}(X)},\quad\forall t\in I:=[0,T],
\end{equation}
where $\|\cdot\|_{\dot{\mathcal{B}}^s_{p,q,\sigma}(X)}$ denotes the homogeneous Besov norm in \eqref{Besov}.
\end{theorem}
\begin{remark}
We shall not present the details of the proof of Theorem \ref{thm:F} since one can completely mimic the steps of \cite{BDDM19} for abstract nonnegative self-adjoint operators on metric measure spaces once the Schr\"odinger dispersive estimate in Theorem \ref{thm:S} and the Gaussian heat kernel estimate were established.
\end{remark}
One of the immediate applications of the dispersive estimate is to derive the corresponding Strichartz estimates via Keel-Tao \cite{KT98}. For instance, it is well-known that the dispersive estimate for the free Schr\"odinger equation
\begin{equation*}
\|e^{it\Delta}f\|_{L^\infty(\mathbb{R}^2)}\lesssim|t|^{-1}\|f\|_{L^1(\mathbb{R}^2)},
\end{equation*}
together with the conservation of mass
\begin{equation*}
\|e^{it\Delta}f\|_{L^2(\mathbb{R}^2)}=\|f\|_{L^2(\mathbb{R}^2)},
\end{equation*}
yields, via the standard $TT^*$ argument of Keel-Tao \cite{KT98}, the following Strichartz estimates
\begin{equation}\label{Str:free}
\|e^{it\Delta}f(x)\|_{L^p_t(\mathbb{R};L_x^q(\mathbb{R}^2))}\lesssim\|f\|_{L^2(\mathbb{R}^2)},
\end{equation}
where the pair $(p,q)$ fulfills the scaling condition
\begin{equation}\label{adm:pair}
\frac{1}{p}=\frac{1}{2}-\frac{1}{q},\quad(p,q)\in[2,+\infty]\times[2,+\infty).
\end{equation}
For convenience, we shall denote by $(q,p)\in\Lambda^S$ if $(q,p)$ satisfies the condition \eqref{adm:pair} in what follows.
Strichartz estimates for the free wave equation are given by
\begin{equation}\label{Str:free-wave}
\left\|\cos(t\sqrt{-\Delta})u_0+\frac{\sin(t\sqrt{-\Delta})}{\sqrt{-\Delta}}u_1\right\|_{L_t^q(\mathbb{R};L_x^p(\mathbb{R}^2))}
\lesssim\|u_0\|_{\dot{H}^s(\mathbb{R}^2)}+\|u_1\|_{\dot{H}^{s-1}(\mathbb{R}^2)},
\end{equation}
where $\|\cdot\|_{\dot{H}^s(\mathbb{R}^2)}$ is the classical homogeneous Sobolev norm associated to the Laplacian $-\Delta$ and the pair $(q,p)$ satisfies the wave-admissible condition
\begin{equation}\label{adm:wave}
\frac{1}{q}\leq\frac{1}{2}\left(\frac{1}{2}-\frac{1}{p}\right),\quad (q,p)\in[2,\infty]\times[2,\infty).
\end{equation}
Throughout the paper, we shall write for simplicity $(q,p)\in\Lambda_s^W$ if $(q,p)$ satisfies \eqref{adm:wave} and
\begin{equation}\label{gap}
\frac{1}{q}+\frac{2}{p}=1-s,\quad s\in[0,1).
\end{equation}
Hence, as a direct consequence of results in Theorem \ref{thm:S} and Theorem \ref{thm:W}, we can obtain similar Strichartz estimates for the Schr\"odinger and wave equations associated to the operator $H_{\alpha,B_0,\sigma}$ on the product cone $(X,g)$.
\begin{corollary}[Strichartz estimates]\label{cor:S-W}
Under the conditions of Theorem \ref{thm:S}, there hold for $(q,p)\in\Lambda^S$ the Strichartz estimates of the Schr\"odinger equation
\begin{equation}\label{str:S}
\|e^{itH_{\alpha,B_0,\sigma}}f\|_{L_t^q((0,\pi/B_0);L_x^p(X))}\lesssim\|f\|_{L^2(X)}
\end{equation}
and for $(q,p)\in\Lambda_s^W$ the Strichartz estimates of the wave equation
\begin{equation}\label{str:W}
\begin{split}
&\left\|\cos(t\sqrt{H_{\alpha,B_0,\sigma}})u_0+\frac{\sin(t\sqrt{H_{\alpha,B_0,\sigma}})}{\sqrt{H_{\alpha,B_0,\sigma}}}u_1\right\|_{L_t^q(I;L_x^p(X))} \\
&\qquad\qquad \lesssim\|u_0\|_{\dot{\mathcal{H}}^s_{\alpha,B_0,\sigma}(X)}+\|u_1\|_{\dot{\mathcal{H}}^{s-1}_{\alpha,B_0,\sigma}(X)},\quad I=[0,T].
\end{split}
\end{equation}
\end{corollary}
\begin{remark}
The special case $\sigma=1,\alpha=0$ of the Schr\"odinger Strichartz estimate \eqref{str:S} corresponds to \cite[Theorem 2]{KL21} for the Landau Hamiltonian and the Euclidean case $\sigma=1$ of \eqref{str:S} corresponds to the recent result of \cite[(1.9)]{WZZ25-AHP}. In the special case $\sigma=1,\alpha=0$, \eqref{str:W} is reduced to Strichartz estimate for the wave equation associated to the Landau Hamiltonian and the Euclidean case $\sigma=1$ of the wave Strichartz estimate \eqref{str:W} corresponds to the recent result of \cite[(1.13)]{WZZ25-JDE}.
\end{remark}
This paper is organized as follows. In Section \ref{sec:pre}, as a foundation of the whole paper, we will describe the explicit spectral properties of the Friedrichs extension of the operator \eqref{operator} and then give a proof of the equivalence between the associated Sobolev norm and a special Besov norm by the spectral theorem. In Section \ref{sec:const}, we firstly present the steps to obtain the weighted dispersive estimate \eqref{decay:S} for the Schr\"odinger equation via the weightless dispersive estimate \eqref{dis:S} and then the efforts will be devoted to obtain the dispersive estimate \eqref{dis:S} via three different approaches. In Section \ref{sec:heat}, we mainly obtain an explicit representation formula for the related heat kernel by the universal covering space technique and then establish its Gaussian boundedness result, thereby yielding the corresponding Bernstein inequalities and the square function estimates via standard arguments.  Finally, in Section \ref{sec:decay}, we prove the dispersive estimate \eqref{dis:W} for the wave equation with the aid of M\"uller-Seeger's subordination formula and also sketch the proof of the related Strichartz estimate \eqref{str:W}.
\vspace{0.2cm}

{\bf Acknowledgments:}\quad  This work is supported by High-level Talent Research Start-up Project Funding of Henan Academy of Sciences (No. 20251819008).
\vspace{0.2cm}

\section{preliminaries} \label{sec:pre}

In this section, we describe the spectral properties of the Friedrichs extension of the concerned operator $H_{\alpha,B_0,\sigma}$. Next, we apply the spectral theory to show the equivalence between the Sobolev norm and a special Besov norm.

\subsection{Quadratic form and Friedrichs extension}

Let us recall the operator $H_{\alpha,B_0,\sigma}$ on the product cone $X=(0,+\infty)_r\times\mathbb{S}_\sigma^1$ with the metric $g=dr^2+r^2d\theta^2$
\begin{equation}\label{op}
H_{\alpha,B_0,\sigma}=-\partial_{rr}-\frac{1}{r}\partial_r+\frac{1}{r^2}\left(-i\nabla_{\mathbb{S}_\sigma^1}+\alpha+\frac{B_0r^2}{2}\right)^2,\quad \alpha\in(0,1/\sigma),B_0>0,
\end{equation}
where $\mathbb{S}_\sigma^1$ denotes the circle of radius $\sigma\geq1$ in the Euclidean plane $\mathbb{R}^2$.
Associated to the operator \eqref{op}, the Sobolev space $\mathcal{H}_{\alpha,B_0,\sigma}^1(X)$ is defined as the completion of the set $C_c^\infty(\bar{X}\setminus P;\mathbb{C})$ with respect to the norm
\begin{equation*}
\|f\|_{\mathcal{H}_{\alpha,B_0,\sigma}^1(X)}=\left(\int_0^\infty\int_0^{2\sigma\pi}H_{\alpha,B_0,\sigma}f(r,\theta)\overline{f(r,\theta)}\mathrm{d}\theta r\mathrm{d}r\right)^{1/2}
\end{equation*}
and the quadratic form $Q_{\alpha,B_0,\sigma}$ is given by
\begin{equation*}
\begin{split}
Q_{\alpha,B_0,\sigma}: & \quad \quad \mathcal{H}_{\alpha,B_0,\sigma}^1\rightarrow\mathbb{R},\\
Q_{\alpha,B_0,\sigma}(f)&=\|f\|^2_{\mathcal{H}_{\alpha,B_0,\sigma}^1(X)},
\end{split}
\end{equation*}
where $P$ denotes the tip of the complete cone $\bar{X}$ and $\overline{f(r,\theta)}$ stands for the complex conjugate of the function $f(r,\theta)$.
Clearly, the quadratic form $Q_{\alpha,B_0,\sigma}$ is positive definite and the operator $H_{\alpha,B_0,\sigma}$ is symmetric. Thus, $H_{\alpha,B_0,\sigma}$ admits a self-adjoint extension (particularly, Friedrichs extension) $H^F_{\alpha,B_0,\sigma}$ on $L^2(X)$ with a natural form domain
\begin{equation}\label{domain}
\mathcal{D}:=\mathcal{D}(H_{\alpha,B_0,\sigma})=\left\{f\in\mathcal{H}_{\alpha,B_0,\sigma}^1(X):H_{\alpha,B_0,\sigma}f\in L^2(X)\right\}.
\end{equation}
We point out that the domain \eqref{domain} is strictly contained in the Sobolev space $H_{0,0,\sigma}^1(X)$ even in the Euclidean case $\sigma=1$. Although the operator $H_{\alpha,B_0,\sigma}$ admits many other self-adjoint extensions (see e.g. \cite{ESV02} for the Euclidean case $\sigma=1$) by the von Neumann self-adjoint extension theory, we are interested in the simplest Friedrichs extension and still denote by $H_{\alpha,B_0,\sigma}$ its Friedrichs extension $H^{F}_{\alpha,B_0,\sigma}$ in what follows.
Now we give the spectral properties of the Hamiltonian $H_{\alpha,B_0,\sigma}$.
\begin{proposition}\label{prop:spect}
Let $H_{\alpha,B_0,\sigma}$ denote the Friedrichs extension of the operator \eqref{op} with $\alpha\in\mathbb{R},B_0,\sigma>0$, then the spectrum of $H_{\alpha,B_0,\sigma}$ consists of purely discrete eigenvalues of finite multiplicity
\begin{equation}\label{eigen-v}
\lambda_{k,m}=\left(2m+1+\left|\frac{k}{\sigma}+\alpha\right|+\frac{k}{\sigma}+\alpha\right)B_0,\quad k\in\mathbb{Z},\quad m\in\mathbb{N}:=\{0,1,2,\cdots\}
\end{equation}
and the eigenfunction $V_{k,m}$ corresponding to each eigenvalue $\lambda_{k,m}$ is given by
\begin{equation}\label{eigen-f}
V_{k,m}(r,\theta)=\frac{1}{\sqrt{2\sigma\pi}}r^{|k/\sigma+\alpha|}e^{-\frac{B_0r^2}{4}}P_{k,m}\Bigg(\frac{B_0r^2}{2}\Bigg)e^{i\frac{k}{\sigma}\theta},
\quad\theta\in[0,2\sigma\pi],
\end{equation}
where $P_{k,m}$ is the polynomial of degree $m$
\begin{equation*}
P_{k,m}(r)=\sum_{n=0}^m\frac{(-m)_n}{(1+|k/\sigma+\alpha|)_n}\frac{r^n}{n!}
\end{equation*}
with $(a)_n(a\in\mathbb{R})$ being the Pochhammer symbol
\begin{align*}
(a)_n=
\begin{cases}
1,&n=0,\\
a(a+1)\cdots(a+n-1),&n=1,2,\ldots.
\end{cases}
\end{align*}
\end{proposition}
\begin{remark}
The condition $\theta\in[0,2\sigma\pi]$ in \eqref{eigen-f} comes from the fact that functions defined on the circle $\mathbb{S}_\sigma^1$ can be identified with $2\sigma\pi$-periodic functions on $\mathbb{R}$. Let $L^\alpha_m(t)$ be the generalized Laguerre polynomial
\begin{equation*}
L^\alpha_m(t)=\sum_{n=0}^m (-1)^n \Bigg(
  \begin{array}{c}
    m+\alpha \\
    m-n \\
  \end{array}
\Bigg)\frac{t^n}{n!},
\end{equation*}
then one can verify the following relation
\begin{equation}\label{P-L}
\begin{split}
P_{k,m}(\tilde{r})&=\sum_{n=0}^m\frac{(-1)^n m(m-1)\cdots(m-(n-1))}{(\alpha_k +1)(\alpha_k+2)\cdots(\alpha_k+n)}\frac{\tilde{r}^n}{n!}\\
&=\Bigg(
  \begin{array}{c}
    m+\alpha_k \\
    m \\
  \end{array}
\Bigg)^{-1}L^{\alpha_k}_m(\tilde{r}),\quad \alpha_k:=\left|\frac{k}{\sigma}+\alpha\right|.
\end{split}
\end{equation}
Moreover, by the well-known orthogonality of the generalized Laguerre polynomials
\begin{equation}\label{rel:orth}
\int_0^\infty x^{\alpha}e^{-x}L^\alpha_m(x)L^\alpha_n(x)\mathrm{d}x=\frac{\Gamma(n+\alpha+1)}{n!}\delta_{n,m},
\end{equation}
one can easily verify the following orthogonality
\begin{equation*}
\int_0^\infty\int_0^{2\sigma\pi}V_{k_1,m_1}(r,\theta)\overline{V_{k_2,m_2}(r,\theta)}\mathrm{d}\theta r\mathrm{d}r=0,\quad \text{if}\quad (k_1, m_1)\neq (k_2, m_2)
\end{equation*}
and calculate the $L^2$-norm of the eigenfunction $V_{k,m}(r,\theta)$
\begin{equation}\label{V-km-2}
\begin{split}
\|V_{k,m}\|^2_{L^2(X)}&=\int_0^\infty\int_0^{2\sigma\pi}|V_{k,m}(r,\theta)|^2\mathrm{d}\theta r\mathrm{d}r\\
&=\Big(\frac{2}{B_0}\Big)^{\alpha_k+1}\Gamma(1+\alpha_k)\Bigg(
  \begin{array}{c}
    m+\alpha_k \\
    m \\
  \end{array}
\Bigg)^{-1},
\end{split}
\end{equation}
where $\delta_{n,m}$ denotes the Kronecker delta symbol.
A significant fact we would like to stress is that, in view of \eqref{rel:orth}, the family of the eigenfunctions $\{V_{k,m}\}_{k\in\mathbb{Z},m\in\mathbb{N}}$ in \eqref{eigen-f} forms a complete orthogonal basis for $L^2(X)$.
\end{remark}
\begin{proof}[Proof of Proposition \ref{prop:spect}]
Consider the eigenvalue problem for the operator $H_{\alpha,B_0,\sigma}$
\begin{equation}\label{eigen-p}
H_{\alpha,B_0,\sigma}g(r,\theta)=\lambda g(r,\theta),\quad (r,\theta)\in(0,+\infty)_r\times\mathbb{S}_\sigma^1,\quad g\in\mathcal{D},
\end{equation}
where $\mathcal{D}$ denotes the domain of the self-adjoint operator $H_{\alpha,B_0,\sigma}$ given by \eqref{domain}.

Define the projection $P_k$ by
\begin{equation*}
P_kg(r,\theta)=\frac{1}{2\sigma\pi}\int_0^{2\sigma\pi}e^{i\frac{k}{\sigma}(\theta-\theta')}g(r,\theta')\mathrm{d}\theta',\quad k\in\mathbb{Z},
\end{equation*}
then the operator $H_{\alpha,B_0,\sigma}$ commutes with the projection $P_k$ and the problem \eqref{eigen-p} is reduced to
\begin{equation}\label{eq:gk}
g''_k(r)+\frac{1}{r}g'_k(r)-\frac{1}{r^2}\left(\frac{k}{\sigma}+\alpha+\frac{B_0r^2}{2}\right)^2g_k(r)=-\lambda g_k(r),\quad k\in\mathbb{Z},
\end{equation}
where
\begin{equation*}
g_k(r)=\frac{1}{\sqrt{2\sigma\pi}}\int_0^{2\sigma\pi}e^{-i\frac{k}{\sigma}\theta}g(r,\theta)\mathrm{d}\theta.
\end{equation*}
Let $\alpha_k:=|k/\sigma+\alpha|$ and $\phi_k(s)=\Big(\frac{2s}{B_0}\Big)^{-\frac{\alpha_k}{2}}e^{\frac{s}{2}}g_k\Big(\sqrt{\frac{2s}{B_0}}\Big)$, then it is not difficult to verify that $\phi_k$ satisfies
\begin{equation}\label{eq:phik}
s\phi''_k(s)+(1+\alpha_k-s)\phi'_k(s)-\frac{1}{2}\left(1+\alpha_k+\frac{k}{\sigma}+\alpha-\frac{\lambda}{B_0}\right)\phi_k(s)=0.
\end{equation}
This is the standard Kummer Confluent Hypergeometric equation, which is a solvable model in quantum mechanics (see e.g. \cite{AS65}).
\begin{lemma}\label{lem:KCH}
The Kummer Confluent Hypergeometric equation
\begin{equation*}
s\phi''(s)+(b-s)\phi'(s)-a\phi(s)=0,\quad s>0
\end{equation*}
admits at most two linearly independent solutions given by the Kummer function (or the confluent hypergeometric function of the first kind)
\begin{equation*}
M(a,b,s)=\sum_{n=0}^\infty\frac{(a)_n}{(b)_n}\frac{s^n}{n!},\quad b\neq0,-1,-2,\cdots
\end{equation*}
and the Tricomi function (or the confluent hypergeometric function of the second kind)
\begin{equation*}
U(a,b,s)=\frac{\Gamma(1-b)}{\Gamma(a-b+1)}M(a,b,s)+\frac{\Gamma(b-1)}{\Gamma(a)}s^{1-b}M(a-b+1,2-b,s).
\end{equation*}
\end{lemma}
Hence, the ODE \eqref{eq:gk} admits at most two linearly independent solutions
\begin{align*}
g_k^1(\lambda;r)&=r^{\alpha_k}M\Big(\beta(k,\lambda),1+\alpha_k,\frac{B_0r^2}{2}\Big)e^{-\frac{B_0r^2}{4}}\\
g_k^2(\lambda;r)&=r^{\alpha_k}U\Big(\beta(k,\lambda),1+\alpha_k,\frac{B_0r^2}{2}\Big)e^{-\frac{B_0r^2}{4}},
\end{align*}
where
\begin{equation*}
\beta(k,\lambda)=\frac{1}{2}\left(1+\alpha_k+\frac{k}{\sigma}+\alpha-\frac{\lambda}{B_0}\right).
\end{equation*}
The general solution of \eqref{eq:gk} is then given by
\begin{align}\label{sol:gk}
g_k(r)&=A_kg_k^1(\lambda;r)+B_kg_k^2(\lambda;r)=r^{\alpha_k}e^{-\frac{B_0r^2}{4}}\times\\
&\times\Bigg(A_kM\Big(\beta(k,\lambda),1+\alpha_k,\frac{B_0r^2}{2}\Big)+B_kU\Big(\beta(k,\lambda),1+\alpha_k,\frac{B_0r^2}{2}\Big)\Bigg),\nonumber
\end{align}
where $A_k,B_k$ are two constants. The determination of the constants $A_k,B_k$ requires the asymptotic properties of the confluent hypergeometric functions, which are summarised in the next lemma.
\begin{lemma}\label{lem:asy}
\cite[Chap.13]{AS65}
\begin{itemize}

\item The confluent hypergeometric functions $M(a,b,z)$ and $U(a,b,z)$ are linearly dependent if and only if $-a\in\mathbb{N}$.

\item $M(a,b,z)$ is an entire function with respect to $z$ and is regular at $z=0$, while $U(a,b,z)$ is singular at $z=0$ provided $b>1$ and $-a\notin\mathbb{N}$.
Moreover, it holds
    \begin{equation}\label{U-0}
    \lim_{z\rightarrow0^+}z^{b-1}U(a,b,z)=\frac{\Gamma(b-1)}{\Gamma(a)}.
    \end{equation}
    If $b\in(1,2)$, the asymptotic expression of $U(a,b,z)$ as $z\rightarrow0^+$ reads
    \begin{equation*}
    U(a,b,z)=\frac{\Gamma(1-b)}{\Gamma(a-b+1)}+\frac{\Gamma(b-1)}{\Gamma(a)}z^{1-b}+O(z^{2-b}).
    \end{equation*}
\item If $-a\notin\mathbb{N}$, then one has, as $z\rightarrow+\infty$,
    \begin{equation}\label{M-inf}
    M(a,b,z)=\frac{\Gamma(b)}{\Gamma(b-a)}(-z)^{-a}(1+O(z^{-1}))+\frac{\Gamma(b)}{\Gamma(a)}e^z z^{a-b}(1+O(z^{-1}))
    \end{equation}
    and
    \begin{equation*}
    U(a,b,z)=z^{-a}(1+O(z^{-1})).
    \end{equation*}
\end{itemize}
\end{lemma}
Now we apply Lemma \ref{lem:asy} to determine the constants $A_k,B_k$ in \eqref{sol:gk}.

For the moment, we denote
\begin{equation*}
m:=-\beta(k,\lambda)=-\frac{1}{2}\Big(1+\alpha_k+\frac{k}{\sigma}+\alpha-\frac{\lambda}{B_0}\Big).
\end{equation*}
If $m\notin\mathbb{N}$, it follows from \eqref{M-inf} that
\begin{equation*}
M\big(-m,1+\alpha_k,\tilde{r}\big)\sim\tilde{r}^{-m-1-\alpha_k}e^{\tilde{r}},\quad \tilde{r}\rightarrow+\infty,
\end{equation*}
which implies that $M\big(-m,1+\alpha_k,\tilde{r}\big)$ is not integrable near infinity.
If $m\in\mathbb{N}$, then $M(-m,1+\alpha_k,\tilde{r})$ is exactly a polynomial of degree $m$ with respect to $\tilde{r}$
\begin{equation*}
P_{k,m}(\tilde{r}):=M(-m,1+\alpha_k,\tilde{r})=\sum_{n=0}^m\frac{(-m)_n}{(1+\alpha_k)_n}\frac{\tilde{r}^n}{n!}.
\end{equation*}
In view of \eqref{U-0} and $\alpha_k=|k/\sigma+\alpha|\geq0$, we obtain
\begin{equation*}
U\big(\beta(k,\lambda),1+\alpha_k,\tilde{r}\big)\sim \tilde{r}^{-\alpha_k},\quad \tilde{r}\rightarrow0+.
\end{equation*}
Therefore, we have
\begin{equation*}
g_k(r)\sim B_k r^{\alpha_k}e^{-\frac{B_0r^2}{4}}
U\Big(\beta(k,\lambda),1+\alpha_k,\frac{B_0r^2}{2}\Big)\sim B_kr^{-\alpha_k}\quad\text{as}\quad r\rightarrow0+
\end{equation*}
and
\begin{equation*}
g_k(r)\sim A_k r^{\alpha_k}e^{-\frac{B_0r^2}{4}}
M\Big(\beta(k,\lambda),1+\alpha_k,\frac{B_0r^2}{2}\Big)\sim A_k e^{\frac{B_0r^2}{4}}r^{-1+\frac{k}{\sigma}+\alpha-\frac{\lambda}{B_0}}\quad\text{as}\quad r\rightarrow \infty.
\end{equation*}
Now we claim that $B_k\equiv0$. If otherwise, due to $g\in\mathcal{D}$, we will get
\begin{equation*}
\int_0^\infty g_k^2(r)\frac{\mathrm{d}r}{r}<+\infty,
\end{equation*}
which is a contradiction since
\begin{equation*}
\int_0^\infty r^{-2\alpha_k-1}\mathrm{d}r=+\infty,\quad \forall k\in\mathbb{Z}.
\end{equation*}
Next, we claim that $A_k\equiv0$ unless $m\in\mathbb{N}$. In fact, due to $g\in\mathcal{D}$, we also have
\begin{equation*}
\int_0^\infty g_k^2(r)r\mathrm{d}r<+\infty,
\end{equation*}
which is a contradiction as well since, in view of $B_0>0$, one has
\begin{equation*}
\int_0^\infty g_k^2(r)rdr\geq\int_1^\infty e^{\frac{B_0r^2}{4}}rdr=+\infty.
\end{equation*}
Above all, we conclude that
\begin{equation*}
\mathbb{N}\ni m=-\frac{1}{2}\left(1+\alpha_k+\frac{k}{\sigma}+\alpha-\frac{\lambda}{B_0}\right)
\end{equation*}
and
\begin{equation*}
g_k(r)=r^{\alpha_k}e^{-\frac{B_0r^2}{4}}P_{k,m}\left(\frac{B_0r^2}{2}\right).
\end{equation*}
Therefore, we obtain the eigenfunction of the operator $H_{\alpha,B_0,\sigma}$
\begin{equation*}
V_{k,m}(r,\theta)=\frac{1}{\sqrt{2\sigma\pi}}r^{\alpha_k}e^{-\frac{B_0r^2}{4}}P_{k,m}\Big(\frac{B_0r^2}{2}\Big)e^{i\frac{k}{\sigma}\theta},\quad \alpha_k=\left|\frac{k}{\sigma}+\alpha\right|.
\end{equation*}
From $-m=\frac{1}{2}\Big(1+\alpha_k+\frac{k}{\sigma}+\alpha-\frac{\lambda}{B_0}\Big)$, we obtain the corresponding eigenvalue
\begin{equation*}
\lambda_{k,m}=\left(2m+1+\alpha_k+\frac{k}{\sigma}+\alpha\right)B_0,\quad k\in\mathbb{Z},\quad m\in\mathbb{N}.
\end{equation*}
\end{proof}

\subsection{Sobolev spaces associated to $H_{\alpha,B_0,\sigma}$}

In this subsection, we will prove the equivalence of two norms.
\begin{proposition}[Equivalent norms]\label{prop:E-S-B}
Let the Sobolev norm and Besov norm be defined by \eqref{Sobolev2} and \eqref{Besov} respectively, then, for $s\in\mathbb{R}$, there exist two positive constants $c,C>0$ such that
\begin{equation}\label{S-B-H}
c\|f\|_{\dot \HH^s_{\alpha,B_0,\sigma}(X)}\leq \|f\|_{\dot{\mathcal{B}}^s_{2,2,\sigma}(X)}\leq C \|f\|_{\dot \HH^s_{\alpha,B_0,\sigma}(X)}
\end{equation}
and
\begin{equation}\label{S-B-H'}
c\|f\|_{\HH^s_{\alpha,B_0,\sigma}(X)}\leq \|f\|_{{\mathcal{B}}^s_{2,2,\sigma}(X)}\leq C \|f\|_{ \HH^s_{\alpha,B_0,\sigma}(X)}.
\end{equation}
\end{proposition}
\begin{proof}
Let $\tilde{V}_{k,m}$ be the $L^2$-normalization of the eigenfunction $V_{k,m}$ in \eqref{eigen-f}, then it follows by direct verification that the family of these eigenfunctions $\Big\{\tilde{V}_{k,m}\Big\}_{k\in\mathbb{Z}, m\in\mathbb{N}}$ forms a complete orthonormal basis for $L^2(X)$.

By functional calculus, for any well-behaved functions $F$ (e.g. bounded Borel measurable function) and $f\in L^2(X)$, we can write
\begin{equation*}
\big(F(H_{\alpha,B_0,\sigma})f\big)(r,\theta)=\sum_{k\in\mathbb{Z} \atop m\in\mathbb{N}}F(\lambda_{k,m})c_{k,m}\tilde{V}_{k,m}(r,\theta),
\end{equation*}
where
\begin{equation*}
c_{k,m}=\int_0^\infty\int_0^{2\sigma\pi}f(r,\theta)\overline{\tilde{V}_{k,m}(r,\theta)}\mathrm{d}\theta r\mathrm{d}r.
\end{equation*}
Then, it follows by the $L^2$-orthogonality between different eigenfunctions
\begin{equation}\label{norm-f}
\|F(H_{\alpha,B_0,\sigma})f\|_{L^2(X)}=\left(\sum_{k\in\mathbb{Z}\atop m\in\mathbb{N}}\big|F(\lambda_{k,m})c_{k,m}\big|^2\right)^{1/2}.
\end{equation}
In particular, it holds
\begin{equation*}
\|f\|_{\dot \HH^s_{\alpha,B_0,\sigma}(X)}=\|H^{s/2}_{\alpha,B_0,\sigma} f\|_{L^2(X)}=\left(\sum_{k\in\mathbb{Z} \atop m\in\mathbb{N}} \big|\lambda^{s/2}_{k,m} c_{k,m}\big|^2\right)^{1/2}.
\end{equation*}
Take $\varphi\in C_c^\infty(\mathbb{R}_+)$ as in \eqref{LP-dp}. On the one hand, we have by the definition of Besov norm and \eqref{norm-f}
\begin{equation*}
\begin{split}
\|f\|_{\dot{\mathcal{B}}^s_{2,2,\sigma}(X)}&=\left(\sum_{j\in\mathbb{Z}}2^{2js}\|\varphi_j(\sqrt{H_{\alpha,B_0,\sigma}})f\|_{L^2(X)}^2\right)^{1/2}\\
&=\left(\sum_{j\in\mathbb{Z}}\sum_{k\in\mathbb{Z} \atop m\in\mathbb{N}}2^{2js}\left|\varphi\left(\sqrt{\lambda_{k,m}}/2^j\right)c_{k,m}\right|^2\right)^{1/2}\\
&\approx\left(\sum_{j\in\mathbb{Z}}\sum_{k\in\mathbb{Z} \atop m\in\mathbb{N}}\lambda_{k,m}^s\left|\varphi\left(\sqrt{\lambda_{k,m}}/2^j\right)c_{k,m}\right|^2\right)^{1/2}\\
&\lesssim\left(\sum_{k\in\mathbb{Z} \atop m\in\mathbb{N}}|\lambda_{k,m}^{s/2}c_{k,m}|^2\sum_{j\in\mathbb{Z}}\left|\varphi\left(\sqrt{\lambda_{k,m}}/2^j\right)\right|^2\right)^{1/2}\\
&\lesssim\left(\sum_{k\in\mathbb{Z} \atop m\in\mathbb{N}}|\lambda_{k,m}^{s/2}c_{k,m}|^2\right)^{1/2}=\|f\|_{\dot \HH^s_{\alpha,B_0,\sigma}(X)}.
\end{split}
\end{equation*}
On the other hand, we have
\begin{equation*}
\begin{split}
\|f\|_{\dot \HH^s_{\alpha,B_0,\sigma}(X)}&=\left(\sum_{k\in\mathbb{Z} \atop m\in\mathbb{N}}|\lambda_{k,m}^{s/2}c_{k,m}|^2\right)^{1/2}\\
&=\left(\sum_{k\in\mathbb{Z} \atop m\in\mathbb{N}}\left|\sum_{j\in\mathbb{Z}}\varphi\left(\sqrt{\lambda_{k,m}}/2^j\right)\lambda_{k,m}^{s/2}c_{k,m}\right|^2\right)^{1/2}\\
&\lesssim\left(\sum_{j\in\mathbb{Z}}\sum_{k\in\mathbb{Z} \atop m\in\mathbb{N}}2^{2js}\left|\varphi\left(\sqrt{\lambda_{k,m}}/2^j\right)c_{k,m}\right|^2\right)^{1/2}\\
&\lesssim\left(\sum_{j\in\mathbb{Z}}2^{2js}\|\varphi_j(\sqrt{H_{\alpha,B_0,\sigma}})f\|_{L^2(X)}^2\right)^{1/2}=\|f\|_{\dot{\mathcal{B}}^s_{2,2,\sigma}(X)},
\end{split}
\end{equation*}
where we have used the fact that, for a fixed $\lambda>0$, there are only finite many non-vanishing terms in the partition of unity
\begin{equation*}
1=\sum_{j\in\mathbb{Z}}\varphi\left(\lambda/2^j\right).
\end{equation*}
Therefore, the proof of \eqref{S-B-H} is finished. One can similarly verify \eqref{S-B-H'}.
\end{proof}

\section{Weight dispersive estimate for Schr\"odinger equation}\label{sec:const}

The aim of this section is to establish the weighted estimate \eqref{decay:S} for the Schr\"odinger equation \eqref{eq:S} associated to the operator \eqref{operator} on the product cone $(X,g)$. Since the weighted estimate \eqref{decay:S} essentially follows by the usual dispersive estimate \eqref{dis:S}, the main efforts will be devoted to the proof of the dispersive estimate \eqref{dis:S} and we shall achieve this goal via three different methods. The idea of the first and second approaches is to construct an explicit representation formula for the kernel of the associated Schr\"odinger propagator $e^{itH_{\alpha,B_0,\sigma}}$ and these two approaches are related by applying the universal covering space technique originally developed by Schulman \cite{Sch68} in connection with the Feynman path integral on multiply connected spaces. The third approach is based on directly analysis of the uniform convergence of the series form of the associated Schr\"odinger kernel obtained via the spectral theorem.

Let us briefly introduce the background of applying Schulman's ansatz. It is easy to see that a natural universal covering space of the product cone $X=(0,+\infty)_r\times\mathbb{S}_\sigma^1$ is given by the right-half Euclidean plane $\tilde{X}=(0,+\infty)_r\times\mathbb{R}$ since there exists a natural bijection between the planar circle $\mathbb{S}_\sigma^1$ and the real line $\mathbb{R}$. Let $\Gamma=2\sigma\pi\mathbb{Z}$ be the group acting on the second entry of the Cartesian product $(0,+\infty)_r\times\mathbb{S}_\sigma^1$, then we have $X\simeq\tilde{X}/\Gamma$ (one knows $\mathbb{S}_\sigma^1\cong\mathbb{T}_\sigma^1\cong\mathbb{R}/2\sigma\pi\mathbb{Z}$). Schulman's ansatz leads us to obtain a connection formula for the Schr\"odinger kernel of $H_{\alpha,B_0,\sigma}$ on $X$ and on $\tilde{X}$. More precisely, we shall present the first two approaches by verifying the following formula (see \cite[(1)]{Sto89} for the abstract form)
\begin{equation}\label{SS}
e^{itH_{\alpha,B_0,\sigma}}(r_1,\theta_1;r_2,\theta_2)=\sum_{j\in\mathbb{Z}}e^{it\tilde{H}_{\alpha,B_0,\sigma}}(r_1,\theta_1-2j\sigma\pi;r_2,\theta_2),
\end{equation}
where $\tilde{H}_{\alpha,B_0,\sigma}$ denotes the operator having the same expression as $H_{\alpha,B_0,\sigma}$ but acting on $\tilde{X}$.
The formula \eqref{SS} is quite similar to the one of the wave kernels on $\mathbb{R}^n$ and $\mathbb{T}^n$ (see \cite[(3.5.12)]{sogge}).

All the three approaches are based on the series representation of the Schr\"odinger kernel obtained via the spectral theorem.
\begin{proposition}\label{prop:S-series}
Let $H_{\alpha,B_0,\sigma}$ be the Hamiltonian in \eqref{op} and let $K_t^S(r_1,\theta_1;r_2,\theta_2)$ denote the kernel of the Schr\"odinger propagator $e^{itH_{\alpha,B_0,\sigma}}$. If we identify functions defined on the planar circle $\mathbb{S}_\sigma^1$ of radius $\sigma>0$ with $2\sigma\pi$-periodic functions defined on $\mathbb{R}$ and denote $\alpha_k=|k/\sigma+\alpha|$, we will have
\begin{equation}\label{Sch:series}
\begin{split}
K_t^S(r_1,\theta_1;r_2,\theta_2)&=\frac{iB_0e^{itB_0\alpha}}{8\pi\sigma\sin(tB_0)}e^{\frac{B_0(r_1^2+r_2^2)}{4i\tan(tB_0)}}\\
&\qquad \times\sum_{k\in\mathbb{Z}}e^{i\frac{k}{\sigma}(tB_0-(\theta_1-\theta_2))}I_{\alpha_k}\left(\frac{iB_0r_1r_2}{2\sin(tB_0)}\right),
\end{split}
\end{equation}
where $I_\nu(z)$ denotes the modified Bessel function of the first kind of order $\nu$ with respect to the variable $z$.
\end{proposition}
\begin{remark}
In view of the periodicity feature of the sine factor in \eqref{Sch:series}, we shall restrict ourselves to consider the time interval $t\in(0,\pi/B_0)$ in what follows if no otherwise specified. We stress that the expression \eqref{Sch:series} is valid for any $\sigma>0$ but an explicit formula for the Schr\"odinger kernel is available only for the restricted case $\sigma\geq1$. In addition, the series in \eqref{Sch:series} is actually of the following form (in view of the relation $I_\nu(r)=i^{-\nu}J_\nu(ir)$)
\begin{equation}\label{series}
K(r,\delta):=\sum_{k\in\mathbb{Z}}i^{-\alpha_k}e^{i\frac{k}{\sigma}\delta}J_{\alpha_k}(r),\quad \alpha_k=|k/\sigma+\alpha|
\end{equation}
and the series \eqref{series} should be understood in the sense that there exists some integer $k_0\in\mathbb{N}$ such that
\begin{equation}\label{tail}
K_{\geq}(r,\delta):=\sum_{k\in\mathbb{Z}:|k|\geq1+k_0}i^{-\alpha_k}e^{i\frac{k}{\sigma}\delta}J_{\alpha_k}(r)
\end{equation}
is uniformly convergent for any $r\in\mathbb{R}_+$ and
\begin{equation*}
K(r,\delta)-\sum_{k\in\mathbb{Z}:|k|\leq k_0}i^{-\alpha_k}e^{i\frac{k}{\sigma}\delta}J_{\alpha_k}(r)\in L_{loc}^\infty(\mathbb{R}_+,\mathbb{C}),
\end{equation*}
where $J_\nu(r)$ is the Bessel function of the first kind. In fact, it is not difficult to see that $i^{-\alpha_k}e^{i\frac{k}{\sigma}\delta}J_{\alpha_k}(r)$ belongs to $L_{loc}^\infty(\mathbb{R}_+,\mathbb{C})$ for suitably large $|k|$ and for some fixed $R<+\infty$, it holds $|er|\leq1+\alpha_k$ for any $r\leq R$ and suitably large $|k|$.
Hence, for all $r\leq R$ and suitably large $|k|$, we have
\begin{align*}
|i^{-\alpha_k}e^{i\frac{k}{\sigma}\delta}J_{\alpha_k}(r)|&\lesssim\sum_{m=0}^\infty\frac{(r/2)^{2m+\alpha_k}}{m!\Gamma(m+1+\alpha_k)}\\
&\lesssim\frac{(r/2)^{\alpha_k}}{\Gamma(1+\alpha_k)}e^{(r/2)^2}\\
&\lesssim\left(\frac{er}{2(1+\alpha_k)}\right)^{\alpha_k}e^{(r/2)^2}\\
&\lesssim2^{-c|k|}=:M_k,\quad \text{for some }\quad c:=c(\sigma)>0
\end{align*}
and thereby the Weierstrass M-test gives the desired uniform convergence of the series \eqref{tail} (since $\sum_{k=0}^\infty2^{-ck}$ is obviously convergent), where the implicit constant depends on $R$ but not on $k$.
\end{remark}
\begin{proof}[Proof of Proposition \ref{prop:S-series}]
Let $\tilde{V}_{k,m}$ be the $L^2$-normalization of $V_{k,m}$ in \eqref{eigen-f}, then it is not hard to verify that (using the orthogonal relation \eqref{rel:orth}) the family of the eigenfunctions $\Big\{\tilde{V}_{k,m}(r,\theta)\Big\}_{k\in\mathbb{Z}, m\in\mathbb{N}}$ forms a complete orthonormal basis for $L^2(X)$.
We express $f(r_1,\theta_1)\in L^2(X)$ in terms of the basis as
\begin{equation*}
f(r_1,\theta_1)=\sum_{k\in\mathbb{Z} \atop m\in\mathbb{N}}c_{k,m}\tilde{V}_{k,m}(r_1,\theta_1),
\end{equation*}
where
\begin{equation}\label{cmk1}
c_{k,m}=\int_0^\infty\int_0^{2\sigma\pi}f(r_2,\theta_2)\overline{\tilde{V}_{k,m}(r_2,\theta_2)}\mathrm{d}\theta_2r_2\mathrm{d}r_2.
\end{equation}
By the spectral theorem, we can represent the Schr\"odinger propagator $e^{itH_{\alpha,B_0,\sigma}}$ via
\begin{equation*}
v(t,r_1,\theta_1):=\left(e^{itH_{\alpha,B_0,\sigma}}f\right)(r_1,\theta_1)=\sum_{k\in\mathbb{Z}\atop m\in\mathbb{N}}c_{k,m}e^{it\lambda_{k,m}}\tilde{V}_{k,m}(r_1,\theta_1).
\end{equation*}
Inserting \eqref{cmk1} into the above expression yields
\begin{equation*}
v(t,r_1,\theta_1)=\sum_{k\in\mathbb{Z}\atop m\in\mathbb{N}}e^{it\lambda_{k,m}}\left(\int_0^\infty\int_0^{2\sigma\pi}f(r_2,\theta_2)
\overline{\tilde{V}_{k,m}(r_2,\theta_2)}\mathrm{d}\theta_2r_2\mathrm{d}r_2\right)\tilde{V}_{k,m}(r_1,\theta_1).
\end{equation*}
A function $f\in L^2(X)$ can be decomposed alternatively as
\begin{equation*}
f(r_2,\theta_2)=\sum_{k\in\mathbb{Z}}f_k(r_2)e^{i\frac{k}{\sigma}\theta_2},
\end{equation*}
where
\begin{equation}\label{f-k}
f_k(r_2)=\frac{1}{2\sigma\pi}\int_0^{2\sigma\pi}f(r_2,\theta'_2)e^{-i\frac{k}{\sigma}\theta'_2}\mathrm{d}\theta'_2.
\end{equation}
Hence, by \eqref{eigen-v}, \eqref{eigen-f} and \eqref{V-km-2}, we obtain
\begin{align*}
v(t,r_1,\theta_1)&=\sum_{k\in\mathbb{Z} \atop m\in\mathbb{N}}e^{it\lambda_{k,m}}\frac{V_{k,m}(r_1,\theta_1)}{\|V_{k,m}\|^2_{L^2(X)}}\Bigg(\int_0^\infty f_k(r_2)e^{-\frac{B_0r_2^2}{4}}P_{k,m}\Big(\frac{B_0r_2^2}{2}\Big)r_2^{1+\alpha_k}\mathrm{d}r_2\Bigg)\\
&=\frac{B_0}{2}\sum_{k\in\mathbb{Z}}e^{i\frac{k}{\sigma}\theta_1}\frac{B_0^{\alpha_k}e^{it\beta_k}}{2^{\alpha_k}\Gamma(1+\alpha_k)}\Bigg[\sum_{m=0}^\infty
\Bigg(
  \begin{array}{c}
    m+\alpha_k \\
    m \\
  \end{array}
\Bigg)e^{2itmB_0}\\
&\times\Bigg(\int_0^\infty f_k(r_2)(r_1r_2)^{\alpha_k}e^{-\frac{B_0(r_1^2+r_2^2)}{4}}P_{k,m}\left(\frac{B_0r_2^2}{2}\right)P_{k,m}\left(\frac{B_0r_1^2}{2}\right)r_2 \mathrm{d}r_2\Bigg)\Bigg],
\end{align*}
where
\begin{equation}\label{beta-k}
\alpha_k=\left|\frac{k}{\sigma}+\alpha\right|,\quad \lambda_{k,m}=2mB_0+\beta_k,\quad \beta_k=\left(1+\alpha_k+\frac{k}{\sigma}+\alpha\right)B_0.
\end{equation}
By the integral representation of $P_{k,m}$ in terms of the Bessel function $J_\nu$ of order $\nu$ (see e.g. \cite[(6.2.15)]{AAR01})
\begin{equation}\label{P-km}
P_{k,m}\left(\frac{r^2}{2}\right)=\frac{\Gamma(1+\alpha_k)}{\Gamma(1+\alpha_k+m)}e^{\frac{r^2}{2}}r^{-\alpha_k}2^{\frac{\alpha_k}{2}}\int_0^\infty
e^{-s}s^{m+\frac{\alpha_k}{2}}J_{\alpha_k}(\sqrt{2s}r)\mathrm{d}s,
\end{equation}
we further get
\begin{align*}
v(t,r_1,\theta_1)=&\frac{B_0}{2}\sum_{k\in\mathbb{Z}}e^{i\frac{k}{\sigma}\theta_1}\frac{B_0^{\alpha_k}e^{it\beta_k}}{2^{\alpha_k}\Gamma(1+\alpha_k)}
\Bigg[\sum_{m\in\mathbb{N}}
\left(
  \begin{array}{c}
    m+\alpha_k \\
    m \\
  \end{array}
\right)e^{2itmB_0}\\
&\times\left(\frac{\Gamma(1+\alpha_k)}{\Gamma(1+\alpha_k+m)}\right)^2\Bigg(\int_0^\infty f_k(r_2)(r_1r_2)^{\alpha_k} e^{-\frac{B_0(r_1^2+r_2^2)}{4}}e^{\frac{B_0(r_1^2+r_2^2)}{2}}\Big(\frac{2}{B_0r_2^2}\Big)^{\frac{\alpha_k}{2}}\\
&\times\Big(\frac{2}{B_0r_1^2}\Big)^{\frac{\alpha_k}{2}}\Big(\int_0^\infty\int_0^\infty e^{-s_1-s_2}(s_1s_2)^{m+\frac{\alpha_k}{2}}
J_{\alpha_k}(\sqrt{2B_0s_1}r_1)J_{\alpha_k}(\sqrt{2B_0s_2}r_2)\mathrm{d}s_1\mathrm{d}s_2\Big)r_2\mathrm{d}r_2\Bigg)\Bigg]\\
=&\frac{B_0}{2}\sum_{k\in\mathbb{Z}}e^{i\frac{k}{\sigma}\theta_1}e^{it\beta_k}\Gamma(1+\alpha_k)\Bigg[\sum_{m\in\mathbb{N}}
\left(
  \begin{array}{c}
    m+\alpha_k \\
    m \\
  \end{array}
\right)\frac{e^{2itmB_0}}{(\Gamma(1+\alpha_k+m))^2}\\
&\times\Bigg(\int_0^\infty f_k(r_2)e^{\frac{B_0(r_1^2+r_2^2)}{4}}\Big(\int_0^\infty\int_0^\infty e^{-s_1-s_2}(s_1s_2)^{m+\frac{\alpha_k}{2}}\\
&\times J_{\alpha_k}(\sqrt{2B_0s_1}r_1)J_{\alpha_k}(\sqrt{2B_0s_2}r_2)\mathrm{d}s_1\mathrm{d}s_2\Big)r_2\mathrm{d}r_2\Bigg)\Bigg]\\
=&2B_0\sum_{k\in\mathbb{Z}}e^{i\frac{k}{\sigma}\theta_1}e^{it\beta_k}\Bigg[\int_0^\infty f_k(r_2)e^{\frac{B_0(r_1^2+r_2^2)}{4}}e^{-i\alpha_k(tB_0+\frac{\pi}{2})}\\
&\times\Bigg(\int_0^\infty\int_0^\infty\frac{s_1s_2}{e^{s_1^2+s_2^2}}\Bigg(\sum_{m\in\mathbb{N}}
\frac{(-1)^me^{i(tB_0+\frac{\pi}{2})(2m+\alpha_k)}}{\Gamma(1+m)\Gamma(1+\alpha_k+m)}(s_1s_2)^{2m+\alpha_k}\Bigg)\\
&\times J_{\alpha_k}\big(\sqrt{2B_0}s_1r_1\big)J_{\alpha_k}\big(\sqrt{2B_0}s_2 r_2\big)\mathrm{d}s_1\mathrm{d}s_2\Bigg)r_2\mathrm{d}r_2\Bigg].
\end{align*}
Since (see e.g. \cite[(4.5.2)]{AAR01})
\begin{equation}\label{J-alpha}
\sum_{m=0}^\infty\frac{(-1)^me^{i(tB_0+\frac{\pi}{2})(2m+\alpha_k)}}{\Gamma(1+\alpha_k+m)\Gamma(1+m)}(s_1s_2)^{2m+\alpha_k}
=J_{\alpha_k}(2s_1s_2e^{i(tB_0+\frac{\pi}{2})}),
\end{equation}
we obtain
\begin{align*}
v(t,r_1,\theta_1)=2B_0\sum_{k\in\mathbb{Z}}e^{i\frac{k}{\sigma}\theta_1}e^{-it(\alpha_kB_0-\beta_k)-i\frac{\pi}{2}\alpha_k}
\int_0^\infty e^{\frac{B_0(r_1^2+r_2^2)}{4}}f_k(r_2)G_{k,t}(r_1,r_2)r_2\mathrm{d}r_2,
\end{align*}
where
\begin{equation*}
\begin{split}
G_{k,t}(r_1,r_2)=\int_0^\infty\int_0^\infty&\frac{s_1s_2}{e^{s_1^2+s_2^2}}J_{\alpha_k}(2s_1s_2e^{i(tB_0+\frac{\pi}{2})})
\\&\times J_{\alpha_k}(\sqrt{2B_0}r_1s_1)J_{\alpha_k}(\sqrt{2B_0}r_2s_2)\mathrm{d}s_1\mathrm{d}s_2.
\end{split}
\end{equation*}
Applying the following identity (see \cite[formula (1), P.395]{Wat44})
\begin{equation}\label{formula}
\begin{split}
\int_0^\infty&e^{-p^2t^2}J_\nu(at)J_\nu(bt)t\mathrm{d}t=\frac{1}{2p^2}e^{-\frac{a^2+b^2}{4p^2}}I_\nu\left(\frac{ab}{2p^2}\right),\\
&\qquad \Re\nu>-1,\quad |\arg p|<\frac{\pi}{4},\quad I_\nu(r)=e^{-\frac{1}{2}\nu\pi i}J_\nu(re^{i\frac{\pi}{2}})
\end{split}
\end{equation}
with $t=s_2,p=1,a=\sqrt{2B_0}r_2,b=2s_1e^{i(tB_0+\frac{\pi}{2})},\nu=\alpha_k$, we get
\begin{align*}
\int_0^\infty e^{-s_2^2}J_{\alpha_k}\big(\sqrt{2B_0}r_2s_2\big)&J_{\alpha_k}\Big(2s_1s_2e^{i(tB_0+\frac{\pi}{2})}\Big)s_2\mathrm{d}s_2\\
&=\frac{1}{2}e^{-\frac{B_0r_2^2+2s_1^2e^{i(2tB_0+\pi)}}{2}}I_{\alpha_k}\left(\sqrt{2B_0}r_2s_1e^{i(tB_0+\frac{\pi}{2})}\right),
\end{align*}
where $I_\nu$ denotes the modified Bessel function of the first kind of order $\nu$.

Now we can apply the formula \eqref{formula} to simplify $G_{k,t}$ as follows
\begin{align*}
G_{k,t}(r_1,r_2)&=\frac{1}{2}\int_0^\infty e^{-s_1^2}J_{\alpha_k}\left(\sqrt{2B_0}r_1s_1\right)e^{-\frac{B_0r_2^2+2s_1^2e^{i(2tB_0+\pi)}}{2}}
I_{\alpha_k}\left(\sqrt{2B_0}r_2s_1e^{i(tB_0+\frac{\pi}{2})}\right)s_1\mathrm{d}s_1\\
&=\frac{1}{4B_0}\int_0^\infty e^{-\frac{s_1^2}{2B_0}}J_{\alpha_k}(r_1s_1)e^{-\frac{B_0r_2^2+\frac{s_1^2}{B_0}e^{i(2tB_0+\pi)}}{2}}
I_{\alpha_k}\left(r_2s_1e^{i(tB_0+\frac{\pi}{2})}\right)s_1\mathrm{d}s_1\\
&=\frac{1}{4B_0}e^{-i\alpha_k\frac{\pi}{2}}e^{-\frac{B_0r_2^2}{2}}\int_0^\infty e^{-\frac{s_1^2}{2B_0}\left(1+e^{i(2tB_0+\pi)}\right)}J_{\alpha_k}(r_1s_1)J_{\alpha_k}(r_2s_1e^{itB_0})s_1\mathrm{d}s_1\\
&=\frac{1}{4B_0}e^{-i\alpha_k\frac{\pi}{2}}e^{-\frac{B_0r_2^2}{2}}\frac{B_0}{1+e^{i(2tB_0+\pi)}}
e^{-\frac{B_0r_1^2+B_0r_2^2e^{2itB_0}}{2\left(1+e^{i(2tB_0+\pi)}\right)}}I_{\alpha_k}\left(\frac{B_0r_1 r_2e^{itB_0}}{1+e^{i(2tB_0+\pi)}}\right)\\
&=\frac{1}{4(1+e^{i(2tB_0+\pi)})}e^{-i\alpha_k\frac{\pi}{2}}e^{-\frac{B_0(r_1^2+r_2^2)}{2\left(1+e^{i(2tB_0+\pi)}\right)}}I_{\alpha_k}\left(\frac{B_0r_1 r_2e^{itB_0}}{1+e^{i(2tB_0+\pi)}}\right),
\end{align*}
where we use $I_\nu(r)=e^{-i\frac{\pi}{2}\nu}J_\nu(re^{i\frac{\pi}{2}})$ in the third line and the formula \eqref{formula} with $t=s_1,p^2=\frac{1+e^{i(2tB_0+\pi)}}{2B_0},a=r_1,b=r_2e^{itB_0},\nu=\alpha_k$ in the last line.
In view of
\begin{equation*}
\frac{1}{1+e^{i(2tB_0+\pi)}}=\frac{ie^{-itB_0}}{2\sin(tB_0)},
\end{equation*}
we obtain
\begin{align*}
v(t,r_1,\theta_1)=&2B_0\frac{ie^{-itB_0}}{8\sin(tB_0)}\sum_{k\in\mathbb{Z}}e^{i\frac{k}{\sigma}\theta_1}e^{-it(\alpha_kB_0-\beta_k)}\\
&\times\Bigg[\int_0^\infty e^{\frac{B_0(r_1^2+r_2^2)}{4}}f_k(r_2)e^{-\frac{B_0(r_2^2+r_1^2)}{4}\cdot\frac{ie^{-itB_0}}{\sin(tB_0)}}I_{\alpha_k}\Bigg(\frac{iB_0r_1r_2}{2\sin (tB_0)}\Bigg)r_2\mathrm{d}r_2\Bigg]\\
=&\frac{iB_0e^{itB_0\alpha}}{4\sin(tB_0)}\sum_{k\in\mathbb{Z}}e^{i\frac{k}{\sigma}\theta_1}e^{itB_0k}\Bigg(\int_0^\infty
f_k(r_2)e^{\frac{B_0(r_1^2+r_2^2)}{4i\tan(tB_0)}}I_{\alpha_k}\Bigg(\frac{iB_0r_1r_2}{2\sin(tB_0)}\Bigg)r_2\mathrm{d}r_2\Bigg).
\end{align*}
Inserting $f_k$ given by \eqref{f-k} into the above expression, we get
\begin{align*}
v(t,r_1,\theta_1)=&\frac{iB_0e^{itB_0\alpha}}{8\pi\sigma\sin(tB_0)}\int_0^\infty\int_0^{2\sigma\pi}e^{\frac{B_0(r_1^2+r_2^2)}{4i\tan(tB_0)}}\\
&\times\sum_{k\in\mathbb{Z}}\Bigg(e^{i\frac{k}{\sigma}(tB_0-(\theta_1-\theta_2))}I_{\alpha_k}\bigg(\frac{iB_0r_1r_2}{2\sin(tB_0)}\bigg)\Bigg)
f(r_2,\theta_2)r_2\mathrm{d}r_2\mathrm{d}\theta_2,
\end{align*}
which thereby yields the desired form of the Schr\"odinger kernel \eqref{Sch:series}.
\end{proof}
We briefly describe here how to deduce the expected weighted estimate \eqref{decay:S} from the usual dispersive estimate \eqref{dis:S} (corresponding to \eqref{decay:S} with $\gamma=0$). It is known that the dispersive estimate \eqref{dis:S} reduces to the following uniform bound
\begin{equation}\label{bd:reduce}
\sup_{r_1,r_2\geq0 \atop \theta_1,\theta_2\in\mathbb{S}_\sigma^1}\left|\sum_{k\in\mathbb{Z}}e^{i\frac{k}{\sigma}(tB_0-(\theta_1-\theta_2))}
I_{\alpha_k}\bigg(\frac{iB_0r_1r_2}{2\sin(tB_0)}\bigg)\right|<+\infty.
\end{equation}
If the bound \eqref{bd:reduce} is true, then the weighted estimate \eqref{decay:S} for $\gamma\in[0,\kappa_\sigma]$ follows by \eqref{decay:S} with $\gamma=\kappa_\sigma$ and a standard interpolation. More precisely, we are reduced to prove
\begin{equation}\label{reduce:Sch}
\||\cdot|^{-\kappa_\sigma}e^{itH_{\alpha,B_0,\sigma}}|\cdot|^{-\kappa_\sigma}\|_{L^1(X)\rightarrow L^\infty(X)}\lesssim|\sin(tB_0)|^{-1-\kappa_\sigma}.
\end{equation}
By the integral representation for the modified Bessel function $I_\nu(z)$ (see \cite[(9.6.18)]{AS65})
\begin{equation*}
I_\nu(z)=\frac{(z/2)^\nu}{\pi\Gamma(1/2+\nu)}\int_{-1}^1(1-s^2)^{\nu-\frac{1}{2}}e^{zs}\mathrm{d}s,\quad z\in\mathbb{C},
\end{equation*}
we obtain an upper bound
\begin{equation*}
|I_\nu(i\rho)|\lesssim\frac{(|\rho|/2)^\nu}{\Gamma(1/2+\nu)},\quad \forall\rho\in\mathbb{R},\quad \nu\geq0.
\end{equation*}
We decompose the domain $(0,+\infty)^2$ into two parts $\Omega_1,\Omega_2$, where
\begin{equation*}
  \Omega_1:=\left\{(r_1,r_2)\in(0,+\infty)^2:r_1r_2\geq\frac{2\sin(tB_0)}{B_0}\right\}
\end{equation*}
and
\begin{equation*}
  \Omega_2:=\left\{(r_1,r_2)\in(0,+\infty)^2:r_1r_2<\frac{2\sin(tB_0)}{B_0}\right\}.
\end{equation*}
Due to $\kappa_\sigma=\text{dist}(\alpha,\sigma^{-1}\mathbb{Z})>0$, it follows from the uniform bound \eqref{bd:reduce} that
\begin{equation*}
K_1:=\sup_{r_1,r_2\in\Omega_1}\left(\frac{B_0r_1r_2}{2\sin(tB_0)}\right)^{-\kappa_\sigma}
\left|\sum_{k\in\mathbb{Z}}e^{i\frac{k}{\sigma}(tB_0-(\theta_1-\theta_2))}I_{\alpha_k}\left(\frac{iB_0r_1r_2}{2\sin(tB_0)}\right)\right|<+\infty
\end{equation*}
and
\begin{align*}
K_2:&=\sup_{r_1,r_2\in\Omega_2}\left(\frac{B_0r_1r_2}{2\sin(tB_0)}\right)^{-\kappa_\sigma}
\left|\sum_{k\in\mathbb{Z}}e^{i\frac{k}{\sigma}(tB_0-(\theta_1-\theta_2))}I_{\alpha_k}\left(\frac{iB_0r_1r_2}{2\sin(tB_0)}\right)\right|\\
&\leq\sup_{r_1,r_2\in\Omega_2}\left(\frac{B_0r_1r_2}{2\sin(tB_0)}\right)^{-\kappa_\sigma}
\sum_{k\in\mathbb{Z}}\left|I_{\alpha_k}\left(\frac{iB_0r_1r_2}{2\sin(tB_0)}\right)\right|\\
&\lesssim\sup_{r_1,r_2\in\Omega_2}\sum_{k\in\mathbb{Z}}\frac{1}{2^{\alpha_k}\Gamma(1/2+\alpha_k)}\left(\frac{B_0r_1r_2}{2\sin(tB_0)}\right)^{\alpha_k-\kappa_\sigma}.
\end{align*}
Note that the power of $\frac{B_0r_1r_2}{2\sin(tB_0)}$ in the above series is always positive in view of the definition of $\kappa_\sigma$ and $\alpha_k$, we obtain $K_2<+\infty$. Collecting the bounds of $K_1,K_2$, we have
\begin{align*}
\sup_{r_1,r_2\geq0 \atop \theta_1,\theta_2\in\mathbb{S}_\sigma^1}\left(\frac{B_0r_1r_2}{2\sin(tB_0)}\right)^{-\kappa_\sigma}
&\left|\sum_{k\in\mathbb{Z}}e^{i\frac{k}{\sigma}(tB_0-(\theta_1-\theta_2))}I_{\alpha_k}\left(\frac{iB_0r_1r_2}{2\sin(tB_0)}\right)\right|\\
&\qquad\qquad \leq(B_0/2)^{\kappa_\sigma}\max\{K_1,K_2\}|\sin(tB_0)|^{-\kappa_\sigma},
\end{align*}
which yields the estimate \eqref{reduce:Sch} and the proof of the desired weighted estimate \eqref{decay:S} is thus finished.

Therefore, it remains to prove the dispersive estimate \eqref{dis:S}.

\subsection{The first and second approaches}

Now we are in the position to verify the formula \eqref{SS} (by performing the first and second approaches subsequently). Writing $\theta=tB_0-(\theta_1-\theta_2)$ and $z=\frac{iB_0r_1r_2}{2\sin(tB_0)}$ for the moment, then we claim that the Schr\"odinger kernel $K_t^S(r_1,\theta_1;r_2,\theta_2)$ at the left-hand side of \eqref{SS} is given by
\begin{equation}\label{S:express}
\begin{split}
K_t^S(r_1,\theta_1;r_2,\theta_2)&=\frac{iB_0e^{itB_0\alpha}}{8\pi\sin(tB_0)}e^{\frac{B_0(r_1^2+r_2^2)}{4i\tan(tB_0)}}\\
&\quad \times\Bigg[\sum_{j\in\mathbb{Z}:|\theta+2j\sigma\pi|\leq\pi}e^{z\cos(\theta+2j\sigma\pi)-i\alpha(\theta+2j\sigma\pi)}\\
&\qquad\qquad -\frac{1}{2\sigma\pi}\int_0^{\infty}e^{\frac{B_0r_1r_2}{2i\sin(tB_0)}\cosh s}A_{\sigma,\alpha}(s,tB_0-\theta)\mathrm{d}s\Bigg],
\end{split}
\end{equation}
where $A_{\sigma,\alpha}(s,\theta)$ is given by
\begin{equation}\label{A-s-1-2}
\begin{split}
A_{\sigma,\alpha}(s,\theta_1-\theta_2):=\sin(\pi|\alpha|)e^{-s|\alpha|}&+\frac{1}{4i}\Bigg(\sum_{\pm}\pm e^{-s\alpha\pm i\alpha\pi}\frac{e^{\mp i\phi_\mp}-e^{-\frac{s}{\sigma}}}{\cosh(s/\sigma)-\cos\phi_\mp}\\
&\quad -\sum_{\pm}\pm e^{s\alpha\pm i\alpha\pi}\frac{e^{\pm i\phi_\mp}-e^{-\frac{s}{\sigma}}}{\cosh(s/\sigma)-\cos\phi_\mp}\Bigg)
\end{split}
\end{equation}
with $\phi_\pm$ denoting
\begin{equation}\label{+-}
\phi_\pm:=\phi_\pm(\theta_1-\theta_2)=\frac{1}{\sigma}(\theta_1-\theta_2\pm\pi).
\end{equation}
To compute the series in \eqref{Sch:series}, we will use the integral representation of the modified Bessel function $I_\nu$ (see \cite{Wat44})
\begin{equation}\label{mBessel:integral}
I_\nu(z)=\frac{1}{\pi}\int_0^\pi e^{z\cos s}\cos(s\nu)\mathrm{d}s-\frac{\sin(\nu\pi)}{\pi}\int_0^\infty e^{-z\cosh s}e^{-s\nu}\mathrm{d}s.
\end{equation}
Substituting the modified Bessel function $I_{\alpha_k}$ in \eqref{Sch:series} by \eqref{mBessel:integral}, we are reduced to address
\begin{equation}\label{equ:term1}
\frac{1}{2\sigma\pi}\sum_{k\in\mathbb{Z}}e^{i\frac{k}{\sigma}\theta}\int_0^\pi e^{z\cos s}\cos(s\alpha_k)\mathrm{d}s
\end{equation}
and
\begin{equation}\label{equ:term2}
\frac{1}{2\sigma\pi}\sum_{k\in\mathbb{Z}}e^{i\frac{k}{\sigma}\theta}\sin(\alpha_k\pi)\int_0^\infty e^{-z\cosh s}e^{-s\alpha_k}\mathrm{d}s.
\end{equation}
Following the proof of \cite[Proposition 3.1]{FZZ22}, we use the Poisson summation formula
\begin{equation*}
\sum_{j\in\mathbb{Z}}\delta(x-2j\sigma\pi)=\frac{1}{2\sigma\pi}\sum_{k\in\mathbb{Z}}e^{i\frac{k}{\sigma}x}
\end{equation*}
to obtain (recall $\alpha_k=|k/\sigma+\alpha|$)
\begin{align*}
\frac{1}{\sigma\pi}\sum_{k\in\mathbb{Z}}e^{i\frac{k}{\sigma}\theta}\cos(s\alpha_k)
&=\frac{1}{2\sigma\pi}\sum_{k\in\mathbb{Z}}e^{i\frac{k}{\sigma}\theta}(e^{is(\frac{k}{\sigma}+\alpha)}+e^{-is(\frac{k}{\sigma}+\alpha)})\\
&=\frac{1}{2\sigma\pi}\left(e^{is\alpha}\sum_{k\in\mathbb{Z}}e^{i\frac{k}{\sigma}(\theta+s)}+e^{-is\alpha}\sum_{j\in\mathbb{Z}}e^{i\frac{k}{\sigma}(\theta-s)}\right)\\
&=\sum_{j\in\mathbb{Z}}\left[e^{is\alpha}\delta(\theta+s+2j\sigma\pi)+e^{-is\alpha}\delta(\theta-s+2j\sigma\pi)\right].
\end{align*}
Hence, the first term \eqref{equ:term1} becomes
\begin{align}
&\frac{1}{2\sigma\pi}\sum_{k\in\mathbb{Z}}e^{i\frac{k}{\sigma}\theta}\int_0^\pi e^{z\cos s}\cos(s\alpha_k)\mathrm{d}s\nonumber\\
&=\frac{1}{2}\sum_{j\in\mathbb{Z}}\int_0^\pi e^{z\cos s}\Bigg(e^{is\alpha}\delta(\theta+s+2j\sigma\pi)+e^{-is\alpha}\delta(\theta-s+2j\sigma\pi)\Bigg)\mathrm{d}s\nonumber\\
&=\frac{1}{2}\sum_{j\in\mathbb{Z}}e^{z\cos(\theta+2j\sigma\pi)-i\alpha(\theta+2j\sigma\pi)}1_{[0,\pi]}(|\theta+2j\sigma\pi|)\nonumber\\
&=\frac{1}{2}\sum_{j\in\mathbb{Z}:|\theta+2j\sigma\pi|\leq\pi}e^{z\cos(\theta+2j\sigma\pi)-i\alpha(\theta+2j\sigma\pi)}.   \label{term1+2}
\end{align}
Here one can confirm the generality of the reduced assumption $\alpha\in(0,1/\sigma)$ more elementarily by the simple fact  $\{\frac{k}{\sigma}:k\in\mathbb{Z}\}\equiv\{\frac{k+1}{\sigma}:k\in\mathbb{Z}\}$. Hence, for the second term \eqref{equ:term2}, we can apply
\begin{align*}
\alpha_k=\left|\alpha+\frac{k}{\sigma}\right|=
\begin{cases}
\frac{k}{\sigma}+\alpha,&k\geq0,\\
-\alpha-\frac{k}{\sigma},&k\leq-1
\end{cases}
\end{align*}
to obtain
\begin{align*}
&\sum_{k\in\mathbb{Z}}e^{i\frac{k}{\sigma}\theta}\sin\left(\pi\left|\alpha+\frac{k}{\sigma}\right|\right)e^{-s\left|\alpha+\frac{k}{\sigma}\right|}\\
&=\sin(\pi\alpha)\sum_{k\geq0}e^{i\frac{k}{\sigma}\theta}\frac{e^{i\frac{k}{\sigma}\pi}+e^{-i\frac{k}{\sigma}\pi}}{2}e^{-s(\frac{k}{\sigma}+\alpha)}
+\cos(\pi\alpha)\sum_{k\geq0}e^{i\frac{k}{\sigma}\theta}\frac{e^{i\frac{k}{\sigma}\pi}-e^{-i\frac{k}{\sigma}\pi}}{2i}e^{-s(\frac{k}{\sigma}+\alpha)}\\
&\qquad -\sin(\pi\alpha)\sum_{k\leq-1}e^{i\frac{k}{\sigma}\theta}\frac{e^{i\frac{k}{\sigma}\pi}+e^{-i\frac{k}{\sigma}\pi}}{2}e^{s(\frac{k}{\sigma}+\alpha)}
-\cos(\pi\alpha)\sum_{k\leq-1}e^{i\frac{k}{\sigma}\theta}\frac{e^{i\frac{k}{\sigma}\pi}-e^{-i\frac{k}{\sigma}\pi}}{2i}e^{s(\frac{k}{\sigma}+\alpha)}\\
&=\frac{e^{-s\alpha}\sin(\pi\alpha)}{2}\sum_{k\geq0}\left(e^{i\frac{k}{\sigma}(\theta+is+\pi)}+e^{i\frac{k}{\sigma}(\theta+is-\pi)}\right)
+\frac{e^{-s\alpha}\cos(\pi\alpha)}{2i}\sum_{k\geq0}\left(e^{i\frac{k}{\sigma}(\theta+is+\pi)}-e^{i\frac{k}{\sigma}(\theta+is-\pi)}\right)\\
&\qquad -\frac{e^{s\alpha}\sin(\pi\alpha)}{2}\sum_{k\geq1}\left(e^{i\frac{k}{\sigma}(is-\theta-\pi)}+e^{i\frac{k}{\sigma}(is-\theta+\pi)}\right)
-\frac{e^{s\alpha}\cos(\pi\alpha)}{2i}\sum_{k\geq1}\left(e^{i\frac{k}{\sigma}(is-\theta-\pi)}-e^{i\frac{k}{\sigma}(is-\theta+\pi)}\right)\\
&=\frac{e^{i\alpha(is+\pi)}}{2i}\left(1+\sum_{k\geq1}e^{i\frac{k}{\sigma}(\theta+is+\pi)}\right)-\frac{e^{i\alpha(is-\pi)}}{2i}\left(1+\sum_{k\geq1}e^{i\frac{k}{\sigma}(\theta+is-\pi)}\right)\\
&\qquad\qquad\qquad\qquad -\frac{e^{\alpha(s+i\pi)}}{2i}\sum_{k\geq1}e^{i\frac{k}{\sigma}(is-\theta-\pi)}+\frac{e^{\alpha(s-i\pi)}}{2i}\sum_{k\geq1}e^{i\frac{k}{\sigma}(is-\theta+\pi)}\\
&=e^{-s\alpha}\sin(\alpha\pi)+\frac{1}{2i}\Bigg(\frac{e^{i\alpha(is+\pi)}e^{\frac{i}{\sigma}(is+\theta+\pi)}}{1-e^{\frac{i}{\sigma}(is+\theta+\pi)}}
-\frac{e^{i\alpha(is-\pi)}e^{\frac{i}{\sigma}(is+\theta-\pi)}}{1-e^{\frac{i}{\sigma}(is+\theta-\pi)}}\\
&\qquad\qquad \qquad-\frac{e^{\alpha(s+i\pi)}e^{\frac{i}{\sigma}(is-\theta-\pi)}}{1-e^{\frac{i}{\sigma}(is-\theta-\pi)}}
+\frac{e^{\alpha(s-i\pi)}e^{\frac{i}{\sigma}(is-\theta+\pi)}}{1-e^{\frac{i}{\sigma}(is-\theta+\pi)}}\Bigg),
\end{align*}
where we have used the following summation formula
\begin{equation*}
\sum_{k=1}^\infty e^{ikz}=\frac{e^{iz}}{1-e^{iz}},\quad \Im z>0.
\end{equation*}
In view of the identity
\begin{equation*}
\frac{e^{\pm i\phi_\pm-\frac{s}{\sigma}}}{1-e^{-\frac{s}{\sigma}\pm i\phi_\pm}}=\frac{e^{\pm i\phi_\pm}-e^{-s/\sigma}}{2(\cosh(s/\sigma)-\cos\phi_\pm)},
\end{equation*}
we obtain the explicit representation of the second term \eqref{equ:term2}
\begin{align*}
\frac{1}{2\sigma\pi}\sum_{k\in\mathbb{Z}}e^{i\frac{k}{\sigma}\theta}\sin(\alpha_k\pi)&\int_0^\infty e^{-z\cosh s}e^{-s\alpha_k}\mathrm{d}s\\
&=\frac{1}{2\sigma\pi}\int_0^\infty e^{-z\cosh s}A_{\sigma,\alpha}(s,\theta_1-\theta_2)\mathrm{d}s,
\end{align*}
where $A_{\sigma,\alpha}(s,\theta_1-\theta_2)$ is given by \eqref{A-s-1-2}.
Collecting the above calculations for the terms \eqref{equ:term1} and \eqref{equ:term2}, we obtain the desired representation formula for the Schr\"odinger kernel at the left-hand side of \eqref{SS}
\begin{equation*}
\begin{split}
K_t^S(r_1,\theta_1;r_2,\theta_2)&=\frac{iB_0e^{itB\alpha}}{8\pi\sin(tB_0)}e^{\frac{B_0(r_1^2+r_2^2)}{4i\tan(tB_0)}}\\
&\qquad\times\Bigg[\sum_{j\in\mathbb{Z}:|\theta+2j\sigma\pi|\leq\pi}e^{\frac{iB_0r_1r_2}{2\sin(tB_0)}\cos(\theta+2j\sigma\pi)-i\alpha(\theta+2j\sigma\pi)}\\
&\qquad \qquad\qquad -\frac{1}{2\sigma\pi}\int_0^\infty e^{\frac{B_0r_1r_2}{2i\sin(tB_0)}\cosh s}A_{\sigma,\alpha}(s,tB_0-\theta)\mathrm{d}s\Bigg].
\end{split}
\end{equation*}
where $\theta=tB_0-(\theta_1-\theta_2)$.

To calculate the series at the right-hand side of \eqref{SS}, we would like to make some preparations.
Recall that the operator $H_{\alpha,B_0,\sigma}$ on the product cone $X=(0,+\infty)_r\times\mathbb{S}_\sigma^1$
\begin{equation}\label{ham}
H_{\alpha,B_0,\sigma}=-\partial_{rr}-\frac{1}{r}\partial_r+\frac{1}{r^2}\Big(-i\partial_\theta+\alpha+\frac{B_0r^2}{2}\Big)^2,\quad (r,\theta)\in X
\end{equation}
acts on $L^2(X)$-functions. For the eigenvalue $\lambda_{k,m}$ in \eqref{eigen-v} and the eigenfunction $V_{k,m}$ in \eqref{eigen-f}, we have
\begin{equation*}
H_{\alpha,B_0,\sigma}V_{k,m}(r,\theta)=\lambda_{k,m}V_{k,m}(r,\theta),\quad (r,\theta)\in X.
\end{equation*}
To obtain $V_{k,m}$, after identifying functions defined on $\mathbb{S}_\sigma^1$ with $2\sigma\pi$-periodic functions on $\mathbb{R}$, we have taken $\psi_k(\theta)=\frac{1}{\sqrt{2\sigma\pi}}e^{i\frac{k}{\sigma}\theta}$ as the $L^2(\mathbb{S}_\sigma^1)$-eigenfunction of the angular operator $\Big(-i\partial_\theta+\alpha+\frac{B_0r^2}{2}\Big)^2$.
In fact, one easily verifies that the function $\psi_k(\theta)$ solves the following problem
\begin{equation*}
\begin{cases}
\Big(-i\partial_\theta+\alpha+\frac{B_0r^2}{2}\Big)^2\psi_k(\theta)=\Big(\frac{k}{\sigma}+\alpha+\frac{B_0r^2}{2}\Big)^2\psi_k(\theta),\\
\psi_k(0)=\psi_k(2\sigma\pi).
\end{cases}
\end{equation*}
To apply Schulman's ansatz, let us consider the operator $H_{\alpha,B_0,\sigma}$ acting on $L^2(\tilde{X},r\mathrm{d}r\mathrm{d}\theta)$-functions with $\tilde{X}=(0,+\infty)_r\times\mathbb{R}_\theta$. To avoid confusion, we shall denote the associated operator by $\tilde{H}_{\alpha,B_0,\sigma}$ (i.e. $\tilde{H}_{\alpha,B_0,\sigma}$ is the same as $H_{\alpha,B_0,\sigma}$ in form but defined on $L^2(\tilde{X},r\mathrm{d}r\mathrm{d}\theta)$).
Here we take $\psi(\theta)=e^{-i(\alpha-\tilde{k}/\sigma)\theta}$ as the $L^2(\mathbb{R}_\theta)$-eigenfunction of the operator $\Big(-i\partial_\theta+\alpha+\frac{B_0r^2}{2}\Big)^2$, i.e. $\psi(\theta)$ solves the problem
\begin{equation}\label{eq:ef}
\begin{cases}
\Big(-i\partial_\theta+\alpha+\frac{B_0r^2}{2}\Big)^2\psi(\theta)=\Big(\frac{\tilde{k}}{\sigma}+\frac{B_0r^2}{2}\Big)^2\psi(\theta),\\
\psi(\theta)\in L^2(\mathbb{R}_\theta).
\end{cases}
\end{equation}
Let us calculate the generalised eigenvalues and eigenfunctions of the Hamiltonian $\tilde{H}_{\alpha,B_0,\sigma}$ by suitably modifying the proof of Proposition \ref{prop:spect}.

Consider the eigenvalue problem
\begin{equation}\label{eigen-p'}
\tilde{H}_{\alpha,B_0,\sigma}f(r,\theta)=\lambda f(r,\theta),\quad (r,\theta)\in\tilde{X}=(0,+\infty)_r\times\mathbb{R}_\theta.
\end{equation}
The Fourier transform $F_{\theta\rightarrow\tilde{k}}$ is defined by
\begin{equation}\label{Fourier'}
F_{\theta\rightarrow\tilde{k}}f(r,\tilde{k}/\sigma)=\frac{1}{2\sigma\pi}\int_{\mathbb{R}}e^{i\frac{\tilde{k}}{\sigma}\theta}f(r,\theta)\mathrm{d}\theta
:=\hat{f}(r,\tilde{k}/\sigma),\quad \tilde{k}\in\mathbb{R}.
\end{equation}
Taking the Fourier transform for \eqref{eigen-p'} with respect to $\theta$, we obtain
\begin{equation}\label{eq:gk'}
\hat{f}''(r,\tilde{k}/\sigma)+\frac{1}{r}\hat{f}'(r,\tilde{k}/\sigma)-\frac{1}{r^2}\left(\frac{\tilde{k}}{\sigma}+\frac{B_0r^2}{2}\right)^2\hat{f}(r,\tilde{k}/\sigma)=-\lambda \hat{f}(r,\tilde{k}/\sigma).
\end{equation}
Let $\psi_{\tilde{k}}(s)=\Big(\frac{2s}{B_0}\Big)^{-\frac{|\tilde{k}|}{2\sigma}}e^{\frac{s}{2}}\hat{g}\Big(\sqrt{\frac{2s}{B_0}},\frac{\tilde{k}}{\sigma}\Big)$, then $\psi_{\tilde{k}}(s)$ satisfies
\begin{equation}\label{eq:psik}
s\psi''_{\tilde{k}}(s)+\left(1+\frac{|\tilde{k}|}{\sigma}-s\right)\psi'_{\tilde{k}}(s)-
\frac{1}{2}\left(1+\frac{|\tilde{k}|}{\sigma}+\frac{\tilde{k}}{\sigma}-\frac{\lambda}{B_0}\right)\psi_{\tilde{k}}(s)=0.
\end{equation}
By Lemma \ref{lem:KCH}, we obtain two linearly independent solutions for \eqref{eq:gk'}
\begin{align*}
\hat{f}_1(\lambda;r,\tilde{k}/\sigma)&=r^{\frac{|\tilde{k}|}{\sigma}}M\left(\tilde{\beta}(\tilde{k},\lambda),1+\frac{\tilde{k}}{\sigma},\frac{B_0r^2}{2}\right)e^{-\frac{B_0r^2}{4}}\\
\hat{f}_2(\lambda;r,\tilde{k}/\sigma)&=r^{\frac{|\tilde{k}|}{\sigma}}U\left(\tilde{\beta}(\tilde{k},\lambda),1+\frac{\tilde{k}}{\sigma},\frac{B_0r^2}{2}\right)e^{-\frac{B_0r^2}{4}},
\end{align*}
where
\begin{equation*}
\tilde{\beta}(\tilde{k},\lambda)=\frac{1}{2}\left(1+\frac{\tilde{k}}{\sigma}+\frac{|\tilde{k}|}{\sigma}-\frac{\lambda}{B_0}\right).
\end{equation*}
The general solution of \eqref{eq:gk'} is given by
\begin{equation*}
\hat{f}(r,\tilde{k}/\sigma)=A_{\tilde{k}}\hat{f}_1(\lambda;r,\tilde{k}/\sigma)+B_{\tilde{k}}\hat{f}_2(\lambda;r,\tilde{k}/\sigma),
\end{equation*}
where $A_{\tilde{k}}, B_{\tilde{k}}$ are two constants.
We abuse notations to write for the moment
\begin{equation}\label{cond:m}
m:=-\tilde{\beta}(\tilde{k},\lambda)=-\frac{1}{2}\left(1+\frac{\tilde{k}}{\sigma}+\frac{|\tilde{k}|}{\sigma}-\frac{\lambda}{B_0}\right).
\end{equation}
By Lemma \ref{lem:asy}, reasoning by contradiction as the proof of Proposition \ref{prop:spect}, we conclude again that $m\in\mathbb{N}$ and
\begin{equation*}
\hat{f}(r,\tilde{k}/\sigma)=r^{\frac{|\tilde{k}|}{\sigma}}e^{-\frac{B_0r^2}{4}}P_{\frac{\tilde{k}}{\sigma}-\alpha,m}\left(\frac{B_0r^2}{2}\right);
\end{equation*}
the details will be omitted here to avoid repeating. Let
\begin{equation}\label{eigen-f'}
U_m(r,\theta,\tilde{k})=r^{\frac{|\tilde{k}|}{\sigma}}e^{-\frac{B_0r^2}{4}}P_{\frac{\tilde{k}}{\sigma}-\alpha,m}\left(\frac{B_0r^2}{2}\right)e^{-i(\alpha-\tilde{k}/\sigma)\theta},\quad \tilde{k},\theta\in\mathbb{R},
\end{equation}
then $U_m(r,\theta,\tilde{k})$ provides a generalised eigenfunction for the operator $\tilde{H}_{\alpha,B_0,\sigma}$ (i.e. solves the equation \eqref{eigen-p'}) and in the same fashion as $V_{k,m}(r,\theta)$, the family of these generalized eigenfunctions
\begin{equation*}
\Big\{U_m(r,\theta,\tilde{k}): m\in\mathbb{N}, \tilde{k}\in\mathbb{R}\Big\}
\end{equation*}
forms a complete orthogonal basis for $L^2(\tilde{X},r\mathrm{d}r\mathrm{d}\theta)$.
From \eqref{cond:m}, we obtain the generalised eigenvalues of $\tilde{H}_{\alpha,B_0,\sigma}$
\begin{equation*}
\lambda_m(\tilde{k})=\left(2m+1+\frac{|\tilde{k}|+\tilde{k}}{\sigma}\right)B_0,\quad \tilde{k}\in\mathbb{R},\quad m\in\mathbb{N}.
\end{equation*}
Now we are ready to perform the calculations for the series at the right-hand side of \eqref{SS}. Indeed, let $\tilde{U}_m(r,\theta,\tilde{k})$ be the $L^2$-normalization of $U_m(r,\theta,\tilde{k})$ in \eqref{eigen-f'}, then a function $f(r,\theta)\in L^2(\tilde{X})$ can be expanded as
\begin{equation*}
f(r,\theta)=\sum_{m\in\mathbb{N}}\int_{\mathbb{R}}c_m(\tilde{k})\tilde{U}_m(r,\theta,\tilde{k})\mathrm{d}\tilde{k},
\end{equation*}
where
\begin{equation*}
c_m(\tilde{k})=\int_0^\infty\int_{\mathbb{R}}f(r,\theta)\overline{\tilde{U}_m(r,\theta,\tilde{k})}\mathrm{d}\theta r\mathrm{d}r.
\end{equation*}
Hence, the Schr\"odinger propagator $e^{it\tilde{H}_{\alpha,B_0,\sigma}}$ can be represented via
\begin{equation*}
\begin{split}
&\tilde{u}(t,r_1,\theta_1):=\left(e^{it\tilde{H}_{\alpha,B_0,\sigma}}f\right)(r_1,\theta_1)\\
&=\sum_{m\in\mathbb{N}}\int_{\mathbb{R}}e^{it\lambda_m(\tilde{k})}\left(\int_0^\infty\int_{\mathbb{R}}f(r_2,\theta_2)\overline{\tilde{U}_m(r_2,\theta_2,\tilde{k})}\mathrm{d}\theta_2 r_2\mathrm{d}r_2\right)\tilde{U}_m(r_1,\theta_1,\tilde{k})\mathrm{d}\tilde{k}.
\end{split}
\end{equation*}
Following the proof of Proposition \ref{prop:S-series}, we will obtain
\begin{equation*}
\begin{split}
\tilde{u}(t,r_1,\theta_1)=&\frac{iB_0}{8\pi\sigma\sin(tB_0)}\int_0^\infty\int_{\mathbb{R}}e^{\frac{B_0(r_1^2+r_2^2)}{4i\tan(tB_0)}}\\
&\times\int_{\mathbb{R}}\left(e^{i\left(\frac{\tilde{k}}{\sigma}-\alpha\right)(tB_0-(\theta_1-\theta_2))}
I_{\frac{|\tilde{k}|}{\sigma}}\left(\frac{iB_0r_1r_2}{2\sin(tB_0)}\right)\right)\mathrm{d}\tilde{k}f(r_2,\theta_2)r_2\mathrm{d}r_2\mathrm{d}\theta_2.
\end{split}
\end{equation*}
Hence, the kernel of the Schr\"odinger propagator $e^{it\tilde{H}_{\alpha,B_0,\sigma}}$ is given by
\begin{equation*}
\begin{split}
\tilde{K}_t^S(r_1,\theta_1;r_2,\theta_2)=&\frac{iB_0e^{i\alpha(\theta_1-\theta_2-tB_0)}}{8\pi\sigma\sin(tB_0)}e^{\frac{B_0(r_1^2+r_2^2)}{4i\tan(tB_0)}}\\
&\qquad \times\int_{\mathbb{R}}\left(e^{-i\frac{\tilde{k}}{\sigma}(\theta_1-\theta_2-tB_0)}
I_{\frac{|\tilde{k}|}{\sigma}}\left(\frac{iB_0r_1r_2}{2\sin(tB_0)}\right)\right)\mathrm{d}\tilde{k},
\end{split}
\end{equation*}
where $(r_j,\theta_j)\in\tilde{X},j=1,2$.

Similar as before, we obtain, by \eqref{mBessel:integral},
\begin{equation*}
\begin{split}
\frac{1}{\sigma\pi}\int_{\mathbb{R}}e^{i\frac{\tilde{k}}{\sigma}\theta}\int_0^\pi&e^{z\cos s}\cos(s|\tilde{k}|/\sigma)\mathrm{d}s\mathrm{d}\tilde{k}\\
&=\frac{1}{2\sigma\pi}\int_0^\pi e^{z\cos s}\int_{\mathbb{R}}e^{i\frac{\tilde{k}}{\sigma}\theta}
\left(e^{is\frac{\tilde{k}}{\sigma}}+e^{-is\frac{\tilde{k}}{\sigma}}\right)\mathrm{d}\tilde{k}\mathrm{d}s\\
&=e^{z\cos\theta}\Big(1_{[0,\pi]}(\theta)+1_{[0,\pi]}(-\theta)\Big)\\
&=e^{z\cos\theta}1_{[0,\pi]}(|\theta|)
\end{split}
\end{equation*}
and
\begin{equation*}
\begin{split}
&\frac{1}{\sigma\pi}\int_{\mathbb{R}}e^{i\frac{\tilde{k}}{\sigma}\theta}\sin\left(\pi\frac{|\tilde{k}|}{\sigma}\right)\int_0^\infty e^{-z\cosh s}e^{-s\frac{|\tilde{k}|}{\sigma}}\mathrm{d}s\mathrm{d}\tilde{k}\\
&=\frac{1}{\sigma\pi}\int_0^\infty e^{-z\cosh s}\Bigg(\int_0^\infty e^{i\frac{\tilde{k}}{\sigma}\theta}\frac{e^{i\pi\frac{\tilde{k}}{\sigma}}-e^{-i\pi\frac{\tilde{k}}{\sigma}}}{2i}e^{-\frac{s\tilde{k}}{\sigma}}\mathrm{d}\tilde{k}\\
&\qquad+\int_{-\infty}^0e^{i\frac{\tilde{k}}{\sigma}\theta}\frac{e^{-i\pi\frac{\tilde{k}}{\sigma}}-e^{i\pi\frac{\tilde{k}}{\sigma}}}{2i}e^{s\frac{\tilde{k}}{\sigma}}\mathrm{d}\tilde{k}\Bigg)\mathrm{d}s\\
&=\frac{1}{2\sigma\pi i}\int_0^\infty e^{-z\cosh s}\Bigg(\int_0^\infty\left(e^{-\frac{\tilde{k}}{\sigma}(s-i\theta-i\pi)}-e^{-\frac{\tilde{k}}{\sigma}(s-i\theta+i\pi)}\right)\mathrm{d}\tilde{k}\\
&\qquad-\int_0^{\infty}\left(e^{-\frac{\tilde{k}}{\sigma}(s+i\theta-i\pi)}-e^{-\frac{\tilde{k}}{\sigma}(s+i\theta+i\pi)}\right)\mathrm{d}\tilde{k}\Bigg)\mathrm{d}s\\
&=\frac{1}{2\pi i}\int_0^\infty e^{-z\cosh s}\Bigg(\frac{1}{s-i(\theta+\pi)}-\frac{1}{s-i(\theta-\pi)}\\
&\qquad -\frac{1}{s+i(\theta-\pi)}+\frac{1}{s+i(\theta+\pi)}\Bigg)\mathrm{d}s,
\end{split}
\end{equation*}
where $z=\frac{iB_0r_1r_2}{2\sin(tB_0)}$ and $\theta=tB_0-(\theta_1-\theta_2)$.
Therefore, it follows that
\begin{equation*}
\begin{split}
\tilde{K}_t^S(r_1,\theta_1;r_2,\theta_2)
&=\frac{iB_0e^{-i\alpha\theta}}{8\pi\sin(tB_0)}e^{\frac{B_0(r_1^2+r_2^2)}{4i\tan(tB_0)}}\Bigg(e^{z\cos\theta}1_{[0,\pi]}(|\theta|)\\
&-\frac{1}{2\pi i}\int_0^\infty e^{-z\cosh s}\Bigg(\frac{1}{s-i(\theta+\pi)}-\frac{1}{s-i(\theta-\pi)}\\
&\qquad -\frac{1}{s+i(\theta-\pi)}+\frac{1}{s+i(\theta+\pi)}\Bigg)\mathrm{d}s\Bigg)
\end{split}
\end{equation*}
and the right-hand side of \eqref{SS} becomes
\begin{equation*}
\begin{split}
\sum_{j\in\mathbb{Z}}&e^{it\widetilde{H}_{\alpha,B_0,\sigma}}(r_1,\theta_1-2j\sigma\pi;r_2,\theta_2)
=\frac{iB_0}{8\pi\sin(tB_0)}e^{\frac{B_0(r_1^2+r_2^2)}{4i\tan(tB_0)}}\\
&\qquad \times\sum_{j\in\mathbb{Z}}e^{-i\alpha(\theta+2j\sigma\pi)}\Bigg(e^{z\cos(\theta+2j\sigma\pi)}1_{[0,\pi]}(|\theta+2j\sigma\pi|)\\
&\qquad +\frac{1}{2\pi}\int_0^\infty e^{-z\cosh s}\Bigg(\frac{1}{\theta+2j\sigma\pi+\pi+is}-\frac{1}{\theta+2j\sigma\pi-\pi+is}+\\
&\qquad\qquad +\frac{1}{\theta+2j\sigma\pi-\pi-is}-\frac{1}{\theta+2j\sigma\pi+\pi-is}\Bigg)\mathrm{d}s\Bigg).
\end{split}
\end{equation*}
Obviously, the first term in the big bracket of the above expression can be rewritten as
\begin{align*}
\sum_{j\in\mathbb{Z}}&e^{-i\alpha(\theta+2j\sigma\pi)}e^{z\cos(\theta+2j\sigma\pi)}1_{[0,\pi]}(|\theta+2j\sigma\pi|)\\
&\qquad\qquad =\sum_{j\in\mathbb{Z}:|\theta+2j\sigma\pi|\leq\pi}e^{z\cos(\theta+2j\sigma\pi)-i\alpha(\theta+2j\sigma\pi)}.
\end{align*}
For the second term in the big bracket, we apply the already verified formula \footnote{The reason why we only consider the case $\sigma\geq1$ is that this formula is valid only for $\alpha\sigma\in(0,1)$. The verification of this formula relies on the standard formula
\begin{equation*}
\sum_{j\in\mathbb{Z}}\frac{e^{ijx}}{z+j}=\pi e^{-ixz}(\cot(z\pi)-i),\quad x\in(0,2\pi),\quad z\in\mathbb{C}\setminus\mathbb{Z}.
\end{equation*}}
\begin{equation}\label{identity:sum}
\sum_{j\in\mathbb{Z}}\frac{e^{2ij\alpha\sigma\pi}}{\gamma-2j\sigma\pi}=\frac{ie^{i\alpha\gamma}}{\sigma(e^{i\frac{\gamma}{\sigma}}-1)},\quad \alpha\in(0,1/\sigma)\subseteq(0,1),\quad\gamma\in\mathbb{C}\setminus 2\sigma\pi\mathbb{Z},
\end{equation}
to obtain
\begin{equation*}
\begin{split}
\sum_{j\in\mathbb{Z}}e^{-2ij\alpha\sigma\pi}
&\Bigg(\frac{1}{\theta+2j\sigma\pi+\pi+is}-\frac{1}{\theta+2j\sigma\pi-\pi+is}+\\
&\qquad\qquad +\frac{1}{\theta+2j\sigma\pi-\pi-is}-\frac{1}{\theta+2j\sigma\pi+\pi-is}\Bigg)\\
&=\frac{i}{\sigma}\Bigg(\frac{e^{i\alpha(\theta+\pi+is)}}{e^{\frac{i}{\sigma}(\theta+\pi+is)}-1}
-\frac{e^{i\alpha(\theta-\pi+is)}}{e^{\frac{i}{\sigma}(\theta-\pi+is)}-1}\\
&\qquad+\frac{e^{i\alpha(\theta-\pi-is)}}{e^{\frac{i}{\sigma}(\theta-\pi-is)}-1}
-\frac{e^{i\alpha(\theta+\pi-is)}}{e^{\frac{i}{\sigma}(\theta+\pi-is)}-1}\Bigg)
\end{split}
\end{equation*}
and further get
\begin{equation*}
\begin{split}
\frac{1}{2\pi}\int_0^\infty&e^{-z\cosh s}\sum_{j\in\mathbb{Z}}e^{-i\alpha(\theta+2j\sigma\pi)}\Bigg(\frac{1}{\theta+2j\sigma\pi+\pi+is}\\
&-\frac{1}{\theta+2j\sigma\pi-\pi+is}+\frac{1}{\theta+2j\sigma\pi-\pi-is}-\frac{1}{\theta+2j\sigma\pi+\pi-is}\Bigg)\mathrm{d}s\\
&=\frac{i}{2\sigma\pi}\int_0^\infty e^{-z\cosh s}\Bigg(\frac{e^{i\alpha(is+\pi)}}{e^{\frac{i}{\sigma}(is+\theta+\pi)}-1}-\frac{e^{i\alpha(is-\pi)}}{e^{\frac{i}{\sigma}(is+\theta-\pi)}-1}\\
&\qquad\qquad+\frac{e^{-i\alpha(is+\pi)}}{e^{\frac{i}{\sigma}(\theta-\pi-is)}-1}-\frac{e^{-i\alpha(is-\pi)}}{e^{\frac{i}{\sigma}(\theta+\pi-is)}-1}\Bigg)\mathrm{d}s\\
&=\frac{i}{2\sigma\pi}\int_0^\infty e^{-z\cosh s}\Bigg(\frac{e^{-s\alpha-i\alpha\pi}}{1-e^{-\frac{1}{\sigma}(s-i(\theta-\pi))}}-\frac{e^{-s\alpha+i\alpha\pi}}{1-e^{-\frac{1}{\sigma}(s-i(\theta+\pi))}}\\
&\qquad\qquad+\frac{e^{s\alpha+i\alpha\pi}}{1-e^{\frac{1}{\sigma}(s+i(\theta+\pi))}}-\frac{e^{s\alpha-i\alpha\pi}}{1-e^{\frac{1}{\sigma}(s+i(\theta-\pi))}}\Bigg)\mathrm{d}s.
\end{split}
\end{equation*}
Taking into consideration $z=\frac{iB_0r_1r_2}{2\sin(tB_0)}$ and $\theta=tB_0-(\theta_1-\theta_2)$, we finally obtain
\begin{equation}\label{SS:right}
\begin{split}
\sum_{j\in\mathbb{Z}}&e^{it\widetilde{H}_{\alpha,B_0,\sigma}}(r_1,\theta_1-2j\sigma\pi;r_2,\theta_2)=\frac{iB_0e^{itB_0\alpha}}{8\pi\sin(tB_0)}
e^{\frac{B_0(r_1^2+r_2^2)}{4i\tan(tB_0)}}\\
&\times\Bigg[\sum_{j\in\mathbb{Z}:|tB_0-(\theta_1-\theta_2)+2j\sigma\pi|\leq\pi}
e^{\frac{B_0r_1r_2}{2i\sin(tB_0)}\cos(tB_0-(\theta_1-\theta_2)+2j\sigma\pi)-i\alpha(tB_0-(\theta_1-\theta_2)+2j\sigma\pi)}\\
&\qquad\qquad -\frac{1}{2\sigma\pi}\int_0^\infty e^{\frac{B_0r_1r_2}{2i\sin(tB_0)}\cosh s}\Bigg(\frac{e^{-s\alpha-i\alpha\pi}}{1-e^{-\frac{1}{\sigma}(s-i(\theta-\pi))}}-\frac{e^{-s\alpha+i\alpha\pi}}{1-e^{-\frac{1}{\sigma}(s-i(\theta+\pi))}}\\
&\qquad\qquad\qquad\qquad+\frac{e^{s\alpha+i\alpha\pi}}{1-e^{\frac{1}{\sigma}(s+i(\theta+\pi))}}
-\frac{e^{s\alpha-i\alpha\pi}}{1-e^{\frac{1}{\sigma}(s+i(\theta-\pi))}}\Bigg)\mathrm{d}s\Bigg],
\end{split}
\end{equation}
which, in view of the function $A_{\sigma,\alpha}(s,tB_0-\theta)$ in \eqref{A-s-1-2}, is exactly the same as the expression \eqref{S:express} obtained via the first approach.

Therefore, the connection formula \eqref{SS} is verified and meanwhile we obtain the desire representation formula \eqref{S:express} for the Schr\"odinger kernel $K_t^S(r_1,\theta_1;r_2,\theta_2)$.
\begin{remark}
One may find that, after taking into consideration the function $A_{\sigma,\alpha}(s,tB_0-\theta)$ in \eqref{A-s-1-2}, the expression \eqref{SS:right} does not perfectly coincide with the expression \eqref{S:express} calculated via the first approach. In fact, if one carefully compares the terms in the big bracket of \eqref{SS:right} with the calculations for the second term \eqref{equ:term2} in the first approach, one will realize that the connection formula \eqref{SS} is indeed verified. For our purpose of proving the remained dispersive estimate \eqref{dis:S}, we shall use the representation formula \eqref{S:express} involving $A_{\sigma,\alpha}(s,tB_0-\theta)$ in the following.
\end{remark}
With the explicit representation formula \eqref{S:express} for the kernel of the Schr\"odinger propagator $e^{itH_{\alpha,B_0,\sigma}}$ in hand, we can now prove the Schr\"odinger dispersive estimate
\begin{equation}\label{S-dispersive}
\left\|e^{itH_{\alpha,B_0,\sigma}}f\right\|_{L^\infty(X)}\lesssim_{\sigma}\frac{B_0}{|\sin(tB_0)|}\|f\|_{L^1(X)},\quad t\neq\frac{\pi}{B_0}\mathbb{Z}.
\end{equation}
Indeed, the proof of \eqref{S-dispersive} is finished if we could show the following uniform bounds
\begin{equation}\label{S-dis-1}
\left|\sum_{j\in\mathbb{Z}:|\theta+2j\sigma\pi|\leq\pi}e^{\frac{iB_0r_1r_2}{2\sin(tB_0)}\cos(\theta+2j\sigma\pi)-i\alpha(\theta+2j\sigma\pi)}\right|
\lesssim_\sigma1
\end{equation}
and
\begin{equation}\label{S-dis-2}
\left|\frac{1}{2\sigma\pi}\int_0^\infty e^{\frac{B_0r_1r_2}{2i\sin(tB_0)}\cosh s}A_{\sigma,\alpha}(s,\theta)\mathrm{d}s\right|\lesssim_\sigma1,
\end{equation}
where $A_{\sigma,\alpha}(s,\theta)$ is given by \eqref{A-s-1-2}.

For the first inequality \eqref{S-dis-1}, it is not hard to see that, for a fixed $\sigma\geq1$, there are only finite many integers $j\in\mathbb{Z}$ satisfying $|\theta+2j\sigma\pi|\leq\pi$. More precisely, we have
\begin{equation*}
\begin{split}
&\left|\sum_{j\in\mathbb{Z}:|\theta+2j\sigma\pi|\leq\pi}e^{\frac{iB_0r_1r_2}{2\sin(tB_0)}\cos(\theta+2j\sigma\pi)-i\alpha(\theta+2j\sigma\pi)}\right| \\
&\qquad\qquad\qquad\leq|\#\{j\in\mathbb{Z}:|\theta+2j\sigma\pi|\leq\pi\}|\lesssim1+\sigma^{-1},
\end{split}
\end{equation*}
which yields the expected inequality \eqref{S-dis-1}. For the second inequality \eqref{S-dis-2}, it is sufficient to verify
\begin{equation}\label{S-dis-3}
\int_0^\infty|A_{\sigma,\alpha}(s,\theta)|\mathrm{d}s\lesssim_\sigma1,
\end{equation}
which, in view of \eqref{A-s-1-2}, follows further from
\begin{align}
&\int_0^{\infty}|e^{-s\alpha}\sin(s\alpha)|\mathrm{d}s\lesssim1,\label{ineq:1}\\
&\int_0^\infty\left|\frac{e^{\pm s\alpha}\sin\phi_\pm}{\cosh(s/\sigma)-\cos\phi_\pm}\right|\mathrm{d}s\lesssim_\sigma1,\label{ineq:2}\\
&\int_0^\infty\left|\frac{e^{\pm s\alpha}(\cos\phi_\pm-e^{-s/\sigma})}{\cosh(s/\sigma)-\cos\phi_\pm}\right|\mathrm{d}s\lesssim_\sigma1\label{ineq:3}
\end{align}
with $\phi_\pm$ being given in \eqref{+-}. The inequality \eqref{ineq:1} is trivial. As for \eqref{ineq:2} and \eqref{ineq:3}, by the assumption $\alpha\in(0,1/\sigma)$ (so that $\frac{1}{\sigma}\pm\alpha>0$) and the following identity
\begin{equation*}
\cosh\beta-\cos\gamma=\sinh^2(\beta/2)+\sin^2(\gamma/2),\quad \forall\beta,\gamma\in\mathbb{R},
\end{equation*}
we obtain
\begin{align*}
\int_0^\infty&\left|\frac{e^{\pm s\alpha}\sin\phi_\pm}{\cosh\frac{s}{\sigma}-\cos\phi_\pm}\right|\mathrm{d}s\\
&\lesssim\int_0^1\left|\frac{\sin\phi_\pm}{\sinh^2\frac{s}{2\sigma}+\sin^2\frac{\phi_\pm}{2}}\right|\mathrm{d}s
+\int_1^\infty\left|\frac{e^{\pm s\alpha}}{\sinh^2\frac{s}{2\sigma}+\sin^2\frac{\phi_\pm}{2}}\right|\mathrm{d}s\\
&\lesssim\int_0^1\left|\frac{\sin\phi_\pm}{(\frac{s}{2\sigma})^2+\sin^2\frac{\phi_\pm}{2}}\right|\mathrm{d}s+\int_1^\infty e^{-(\frac{1}{\sigma}\pm\alpha)s}\mathrm{d}s\\
&\lesssim_\sigma1
\end{align*}
and
\begin{align*}
\int_0^\infty&\left|\frac{e^{\pm s\alpha}(\cos\phi_\pm-e^{-s/\sigma})}{\cosh(s/\sigma)-\cos\phi_\pm}\right|\mathrm{d}s\\
&\lesssim\int_0^1\left|\frac{e^{\pm s\alpha}(\cos\phi_\pm-e^{-s/\sigma})}{\sinh^2\frac{s}{2\sigma}+\sin^2\frac{\phi_\pm}{2}}\right|\mathrm{d}s
+\int_1^\infty\left|\frac{e^{\pm s\alpha}(\cos\phi_\pm-e^{-s/\sigma})}{\sinh^2\frac{s}{2\sigma}+\sin^2\frac{\phi_\pm}{2}}\right|\mathrm{d}s\\
&\lesssim\int_0^1\left|\frac{\frac{s}{\sigma}+\sin^2\frac{\phi_\pm}{2}}{(\frac{s}{\sigma})^2+\sin^2\frac{\phi_\pm}{2}}\right|\mathrm{d}s+\int_1^\infty e^{-(\frac{1}{\sigma}\pm\alpha)s}\mathrm{d}s\\
&\lesssim_\sigma1,
\end{align*}
which yield \eqref{ineq:2} and \eqref{ineq:3}, respectively.
Hence, the second inequality \eqref{S-dis-2} follows and so does the desired Schr\"odinger dispersive estimate \eqref{S-dispersive}.

\subsection{The third approach}

Recall that the Schr\"odinger dispersive estimate \eqref{S-dispersive} is reduced to the uniform boundedness of the following series
\begin{equation}\label{ker:sch}
\sum_{k\in\mathbb{Z}}e^{i\frac{k}{\sigma}(tB_0-(\theta_1-\theta_2))}I_{|k/\sigma+\alpha|}\left(\frac{iB_0r_1r_2}{2\sin(tB_0)}\right).
\end{equation}
Define a function $K(\rho,\delta):\mathbb{R}_+\times[a,b]\rightarrow\mathbb{C}$ for any finite real numbers $a<b$
\begin{equation}\label{red:ker}
K(\rho,\delta):=\sum_{k\in\mathbb{Z}}e^{i\frac{k}{\sigma}\delta}I_{|k/\sigma+\alpha|}(i\rho),\quad\alpha\in(0,1/\sigma),\quad \sigma\geq1,
\end{equation}
then one observes that the series \eqref{ker:sch} can be expressed in terms of \eqref{red:ker} as
\begin{equation}\label{Sch:represent}
\begin{split}
 \sum_{k\in\mathbb{Z}}e^{i\frac{k}{\sigma}(tB_0-(\theta_1-\theta_2))}&I_{|k/\sigma+\alpha|}\left(\frac{iB_0r_1r_2}{2\sin(tB_0)}\right)\\
&=K\left(\frac{B_0r_1r_2}{2\sin(tB_0)},tB_0-(\theta_1-\theta_2)\right).
\end{split}
\end{equation}
The Schr\"odinger dispersive estimate \eqref{S-dispersive} now follows from the uniform boundedness result of a recent preprint \cite{GWXZ25}.
\begin{lemma}\label{lem:kernel}
Let $K(\rho,\delta)$ be defined by \eqref{red:ker}, then it holds for any finite $R<+\infty$
\begin{equation*}
  \sup_{\rho\geq0\atop |\delta|\leq R}|K(\rho,\delta)|<+\infty.
\end{equation*}
\end{lemma}
Lemma \ref{lem:kernel} is actually summarized from \cite[Section 2]{GWXZ25} and it is valid only for the range $\sigma\geq1$. We point out that Lemma \ref{lem:kernel} can be obtained alternatively by the perturbation argument in \cite[Section 3]{FFFP15}; we shall omit the details to avoid repeating.

\subsection{Strichartz estimate for Schr\"odinger equation}

The proof of Strichartz estimates for the Schr\"odinger equation \eqref{eq:S} is quite standard by using the dispersive estimate \eqref{S-dispersive} and the abstract mechanism of Keel-Tao \cite{KT98} and we give a quick sketch here for completeness.
\begin{proposition}\label{prop:KT-Sch}
Let $(X,\mathrm{d}\mu)$ be a measurable space and $H$ a Hilbert space. Suppose that for every $t\in\mathbb{R},U(t):H\rightarrow L^2(X)$ satisfies the energy estimate
\begin{equation*}
  \|U(t)\|_{H\rightarrow L^2}\leq C,\quad t\in\mathbb{R}
\end{equation*}
and that for some $\sigma>0$
\begin{equation*}
  \|U(t)U(s)^\ast f\|_{L^\infty(X)}\leq C|t-s|^{-\sigma}\|f\|_{L^1(X)},\quad t\neq s.
\end{equation*}
Then it holds
\begin{equation*}
  \|U(t)f\|_{L_t^q(I;L^p(X))}\lesssim\|f\|_{L^2(X)},
\end{equation*}
where $(q,p)$ is the sharp $\sigma$-admissible pair, i.e. satisfies
\begin{equation*}
(q,p)\in\Lambda_S:=\left\{(q,p)\in[2,+\infty]\times[2,+\infty):\frac{1}{q}=\sigma\left(\frac{1}{2}-\frac{1}{p}\right),(q,p,\sigma)\neq(2,\infty,1)\right\}.
\end{equation*}
\end{proposition}
By the spectral theorem, we trivially have the energy estimate
\begin{equation*}
  \|e^{itH_{\alpha,B_0,\sigma}}\|_{L^2(X)\rightarrow L^2(X)}\lesssim1.
\end{equation*}
Now the dispersive estimate \eqref{S-dispersive}, together with Proposition \ref{prop:KT-Sch} applied to $U(t)=e^{itH_{\alpha,B_0,\sigma}}$, yields the expected Strichartz estimate \eqref{str:S} in Corollary \ref{cor:S-W}.

\section{Littlewood-Paley thery for $H_{\alpha,B_0,\sigma}$}\label{sec:heat}

This section is devoted to establish the Littlewood-Paley theory through the Gaussian boundedness of the heat kernel associated to the operator $H_{\alpha,B_0,\sigma}$ on the product cone $(X,g)$. The key ingredient is an explicit representation formula for the heat kernel via Schulman's ansatz (see e.g. \cite{Sto89,Sto17}).

\subsection{The heat kerenl and its Gaussian boundedness}

In the following proposition, we obtain a series representation for the heat kernel associated to the operator $H_{\alpha,B_0,\sigma}$.
\begin{proposition}\label{prop:H}
Let $H_{\alpha,B_0,\sigma}$ denote the Friedrichs extension of the operator \eqref{op} and let $K_t^H(r_1,\theta_1;r_2,\theta_2)$ be the kernel of the heat semigroup $e^{-tH_{\alpha,B_0,\sigma}}$, then we have
\begin{equation}\label{rep:h}
\begin{split}
K_t^H(r_1,\theta_1;r_2,\theta_2)=&\frac{B_0e^{-tB_0\alpha}}{4\pi\sigma\sinh(tB_0)}e^{-\frac{B_0(r_1^2+r_2^2)}{4\tanh(tB_0)}}\\
    &\qquad \times\sum_{k\in\mathbb{Z}}e^{i\frac{k}{\sigma}(\theta_1-\theta_2+itB_0)}I_{\alpha_k}\left(\frac{B_0r_1r_2}{2\sinh tB_0}\right),\quad t>0,
\end{split}
\end{equation}
where $\alpha_k=|k/\sigma+\alpha|$ and $I_\nu(r)$ is the modified Bessel function of the first kind of order $\nu$ with respect to the variable $r$.
\end{proposition}
\begin{proof}
The representation formula \eqref{rep:h} for the heat kernel can be obtained similar to that for the Schr\"odinger kernel. Indeed, since the family of the $L^2$-normalised eigenfunctions $\Big\{\tilde{V}_{k,m}(r,\theta)\Big\}_{k\in\mathbb{Z}, m\in\mathbb{N}}$ forms a complete orthonormal basis for $L^2(X)$ (see Proposition \ref{prop:spect}), a function $f\in L^2(X)$ can be expanded in terms of this basis as
\begin{equation*}
f(r,\theta)=\sum_{k\in\mathbb{Z} \atop m\in\mathbb{N}}c_{k,m}\tilde{V}_{k,m}(r,\theta),
\end{equation*}
where
\begin{equation}\label{cmk1}
c_{k,m}=\int_0^\infty\int_0^{2\sigma\pi}f(r,\theta)\overline{\tilde{V}_{k,m}(r,\theta)}\mathrm{d}\theta r\mathrm{d}r.
\end{equation}
The heat semigroup $e^{-tH_{\alpha,B_0,\sigma}}$ for $t>0$ can be represented via
\begin{equation*}
\begin{split}
u(t,r_1,\theta_1):&=\left(e^{-tH_{\alpha,B_0,\sigma}}f\right)(r_1,\theta_1)\\
&=\sum_{k\in\mathbb{Z} \atop m\in\mathbb{N}}e^{-t\lambda_{k,m}}\left(\int_0^\infty\int_0^{2\sigma\pi}f(r_2,\theta_2)\overline{\tilde{V}_{k,m}(r_2,\theta_2)}\mathrm{d}\theta_2 r_2\mathrm{d}r_2\right)\tilde{V}_{k,m}(r_1,\theta_1).
\end{split}
\end{equation*}
A function $f\in L^2(X)$ can be alternatively expressed as
\begin{equation*}
f(r_2,\theta_2)=\sum_{k\in\mathbb{Z}}f_k(r_2)e^{i\frac{k}{\sigma}\theta_2},
\end{equation*}
where
\begin{equation}\label{f-k}
f_k(r_2)=\frac{1}{2\sigma\pi}\int_0^{2\sigma\pi}f(r_2,\theta_2)e^{-i\frac{k}{\sigma}\theta_2}\mathrm{d}\theta_2.
\end{equation}
In view of the eigenvalue \eqref{eigen-v}, eigenfunction \eqref{eigen-f} and the normalized constant \eqref{V-km-2}, we obtain
\begin{align*}
u(t,r_1,\theta_1)&=\sum_{k\in\mathbb{Z} \atop m\in\mathbb{N}}e^{-t\lambda_{k,m}}\frac{V_{k,m}(x)}{\|V_{k,m}\|^2_{L^2}}\left(\int_0^\infty f_k(r_2)e^{-\frac{B_0r_2^2}{4}}P_{k,m}\left(\frac{B_0r_2^2}{2}\right)r_2^{1+\alpha_k}\mathrm{d}r_2\right)\\
&=\frac{B_0}{2}\sum_{k\in\mathbb{Z}}e^{i\frac{k}{\sigma}\theta_1}\frac{B_0^{\alpha_k}e^{-t\beta_k}}{2^{\alpha_k}\Gamma(1+\alpha_k)}\Bigg[\sum_{m=0}^\infty
\Bigg(
  \begin{array}{c}
    m+\alpha_k \\
    m \\
  \end{array}
\Bigg)e^{-2tmB_0}\\
&\times\Bigg(\int_0^\infty f_k(r_2)(r_1r_2)^{\alpha_k}e^{-\frac{B_0(r_1^2+r_2^2)}{4}}P_{k,m}\left(\frac{B_0r_2^2}{2}\right)P_{k,m}\left(\frac{B_0r_1^2}{2}\right)r_2 \mathrm{d}r_2\Bigg)\Bigg],
\end{align*}
where $\alpha_k$ and $\beta_k$ are as in \eqref{beta-k}. Taking into consideration the integral form \eqref{P-km} of $P_{k,m}$, the summation formula \eqref{J-alpha} and the product rule \eqref{formula} of the Bessel functions $J_{\alpha_k}$, we will get
\begin{align*}
u(t,r_1,\theta_1)&=\frac{B_0}{4\sigma\pi}\int_0^\infty\int_0^{2\sigma\pi}\sum_{k\in\mathbb{Z}}e^{i\frac{k}{\sigma}(\theta_1-\theta_2)}
\frac{B_0^{\alpha_k}e^{-t\beta_k}}{2^{\alpha_k}\Gamma(1+\alpha_k)}(r_1r_2)^{\alpha_k}\\
&\times e^{-\frac{B_0(r_1^2+r_2^2)}{4}}f(r_2,\theta_2)\frac{2^{\alpha_k}e^{tB_0\alpha_k}}{(B_0r_1r_2)^{\alpha_k}}
e^{-\frac{B_0(r_1^2+r_2^2)}{2(e^{2tB_0}-1)}}I_{\alpha_k}
\left(\frac{B_0r_1r_2}{2\sinh(tB_0)}\right)r_2\mathrm{d}r_2\mathrm{d}\theta_2\\
&=\frac{B_0e^{-tB_0\alpha}}{4\pi\sigma\sinh(tB_0)}\int_0^\infty\int_0^{2\sigma\pi}e^{-\frac{B_0(r_1^2+r_2^2)}{4\tanh(tB_0)}}f(r_2,\theta_2)\\
&\qquad\qquad \times\sum_{k\in\mathbb{Z}}e^{i\frac{k}{\sigma}(\theta_1-\theta_2+itB_0)}I_{\alpha_k}\left(\frac{B_0r_1r_2}{2\sinh(tB_0)}\right)r_2 \mathrm{d}r_2\mathrm{d}\theta_2
\end{align*}
and the heat kernel $K_t^H(r_1,\theta_1;r_2,\theta_2)$ is thus given by
\begin{equation*}
\begin{split}
K_t^H(r_1,\theta_1;r_2,\theta_2)&=\frac{B_0e^{-tB_0\alpha}}{4\pi\sigma\sinh(tB_0)}e^{-\frac{B_0(r_1^2+r_2^2)}{4\tanh(tB_0)}}\\
&\qquad\qquad\times\sum_{k\in\mathbb{Z}}e^{i\frac{k}{\sigma}(\theta_1-\theta_2+itB_0)}I_{\alpha_k}\left(\frac{B_0r_1r_2}{2\sinh(tB_0)}\right).
\end{split}
\end{equation*}
\end{proof}
\begin{remark}
The formula \eqref{rep:h} specialized to $\sigma=1$ corresponds to the Euclidean case in \cite[(3.2)]{WZZ25-JDE}
\begin{equation*}
\begin{split}
K_t^H(r_1,\theta_1;r_2,\theta_2)&=\frac{B_0e^{-tB_0\alpha}}{4\pi\sinh(tB_0)}e^{-\frac{B_0(r_1^2+r_2^2)}{4\tanh(tB_0)}}\\
    &\qquad\qquad\times\sum_{k\in\mathbb{Z}}e^{ik(\theta_1-\theta_2+itB_0)}I_{|k+\alpha|}\left(\frac{B_0r_1r_2}{2\sinh(tB_0)}\right),
\end{split}
\end{equation*}
which, if further specialized to $\alpha=0$, reduces to the heat kernel for the Landau Hamiltonian $(-i\nabla+\frac{B_0}{2}(-x_2,x_1))^2$
\begin{equation*}
\begin{split}
K_t^H(r_1,\theta_1;r_2,\theta_2)
&=\frac{B_0}{4\pi\sinh(tB_0)}e^{-\frac{B_0(r_1^2+r_2^2)}{4\tanh(tB_0)}}\\
&\qquad\qquad\qquad\times \sum_{k\in\mathbb{Z}}e^{ik(\theta_1-\theta_2+itB_0)}I_{|k|}\left(\frac{B_0r_1r_2}{2\sinh(tB_0)}\right)\\
&=\frac{(4\pi)^{-1}B_0}{\sinh(tB_0)}e^{-\frac{B_0(r_1^2+r_2^2)}{4\tanh(tB_0)}+\frac{B_0r_1r_2}{2\sinh(tB_0)}\cosh(tB_0+i(\theta_1-\theta_2))}.
\end{split}
\end{equation*}
\end{remark}
Next, we shall apply again Schulman's universal covering space technique to construct an explicit representation formula for the heat kernel $K_t^H(r_1,\theta_1;r_2,\theta_2)$ associated with the operator $H_{\alpha,B_0,\sigma}$ and then obtain the Gaussian upper bound for this kernel. Recall that the universal covering of $X=(0,+\infty)_r\times\mathbb{S}_\sigma^1$ is given by $\tilde{X}:=(0,+\infty)_r\times\mathbb{R}$ with the structural group $\Gamma=2\sigma\pi\mathbb{Z}$ acting on the second entry of the product cone $X$, then one has $X=\tilde{X}/\Gamma$ and the following connection formula (see e.g. \cite[(1)]{Sto89} and \cite[(3.27)]{WZZ25-JDE})
\begin{equation}\label{HH}
K_t^H(r_1,\theta_1;r_2,\theta_2)=\sum_{j\in\mathbb{Z}}\widetilde{K}_t^H(r_1,\theta_1+2j\sigma\pi;r_2,\theta_2),
\end{equation}
bridging the expected heat kernel $K_t^H(r_1,\theta_1;r_2,\theta_2)$ on $X$ with the heat kernel $\widetilde{K}_t^H(r_1,\theta_1;r_2,\theta_2)$ on $\tilde{X}$. It should be clarified here that $K_t^H(r_1,\theta_1;r_2,\theta_2)$ denotes the heat kernel associated with $H_{\alpha,B_0,\sigma}$ defined on $X$ and $\widetilde{K}_t^H(r_1,\theta_1;r_2,\theta_2)$ denotes the heat kernel associated with $\tilde{H}_{\alpha,B_0,\sigma}$ defined on $\tilde{X}$ (see \eqref{t-H} below).
\begin{proposition}\label{prop:H'}
Let $K_t^H(r_1,\theta_1;r_2,\theta_2)$ be the heat kernel associated to the operator $H_{\alpha,B_0,\sigma}$, then
\begin{equation}\label{heatkernel:2}
\begin{split}
K_t^H&(r_1,\theta_1;r_2,\theta_2)=\frac{B_0}{4\pi\sigma\sinh(tB_0)}e^{-\frac{B_0(r_1^2+r_2^2)}{4\tanh(tB_0)}} \\
&\qquad\times\Bigg(e^{-tB_0\alpha}\sum_{j\in\mathbb{Z}:|\theta+2j\sigma\pi|\leq\pi}e^{\frac{B_0r_1r_2}{2\sinh(tB_0)}\cosh(tB_0+i(\theta+2j\sigma\pi))-i\alpha(\theta+2j\sigma\pi)}\\
&\qquad\qquad\qquad+\frac{1}{2i\sigma\pi}\int_{-\infty}^{\infty}e^{-\frac{B_0r_1r_2}{2\sinh(tB_0)}\cosh s}B_{\sigma,\alpha}(s,\theta)\mathrm{d}s\Bigg),
\end{split}
\end{equation}
where $\theta:=\theta_1-\theta_2$ and $B_{\sigma,\alpha}(s,\theta)$ is given by
\begin{equation}\label{varphi}
B_{\sigma,\alpha}(s,\theta):=e^{s\alpha}\left(\frac{e^{i\alpha\pi}}{e^{\frac{1}{\sigma}(s+tB_0+i(\theta+\pi))}-1}-\frac{e^{-i\alpha\pi}}{e^{\frac{1}{\sigma}(s+tB_0+i(\theta-\pi))}-1}\right).
\end{equation}
\end{proposition}
\begin{proof}
The main idea is again to replace the summation over $k\in\mathbb{Z}$ of the series \eqref{rep:h} by the integration over $\tilde{k}\in\mathbb{R}$ so that the translation invariance and dilation invariance of integration over the whole line can be applied. To this end, let us consider the operator defined on $\tilde{X}:=(0,+\infty)_r\times\mathbb{R}_\theta$
\begin{equation}\label{t-H}
\tilde{H}_{\alpha,B_0,\sigma}=-\partial_{rr}-\frac{1}{r}\partial_r+\frac{1}{r^2}\left(-i\partial_\theta+\alpha+\frac{B_0r^2}{2}\right)^2,
\end{equation}
which acts on functions of $L^2(\tilde{X},r\mathrm{d}r\mathrm{d}\theta)$. Taking $\psi(\theta)=\frac{1}{\sqrt{2\sigma\pi}}e^{i\left(\frac{\tilde{k}}{\sigma}-\alpha\right)\theta}$ as the eigenfunction of the operator $\Big(-i\partial_\theta+\alpha+\frac{B_0r^2}{2}\Big)^2$ on $L^2(\mathbb{R}_\theta)$, we have
\begin{equation}\label{eq:ef}
\left(-i\partial_\theta+\alpha+\frac{B_0r^2}{2}\right)^2\psi(\theta)=\left(\frac{\tilde{k}}{\sigma}+\frac{B_0r^2}{2}\right)^2\psi(\theta),\quad \theta,\tilde{k}\in\mathbb{R}.
\end{equation}
Following the proof of Proposition \ref{prop:H}, we can similarly obtain
\begin{equation}\label{K-tilde}
\begin{split}
\widetilde{K}_t^H(r_1,\theta_1;r_2,\theta_2)=&\frac{B_0}{4\pi\sigma\sinh(tB_0)}e^{-\frac{B_0(r_1^2+r_2^2)}{4\tanh(tB_0)}}\\
&\times\int_{\mathbb{R}}e^{i(\frac{\tilde{k}}{\sigma}-\alpha)(\theta_1-\theta_2-itB_0)}I_{|\tilde{k}/\sigma|}\left(\frac{B_0r_1r_2}{2\sinh(tB_0)}\right)\mathrm{d}\tilde{k},
\end{split}
\end{equation}
where $(r_j,\theta_j)\in\tilde{X},j=1,2$ and $I_\nu(t)$ denotes the modified Bessel function of the first kind of order $\nu$ with respect to the variable $t$.

A key formula in order to calculate the integral in \eqref{K-tilde} comes from \cite[(2.11)]{Sto17}.
\begin{lemma}\label{lem:key2}
For any $z\in\mathbb{C}$, it holds
\begin{equation}\label{id:key2}
\begin{split}
\int_{\mathbb{R}}e^{z\tilde{k}}I_{|\tilde{k}|}(x)\mathrm{d}\tilde{k}=&e^{x\cosh z}H(\pi-|\Im z|)\\
&+\frac{1}{2\pi i}\int_{\mathbb{R}}e^{-x\cosh s}\left(\frac{1}{s+z+i\pi}-\frac{1}{s+z-i\pi}\right)\mathrm{d}s,
\end{split}
\end{equation}
where $H(x)$ denotes the Heaviside step function
\begin{equation*}
H(x)=
\begin{cases}
1,\quad x>0,\\
0,\quad x\leq 0.
\end{cases}
\end{equation*}
\end{lemma}
We refer to \cite[Pages 24-27]{WZZ25-JDE} for a proof of Lemma \ref{lem:key2} and continue to calculate the integral in \eqref{K-tilde}.
Setting $z=tB_0+i(\theta_1-\theta_2)$ for the moment, it follows by Lemma \ref{lem:key2} that
\begin{equation*}
\begin{split}
&\widetilde{K}_t^H(r_1,\theta_1;r_2,\theta_2)=\frac{B_0e^{-z\alpha}}{4\pi\sigma\sinh(tB_0)}e^{-\frac{B_0(r_1^2+r_2^2)}{4\tanh(tB_0)}}\\
&\qquad \qquad\times\Bigg(e^{\frac{B_0r_1r_2}{2\sinh(tB_0)}\cosh z}H\left(\pi-|\theta_1-\theta_2|\right)\\
&\quad+\frac{1}{2\pi i}\int_{\mathbb{R}}e^{-\frac{B_0r_1r_2}{2\sinh(tB_0)}\cosh s}\left(\frac{1}{s+z+i\pi}-\frac{1}{s+z-i\pi}\right)\mathrm{d}s\Bigg).
\end{split}
\end{equation*}
In view of the formula \eqref{HH}, we obtain
\begin{equation}\label{HH'}
\begin{split}
&K_t^H(r_1,\theta_1;r_2,\theta_2)=\sum_{j\in\mathbb{Z}}\widetilde{K}_t^H(r_1,\theta_1+2j\sigma\pi;r_2,\theta_2)\\
&=\frac{B_0}{4\pi\sigma\sinh(tB_0)}e^{-\frac{B_0(r_1^2+r_2^2)}{4\tanh(tB_0)}}\sum_{j\in\mathbb{Z}}e^{-z_j\alpha}\\
&\qquad\qquad\qquad\times\Bigg(e^{\frac{B_0r_1r_2}{2\sinh(tB_0)}\cosh z_j}H\left(\pi-|\theta_1+2j\sigma\pi-\theta_2|\right)\\
&\qquad\qquad+\frac{1}{2\pi i}\int_{\mathbb{R}}e^{-\frac{B_0r_1r_2}{2\sinh(tB_0)}\cosh s}\left(\frac1{s+z_j+i\pi}-\frac1{s+z_j-i\pi}\right)\mathrm{d}s\Bigg),
\end{split}
\end{equation}
where $z_j=tB_0+i(\theta_1+2j\sigma\pi-\theta_2)$.

The first term in the big bracket of \eqref{HH'} can be rewritten as
\begin{equation*}
\begin{split}
&\sum_{j\in\mathbb{Z}}e^{-z_j\alpha}e^{\frac{B_0r_1r_2}{2\sinh(tB_0)}\cosh z_j}H\left(\pi-|\theta+2j\sigma\pi|\right)\\
&=e^{-tB_0\alpha}\sum_{j\in\mathbb{Z}:|\theta+2j\sigma\pi|\leq\pi}
e^{\frac{B_0r_1r_2}{2\sinh(tB_0)}\cosh\big(tB_0+i(\theta+2j\sigma\pi)\big)-i\alpha(\theta+2j\sigma\pi)},
\end{split}
\end{equation*}
where $\theta:=\theta_1-\theta_2$.
As for the second term in the big bracket of \eqref{HH'}, we apply the summation formula \eqref{identity:sum} again to obtain
\begin{equation*}
\begin{split}
&\sum_{j\in\mathbb{Z}}e^{-z_j\alpha}\left(\frac{1}{s+z_j+i\pi}-\frac{1}{s+z_j-i\pi}\right)\\
&=e^{-\alpha(tB_0+i(\theta_1-\theta_2))}\sum_{j\in\mathbb{Z}}\Bigg(\frac{e^{-2ij\alpha\sigma\pi}}{s+tB_0+i(\theta_1-\theta_2+2j\sigma\pi+\pi)}\\
&\qquad\qquad -\frac{e^{-2ij\alpha\sigma\pi }}{s+tB_0+i(\theta_1-\theta_2+2j\sigma\pi-\pi)}\Bigg)\\
&=-ie^{-\alpha(tB_0+i(\theta_1-\theta_2))}\sum_{j\in\mathbb{Z}}\left(\frac{e^{-2ij\alpha\sigma\pi}}{\gamma_++2j\sigma\pi}-\frac{e^{-2ij\alpha\sigma\pi}}{\gamma_-+2j\sigma\pi}\right)\\
&=e^{-\alpha(tB_0+i(\theta_1-\theta_2))}\left(\frac{e^{i\alpha\gamma_+}}{\sigma(e^{i\gamma_+/\sigma}-1)}-\frac{e^{i\alpha\gamma_-}}{\sigma(e^{i\gamma_-/\sigma}-1)}\right)\\
&=\frac{e^{s\alpha}}{\sigma}\left(\frac{e^{i\alpha\pi}}{e^{\frac{1}{\sigma}(s+tB_0+i(\theta_1-\theta_2+\pi))}-1}-\frac{e^{-i\alpha\pi}}{e^{\frac{1}{\sigma}(s+tB_0+i(\theta_1-\theta_2-\pi))}-1}\right),
\end{split}
\end{equation*}
where $i\gamma_\pm=s+tB_0+i(\theta_1-\theta_2\pm\pi)$.
Therefore, we obtain
\begin{equation*}
\begin{split}
&K_t^H(r_1,\theta_1;r_2,\theta_2)=\frac{B_0}{4\pi\sigma\sinh(tB_0)}e^{-\frac{B_0(r_1^2+r_2^2)}{4\tanh(tB_0)}}\Bigg(e^{-tB_0\alpha}\\
&\times\sum_{j\in\mathbb{Z}:|\theta_1+2j\sigma\pi-\theta_2|\leq\pi}
e^{\frac{B_0r_1r_2}{2\sinh(tB_0)}\cosh(tB_0+i(\theta_1+2j\sigma\pi-\theta_2))-i\alpha(\theta_1-\theta_2+2j\sigma\pi)}\\
&\qquad+\frac{1}{2\sigma\pi i}\int_{\mathbb{R}}e^{-\frac{B_0r_1r_2}{2\sinh(tB_0)}\cosh s}\Bigg(\frac{e^{\alpha(s+i\pi)}}{e^{\frac{1}{\sigma}(s+tB_0+i(\theta_1-\theta_2+\pi))}-1}\\
&\qquad\qquad\qquad -\frac{e^{\alpha(s-i\pi)}}{e^{\frac{1}{\sigma}(s+tB_0+i(\theta_1-\theta_2-\pi))}-1}\Bigg)\mathrm{d}s\Bigg)
\end{split}
\end{equation*}
as desired.
\end{proof}
Now we are in the position to prove the Gaussian boundedness of the heat kernel $K_t^H(r_1,\theta_1;r_2,\theta_2)$ by using the explicit representation formula \eqref{heatkernel:2}.
\begin{proposition}[Gaussian boundedness of the heat kernel]\label{prop:heat-Gaussbound}
Let $K_t^H(r_1,\theta_1;r_2,\theta_2)$ be the heat kernel in \eqref{heatkernel:2}, then it holds
\begin{equation}\label{est:heatkernel:2}
\left|K_t^H(r_1,\theta_1;r_2,\theta_2)\right|\lesssim|\sinh(tB_0)|^{-1}e^{-\frac{B_0(r_1^2+r_2^2)}{4\tanh(tB_0)}},\quad t>0.
\end{equation}
\end{proposition}
\begin{proof}
Since there are only finite many integers $j\in\mathbb{Z}$ (depending on the given parameter $\sigma\geq1$) satisfying $|\theta+2j\sigma\pi|\leq\pi$ for any $\theta\in[-4\sigma\pi,4\sigma\pi]$ and
\begin{equation*}
\cosh(tB_0+i\theta)=\cos\theta\cosh(tB_0)+i\sin\theta\sinh(tB_0),
\end{equation*}
the first term of the heat kernel \eqref{heatkernel:2} is controlled by
\begin{equation}\label{est:h1}
\begin{split}
\Bigg|\frac{B_0}{4\pi\sigma\sinh(tB_0)}&e^{-\frac{B_0(r_1^2+r_2^2)}{4\tanh(tB_0)}}e^{-tB_0\alpha}\\
\times\sum_{j\in\mathbb{Z}:|\theta+2j\sigma\pi|\leq\pi}&e^{\frac{B_0r_1r_2}{2\sinh(tB_0)}\cosh(tB_0+i(\theta+2j\sigma\pi))-i\alpha(\theta+2j\sigma\pi)}\Bigg|\\
&\qquad\qquad\lesssim_\sigma|\sinh(tB_0)|^{-1}e^{-\frac{B_0(r_1^2+r_2^2)}{4\tanh(tB_0)}}.
\end{split}
\end{equation}
As for the second term of the heat kernel \eqref{heatkernel:2}, it is sufficient to show
\begin{align}\label{est:h2}
\left|\int_{\mathbb{R}}e^{-\frac{B_0r_1r_2}{2\sinh(tB_0)}\cosh s}\left(\frac{e^{\alpha(s+i\pi)}}{e^{\frac{1}{\sigma}(s+tB_0+i\phi_+)}-1}-\frac{e^{\alpha(s-i\pi)}}{e^{\frac{1}{\sigma}(s+tB_0+i\phi_-)}-1}\right)\mathrm{d}s\right|\lesssim1,
\end{align}
where the implicit constant is independent of $t,r_1,r_2,\phi_\pm$.

\noindent Setting $\tilde{\phi}_\pm=\phi_\pm-\sigma\pi$, we have
\begin{align*}
\frac{e^{i\alpha\pi}}{e^{\frac{1}{\sigma}(s+tB_0+i\phi_+)}-1}&-\frac{e^{-i\alpha\pi}}{e^{\frac{1}{\sigma}(s+tB_0+i\phi_-)}-1}\\
&=\frac{e^{i\pi(\alpha-1)}}{e^{\frac{1}{\sigma}(s+tB_0+i\tilde{\phi}_+)}+1}-\frac{e^{-i\pi(\alpha+1)}}{e^{\frac{1}{\sigma}(s+tB_0+i\tilde{\phi}_-)}+1}
\end{align*}
and the proof of \eqref{est:h2} is thus reduced to show
\begin{equation}\label{est:term1}
\left|\int_{\mathbb{R}}e^{-\frac{B_0r_1r_2}{2\sinh(tB_0)}\cosh s}\frac{e^{s\alpha}}{e^{\frac{1}{\sigma}(s+tB_0+i\tilde{\phi}_\pm)}+1}\mathrm{d}s\right|\lesssim1.
\end{equation}
To verify \eqref{est:term1}, we observe that
\begin{align*}
&\int_{\mathbb{R}}e^{-\frac{B_0r_1r_2}{2\sinh(tB_0)}\cosh s}\frac{e^{s\alpha}}{e^{\frac{1}{\sigma}(s+tB_0+i\tilde{\phi}_\pm)}+1}\mathrm{d}s\\
&=e^{-tB_0\alpha}\int_0^{\infty}\Bigg(e^{-\frac{B_0r_1r_2}{2\sinh(tB_0)}\cosh(s-tB_0)}\frac{e^{s\alpha}}{e^{\frac{1}{\sigma}(s+i\tilde{\phi}_\pm)}+1}\\
&\qquad\quad +e^{-\frac{B_0r_1r_2}{2\sinh(tB_0)}\cosh(-s-tB_0)}\frac{e^{-s\alpha}}{e^{\frac{1}{\sigma}(-s+i\tilde{\phi}_\pm)}+1}\Bigg)\mathrm{d}s\\
&=e^{-tB_0\alpha}\Bigg(\int_0^{\infty}e^{-\frac{B_0r_1r_2}{2\sinh(tB_0)}\cosh(s+tB_0)}\left(\frac{e^{s\alpha}}{e^{\frac{1}{\sigma}(s+i\tilde{\phi}_\pm)}+1}
+\frac{e^{-s\alpha}}{e^{\frac{1}{\sigma}(-s+i\tilde{\phi}_\pm)}+1}\right)\mathrm{d}s\\
&\qquad+\int_0^{\infty}\left(e^{-\frac{B_0r_1r_2}{2\sinh(tB_0)}\cosh(s-tB_0)}
-e^{-\frac{B_0r_1r_2}{2\sinh(tB_0)}\cosh(s+tB_0)}\right)\frac{e^{s\alpha}}{e^{\frac{1}{\sigma}(s+i\tilde{\phi}_\pm)}+1}\mathrm{d}s\Bigg)
\end{align*}
and hence we are further reduced to verify
\begin{equation}\label{check:1}
\int_0^\infty\left|\frac{e^{s\alpha}}{e^{\frac{1}{\sigma}(s+i\tilde{\phi}_\pm)}+1}+\frac{e^{-s\alpha}}{e^{\frac{1}{\sigma}(-s+i\tilde{\phi}_\pm)}+1}\right|\mathrm{d}s\lesssim1
\end{equation}
and
\begin{equation}\label{check:2}
\begin{split}
\Bigg|\int_0^{\infty}\bigg(e^{-\frac{B_0r_1r_2}{2\sinh(tB_0)}\cosh(s-tB_0)}&-e^{-\frac{B_0r_1r_2}{2\sinh(tB_0)}\cosh(s+tB_0)}\bigg)\\
&\times\frac{e^{s\alpha}}{e^{\frac{1}{\sigma}(s+i\tilde{\phi}_\pm)}+1}\mathrm{d}s\Bigg|\lesssim1,
\end{split}
\end{equation}
where the implicit constant is independent of $t,r_1,r_2,\tilde{\phi}_\pm$.

For the former inequality \eqref{check:1}, we use
\begin{align*}
&\frac{e^{-s\alpha}}{1+e^{\frac{1}{\sigma}(-s+i\vartheta)}}+\frac{e^{s\alpha}}{1+e^{\frac{1}{\sigma}(s+i\vartheta)}}\\
&=\frac{\cosh(s\alpha)e^{-i\vartheta/\sigma}+\cosh((1/\sigma-\alpha)s)}{\cos(\vartheta/\sigma)+\cosh(s/\sigma)}\\
&=\frac{\cosh(s\alpha)\cos(\vartheta/\sigma)+\cosh((1/\sigma-\alpha)s)-i\sin(\vartheta/\sigma)\cosh(s\alpha)}{2\left(\cos^2(\frac{\vartheta}{2\sigma})+\sinh^2(\frac{s}{2\sigma})\right)}\\
&=\frac{2\cos^2(\frac{\vartheta}{2\sigma})\cosh(s\alpha)+\big(\cosh((1/\sigma-\alpha)s)-\cosh(s\alpha)\big)-
i\sin(\vartheta/\sigma)\cosh(s\alpha)}{2(\cos^2(\frac{\theta}{2\sigma})+\sinh^2(\frac{s}{2\sigma}))}
\end{align*}
and
\begin{align*}
&\cosh x-1\sim\frac{x^2}2,\quad \sinh x\sim x\quad \text{as}\quad x\to 0, \\
& \cosh x\sim e^x,\quad \sinh x\sim e^x\quad \text{as}\quad x\to \infty
\end{align*}
to compute (remember that $\alpha-\frac{1}{\sigma}<0$)
\begin{equation*}
\begin{split}
\int_0^\infty
&\left|\frac{\cos^2(\frac{\tilde{\phi}_\pm}{2\sigma})\cosh(s\alpha)}{\cos^2(\frac{\tilde{\phi}_\pm}{2\sigma})+\sinh^2(\frac{s}{2\sigma})}\right|\mathrm{d}s\\
&\qquad\lesssim\int_0^1\frac{|\cos(\frac{\tilde{\phi}_\pm}{2\sigma})|^2}{(\frac{s}{2\sigma})^2+(2|\cos(\frac{\tilde{\phi}_\pm}{2\sigma})|)^2}\mathrm{d}s+\int_1^\infty e^{(\alpha-\frac{1}{\sigma})s}\mathrm{d}s\lesssim1
\end{split}
\end{equation*}
and similarly,
\begin{align*}
&\int_0^\infty\left|\frac{\sin(\tilde{\phi}_\pm/\sigma)\cosh(s\alpha)}{\cos^2(\frac{\tilde{\phi}_\pm}{2\sigma})+\sinh^2(\frac{s}{2\sigma})}\right|\mathrm{d}s\lesssim1,\\
&\int_0^\infty\left|\frac{\cosh((1/\sigma-\alpha)s)-\cosh(s\alpha)}{\cos^2(\frac{\theta}{2\sigma})+\sinh^2(\frac{s}{2\sigma})}\right|\mathrm{d}s\\
&\lesssim\int_0^1\frac{\left|\frac{(1/\sigma-\alpha)^2}{2}-\frac{\alpha^2}{2}\right|s^2}{(\frac{s}{2\sigma})^2}\mathrm{d}s\\
&\qquad\qquad+\int_1^\infty\left(e^{-s\alpha}+e^{(\alpha-\frac{1}{\sigma})s}\right)\mathrm{d}s\lesssim 1.
\end{align*}
For the latter inequality \eqref{check:2}, we always have
\begin{equation*}
\left|e^{-\frac{B_0r_1r_2}{2\sinh(tB_0)}\cosh(tB_0\pm s)}\right|\leq1
\end{equation*}
and hence
\begin{align*}
\Bigg|\int_1^\infty\bigg(e^{-\frac{B_0r_1r_2}{2\sinh(tB_0)}\cosh(s-tB_0)}-&e^{-\frac{B_0r_1r_2}{2\sinh(tB_0)}\cosh(s+tB_0)}\bigg)
\frac{e^{s\alpha}}{1+e^{\frac{1}{\sigma}(s+i\tilde{\phi}_\pm)}}\mathrm{d}s\Bigg|\\
&\qquad\qquad\leq\int_1^\infty\left|\frac{2e^{s\alpha}}{1+e^{\frac{1}{\sigma}(s+i\tilde{\phi}_\pm)}}\right|\mathrm{d}s\\
&\qquad\qquad=\int_1^\infty\left|\frac{2e^{(\alpha-\frac{1}{\sigma})s}}{e^{-\frac{s}{\sigma}}+e^{i\frac{\tilde{\phi}_\pm}{\sigma}}}\right|\mathrm{d}s\\
&\qquad\qquad\leq\frac{2}{1-e^{-1/\sigma}}\int_1^\infty e^{(\alpha-\frac{1}{\sigma})s}\mathrm{d}s\lesssim_\sigma1.
\end{align*}
Now the verification of \eqref{check:2} is reduced to prove
\begin{equation*}
\int_0^1|\varphi(s)-\varphi(-s)|\left|\frac{e^{s\alpha}}{1+e^{\frac{1}{\sigma}(s+i\tilde{\phi}_\pm)}}\right|\mathrm{d}s\lesssim1,
\end{equation*}
where $\varphi$ is defined by
\begin{equation*}
\varphi(s)=e^{-\frac{B_0r_1r_2}{2\sinh(tB_0)}\cosh(s-tB_0)}.
\end{equation*}
Since
\begin{equation*}
\varphi'(\pm s)=\mp\frac{B_0r_1r_2}{2}\left(\frac{\pm\sinh s\cosh(tB_0)}{\sinh(tB_0)}-\cosh s\right)\varphi(\pm s),
\end{equation*}
we have
\begin{equation*}
|\varphi'(s)|\lesssim1,\quad \forall|s|\in(0,1)
\end{equation*}
and by the mean value theorem from elementary calculus, it follows that
\begin{equation*}
|\varphi(s)-\varphi(-s)|\lesssim s,\quad \forall s\in(0,1),
\end{equation*}
which yields
\begin{equation*}
\int_0^1|\varphi(s)-\varphi(-s)|\left|\frac{e^{s\alpha}}{1+e^{\frac{1}{\sigma}(s+i\tilde{\phi}_\pm)}}\right|\mathrm{d}s\lesssim\int_0^1\frac{se^{s\alpha}}{e^{s/\sigma}-1}\mathrm{d}s\lesssim1.
\end{equation*}
Therefore, the latter inequality \eqref{check:2} follows and so does the desired Gaussian boundedness of the heat kernel \eqref{est:heatkernel:2}.
\end{proof}

\subsection{Bernstein inequalities and square function estimates }\label{sec:BS}

In this subsection, we prove the Bernstein inequalities and the square function estimates associated with the operator $H_{\alpha,B_0,\sigma}$ with the aid of the Gaussian boundedness of the associated heat kernel \eqref{est:heatkernel:2}.
\begin{proposition}[Bernstein inequalities]
Let $\varphi(\lambda)$ be a smoothing function on $\mathbb{R}_+$ with support in $[\frac{1}{2},2]$, then it holds for any $f\in L^q(X)$ and $j\in\mathbb{Z}$
\begin{equation}\label{est:Bern}
\|\varphi(2^{-j}\sqrt{H_{\alpha,B_0,\sigma}})f\|_{L^p(X)}\lesssim2^{2j\left(\frac{1}{q}-\frac{1}{p}\right)}\|f\|_{L^q(X)},\quad 1\leq q\leq p\leq\infty.
\end{equation}
\end{proposition}
\begin{proof}
Let $\psi(x)=\varphi(\sqrt{x})$ and $\psi_e(x):=\psi(x)e^{2x}$, then $\psi_e$ is a smoothing function on $\mathbb{R}_+$ with support contained in $[\frac{1}{4},4]$ and its Fourier transform $\hat{\psi}_e$ is a Schwartz function. Since
\begin{align*}
\varphi(\sqrt{x})&=\psi(x)=e^{-2x}\psi_{e}(x)\\
&=e^{-2x}\int_{\mathbb{R}}e^{ix\cdot\xi}\hat{\psi}_e(\xi)\mathrm{d}\xi\\
 &=e^{-x}\int_{\mathbb{R}}e^{-x(1-i\xi)}\hat{\psi}_e(\xi)\mathrm{d}\xi,
\end{align*}
it follows by the functional calculus
\begin{equation*}
\varphi(\sqrt{H_{\alpha,B_0,\sigma}})=\psi(H_{\alpha,B_0,\sigma})=e^{-H_{\alpha,B_0,\sigma}}\int_{\mathbb{R}}e^{-(1-i\xi)H_{\alpha,B_0,\sigma}}\hat{\psi}_e(\xi)\mathrm{d}\xi
\end{equation*}
and thus
\begin{equation*}
\varphi(2^{-j}\sqrt{H_{\alpha,B_0,\sigma}})=\psi(2^{-2j}H_{\alpha,B_0,\sigma})=e^{-2^{-2j}H_{\alpha,B_0,\sigma}}\int_{\mathbb{R}} e^{-(1-i\xi)2^{-2j}H_{\alpha,B_0,\sigma}} \hat{\psi}_e(\xi)\mathrm{d}\xi.
\end{equation*}
Let $d_X((r_j,\theta_j),(r_k,\theta_k))$ denote the distance of two different points $(r_j,\theta_j),(r_k,\theta_k)$ on the product cone $X=(0,+\infty)_r\times\mathbb{S}_\sigma^1$, then it is easy to verify that $d_X((r_j,\theta_j),(r_k,\theta_k))\lesssim_\sigma r_1^2+r_2^2$. Hence, for three points $(r_1,\theta_1),(r_2,\theta_2),(r_3,\theta_3)\in X$, we obtain by the Gaussian bound \eqref{est:heatkernel:2} in Proposition \ref{prop:heat-Gaussbound} with $t=2^{-2j}$
\begin{align*}
&\Big|\varphi(2^{-j}\sqrt{H_{\alpha,B_0,\sigma}})(r_1,\theta_1;r_2,\theta_2)\Big|\\
&\lesssim2^{4j}\int_0^\infty\int_0^{2\sigma\pi}e^{-\frac{2^{2j}d_X((r_1,\theta_1),(r_3,\theta_3))^2}{c}}e^{-\frac{2^{2j}d_X((r_2,\theta_2),(r_3,\theta_3))^2}{c}}
\mathrm{d}\theta_3r_3dr_3\int_{\mathbb{R}}\hat{\psi}_e(\xi)\mathrm{d}\xi\\
&\lesssim2^{2j}\int_0^\infty\int_0^{2\sigma\pi}e^{-\frac{2^{2j}d_X((r_1,\theta_1),(r_3,\theta_3))^2}{c}}e^{-\frac{2^{2j}d_X((r_2,\theta_2),(r_3,\theta_3))^2}{c}}
\mathrm{d}\theta_3r_3dr_3\\
&\lesssim2^{2j}(1+2^jd_X((r_1,\theta_1),(r_2,\theta_2)))^{-N},
\end{align*}
where we use the fact
\begin{equation*}
d_X((r_1,\theta_1),(r_3,\theta_3))^2+d_X((r_2,\theta_2),(r_3,\theta_3))^2\geq\frac{1}{2}d_X((r_1,\theta_1),(r_2,\theta_2))^2,\quad \forall(r_j,\theta_j)\in X,\quad j=1,2,3
\end{equation*}
and
\begin{align*}
&\int_0^\infty\int_0^{2\sigma\pi}e^{-\frac{d_X((r_1,\theta_1),(r_3,\theta_3))^2}{c}}
e^{-\frac{d_X((r_2,\theta_2),(r_3,\theta_3))^2}{c}}\mathrm{d}\theta_3r_3\mathrm{d}r_3\\
&\lesssim e^{-\frac{d_X((r_1,\theta_1),(r_2,\theta_2))^2}{4c}}\int_0^\infty\int_0^{2\sigma\pi}
e^{-\frac{d_X((r_1,\theta_1),(r_3,\theta_3))^2}{2c}}e^{-\frac{d_X((r_2,\theta_2),(r_3,\theta_3))^2}{2c}}\mathrm{d}\theta_3r_3\mathrm{d}r_3\\
&\lesssim e^{-\frac{d_X((r_1,\theta_1),(r_2,\theta_2))^2}{2c}}\lesssim(1+d_X((r_1,\theta_1),(r_2,\theta_2)))^{-N}.
\end{align*}
Finally, the inequality \eqref{est:Bern} follows by the Young inequality.
\end{proof}
\begin{remark}
We do not need the explicit information of the distance function $d_X$ defined on the product cone $X$ since it is not hard to verify that $d_X$ is comparable to the Euclidean distance along the lateral orientation of the cone $X$.
\end{remark}
\begin{proposition}[Square function estimates]\label{prop:squarefun}
Let $\{\varphi_j\}_{j\in\mathbb{Z}}$ be a Littlewood-Paley sequence given by\eqref{LP-dp}, then for $1<p<\infty$, there exist two constants $c_p$ and $C_p$ depending on $p$ such that
\begin{equation}\label{square}
c_p\|f\|_{L^p(X)}\leq
\Big\|\Big(\sum_{j\in\mathbb{Z}}|\varphi_j(\sqrt{H_{{\alpha},B_0,\sigma}})f|^2\Big)^{1/2}\Big\|_{L^p(X)}\leq
C_p\|f\|_{L^p(X)}.
\end{equation}
\end{proposition}
\begin{proof}
In view of the Gaussian upper bound for the heat kernel \eqref{est:heatkernel:2} in Proposition \ref{prop:heat-Gaussbound}, the square function estimate \eqref{square} follows standardly from
the Rademacher functions argument in \cite{Ste70} (see also \cite{Ale04}).
\end{proof}

\section{Decay estimate for the wave equation}\label{sec:decay}

In this section, we are devoted to prove the dispersive estimate \eqref{dis:W} for the wave equation.
The key ingredient is the following Proposition about the subordination formula from \cite{MS15,DPR10}.
\begin{proposition}\label{prop:key}
If $\varphi(\lambda)\in C_c^\infty(\mathbb{R}_+)$ is a smoothing function supported in $[\frac{1}{2},2]$, then, for all $j\in\Z, t, x>0$ with $2^jt\geq 1$,  it holds
\begin{equation}\label{key}
\begin{split}
&\varphi(2^{-j}\sqrt{x})e^{it\sqrt{x}}\\
&=\rho\left(\frac{tx}{2^j},2^jt\right)+\varphi(2^{-j}\sqrt{x})\big(2^jt\big)^{\frac{1}{2}}\int_0^\infty\chi(s,2^jt)e^{-\frac{i2^jt}{4s}}e^{i2^{-j}tsx}\mathrm{d}s,
\end{split}
\end{equation}
where $\rho(s,\tau)\in\mathcal{S}(\mathbb{R}\times\mathbb{R})$ is a Schwartz function and $\chi\in C^\infty(\mathbb{R}\times\mathbb{R})$ is a smoothing function with supp$\chi(\cdot,\tau)\subseteq[\frac{1}{16},4]$ such that
\begin{equation}\label{est:chi}
\sup_{\tau\in\mathbb{R}}\big|\partial_s^\alpha\partial_\tau^\beta\chi(s,\tau)\big|\lesssim_{\alpha,\beta}(1+|s|)^{-\alpha},\quad \forall \alpha,\beta\geq0.
\end{equation}
\end{proposition}
\begin{remark}
By the spectral theorem, the subordination formula \eqref{key} leads us to represent the frequency-truncated half-wave propagator corresponding to $H_{\alpha,B_0,\sigma}$ as
\begin{equation}\label{key-operator}
\begin{split}
\varphi(2^{-j}&\sqrt{H_{\alpha,B_0,\sigma}})e^{it\sqrt{H_{\alpha,B_0,\sigma}}}=\rho\left(\frac{tH_{\alpha,B_0,\sigma}}{2^j},2^jt\right)\\
&+\varphi(2^{-j}\sqrt{H_{\alpha,B_0,\sigma}})\big(2^jt\big)^{\frac{1}{2}}\int_0^\infty\chi(s,2^jt)e^{-\frac{i2^jt}{4s}}e^{i2^{-j}tsH_{\alpha,B_0,\sigma}}\mathrm{d}s.
\end{split}
\end{equation}
For a proof of the formula \eqref{key}, we refer to \cite[Pages 32-37]{WZZ25-JDE}.
\end{remark}

\subsection{Decay estimate for the frequency-truncated half-wave propagator}

In this subsection, we obtain the following result:
\begin{proposition}\label{prop:decay-truncated}
Let $\varphi$ be a smoothing cut-off function as in \eqref{LP-dp}, then there exists some constant $C>0$ such that for $2^{-j}|t|\leq\frac{\pi}{2B_0}$ with $j\in\mathbb{Z}$,
\begin{equation}\label{est:mic-decay1}
\begin{split}
\big\|\varphi(2^{-j}\sqrt{H_{\alpha,B_0,\sigma}})&e^{it\sqrt{H_{\alpha,B_0,\sigma}}}f\big\|_{L^\infty(X)}\\
&\leq C2^{2j}\big(1+2^j|t|\big)^{-1/2}\|\varphi(2^{-j}\sqrt{H_{\alpha,B_0,\sigma}}) f\|_{L^1(X)}.
\end{split}
\end{equation}
In particular, there exists a constant $C_T>0$ depending on $T$ such that for any finite $T<+\infty$
\begin{equation}\label{est:mic-decay2}
\begin{split}
&\left\|\varphi(2^{-j}\sqrt{H_{\alpha,B_0,\sigma}})e^{it\sqrt{H_{\alpha,B_0,\sigma}}}f\right\|_{L^\infty(X)}\\
&\qquad\qquad \leq C_T2^{2j}\big(1+2^j|t|\big)^{-1/2}\|\varphi(2^{-j}\sqrt{H_{\alpha,B_0,\sigma}}) f\|_{L^1(X)},\quad 0<|t|<T.
\end{split}
\end{equation}
\end{proposition}
\begin{remark}
The $L^1\rightarrow L^\infty$ estimate \eqref{est:mic-decay2} yields immediately the desired wave dispersive estimate \eqref{dis:W} since it follows by the definition of Besov norm \eqref{Besov} that
\begin{equation*}
\begin{split}
\left\|e^{it\sqrt{H_{\alpha,B_0,\sigma}}}f\right\|_{L^\infty(X)}
&\leq\sum_{j\in\mathbb{Z}}\left\|\varphi(2^{-j}\sqrt{H_{\alpha,B_0,\sigma}})e^{it\sqrt{H_{\alpha,B_0,\sigma}}}f\right\|_{L^\infty(X)}\\
&\leq C_T|t|^{-\frac{1}{2}}\sum_{j\in\mathbb{Z}}2^{\frac{3}{2}j}\|\varphi(2^{-j}\sqrt{H_{\alpha,B_0,\sigma}})f\|_{L^1(X)}\\
&\leq C_T|t|^{-\frac{1}{2}}\|f\|_{\dot{\mathcal{B}}^{3/2}_{1,1,\sigma}(X)}.
\end{split}
\end{equation*}
\end{remark}
\begin{proof}[Proof of Proposition \ref{prop:decay-truncated}]
To obtain \eqref{est:mic-decay1}, we shall consider two cases according to $|t|2^j\geq 1$ (excluding $t=0$) and $|t|2^{j}\lesssim 1$ (including $t=0$).
We may, without loss of generality, assume $t>0$.

{\bf Case 1: $t2^{j}\lesssim 1$.} By the spectral theorem, it trivially follows that
\begin{equation*}
\|e^{it\sqrt{H_{\alpha,B_0,\sigma}}}\|_{L^2(X)\rightarrow L^2(X)}\lesssim1.
\end{equation*}
In fact, by the functional calculus, we can write for $f\in L^2(X)$
\begin{equation*}
\left(e^{it\sqrt{H_{\alpha,B_0,\sigma}}}f\right)(r,\theta)=\sum_{k\in\mathbb{Z} \atop m\in\mathbb{N}}e^{it\sqrt{\lambda_{k,m}}}c_{k,m}\tilde{V}_{k,m}(r,\theta),
\end{equation*}
where
\begin{equation*}
c_{k,m}=\int_0^\infty\int_0^{2\sigma\pi}f(r,\theta)\overline{\tilde{V}_{k,m}(r,\theta)}\mathrm{d}\theta rdr.
\end{equation*}
Then, we have (by the orthogonality between different eigenfunctions $\tilde{V}_{k,m}$)
\begin{equation*}
\begin{split}
\left\|e^{it\sqrt{H_{\alpha,B_0,\sigma}}}f\right\|_{L^2(X)}&=\left(\sum_{k\in\mathbb{Z} \atop m\in\mathbb{N}}\big|e^{it\sqrt{\lambda_{k,m}}}c_{k,m}\big|^2\right)^{1/2}\\
&=\left(\sum_{k\in\mathbb{Z} \atop m\in\mathbb{N}}\big|c_{k,m}\big|^2\right)^{1/2}\\
&=\|f\|_{L^2(X)},
\end{split}
\end{equation*}
which, together with the Bernstein inequality \eqref{est:Bern}, yields
\begin{equation*}
\begin{split}
&\left\|\varphi(2^{-j}\sqrt{H_{\alpha,B_0,\sigma}})e^{it\sqrt{H_{\alpha,B_0,\sigma}}}f\right\|_{L^\infty(X)}\\
&\qquad\qquad\lesssim2^{j}\left\|e^{it\sqrt{H_{\alpha,B_0,\sigma}}}\varphi(2^{-j}\sqrt{H_{\alpha,B_0,\sigma}})f\right\|_{L^2(X)}\\
&\qquad\qquad\lesssim2^{j}\|\varphi(2^{-j}\sqrt{H_{\alpha,B_0,\sigma}})f\|_{L^2(X)}\\
&\qquad\qquad\lesssim2^{2j}\|\varphi(2^{-j}\sqrt{H_{\alpha,B_0,\sigma}})f\|_{L^1(X)}.
\end{split}
\end{equation*}
Hence, for $2^jt\lesssim1$, we have
\begin{equation}\label{<1}
\begin{split}
&\big\|\varphi(2^{-j}\sqrt{H_{\alpha,B_0,\sigma}})e^{it\sqrt{H_{\alpha,B_0,\sigma}}}f\big\|_{L^\infty(X)}\\
&\qquad\lesssim 2^{2j}(1+2^jt)^{-N}\|\varphi(2^{-j}\sqrt{H_{\alpha,B_0,\sigma}})f\|_{L^1(X)},\quad \forall N\geq 0.
\end{split}
\end{equation}

{\bf Case 2: $2^jt\geq 1$.} In view of the representation formula \eqref{key-operator}, we are firstly required to bound by the spectral theorem and the Bernstein inequality \eqref{est:Bern} the following term
\begin{equation*}
\left\|\varphi(2^{-j}\sqrt{H_{\alpha,B_0,\sigma}})\rho\left(\frac{tH_{\alpha,B_0,\sigma}}{2^j},2^jt\right)f\right\|_{L^\infty(X)}.
\end{equation*}
Due to $\rho\in \mathcal{S}(\mathbb{R}_+\times\mathbb{R}_+)$, we see
\begin{equation*}
\left|\rho\left(\frac{t\lambda_{k,m}}{2^j},2^jt\right)\right|\lesssim(1+2^jt)^{-N},\quad \forall N\geq 0
\end{equation*}
and it thus follows by the Bernstein inequality \eqref{est:Bern} that
\begin{equation*}
\begin{split}
\Bigg\|\varphi(2^{-j}\sqrt{H_{\alpha,B_0,\sigma}})\rho&\left(\frac{tH_{\alpha,B_0,\sigma}}{2^j},2^jt\right)f\Bigg\|_{L^\infty(X)}\\
&\lesssim2^{j}\left\|\rho\left(\frac{tH_{\alpha,B_0,\sigma}}{2^j},2^jt\right)\varphi(2^{-j}\sqrt{H_{\alpha,B_0,\sigma}})f\right\|_{L^2(X)}\\
&\lesssim2^{j}(1+2^jt)^{-N}\Big\|\varphi(2^{-j}\sqrt{H_{\alpha,B_0,\sigma}})f\Big\|_{L^2(X)}\\
&\lesssim2^{2j}(1+2^jt)^{-N}\Big\|\varphi(2^{-j}\sqrt{H_{\alpha,B_0,\sigma}})f\Big\|_{L^1(X)}.
\end{split}
\end{equation*}
For the second term of \eqref{key-operator}, we shall apply the Schr\"odinger dispersive estimate \eqref{dis:S} to obtain
\begin{equation*}
\begin{split}
&\left\|\varphi(2^{-j}\sqrt{H_{\alpha,B_0,\sigma}})\big(2^jt\big)^{\frac{1}{2}}\int_0^\infty\chi(s,2^jt)e^{\frac{2^jt}{4is}}e^{i2^{-j}tsH_{\alpha,B_0,\sigma}}f\mathrm{d}s\right\|_{L^\infty(X)}\\
&\quad\lesssim\big(2^jt\big)^{\frac{1}{2}}\int_0^\infty
\chi(s,2^jt)|\sin(2^{-j}tsB_0)|^{-1}\mathrm{d}s\left\|\varphi(2^{-j}\sqrt{H_{\alpha,B_0,\sigma}})f\right\|_{L^1(X)}.
\end{split}
\end{equation*}
In view of $s\in[\frac{1}{16},4]$ (the support of $\chi$) and $\sin t\sim t$ for $t\in(0,\pi/2)$, we have for $2^{-j}t\leq\frac{\pi}{2B_0}$
\begin{equation}\label{>1}
\begin{split}
&\left\|\varphi(2^{-j}\sqrt{H_{\alpha,B_0,\sigma}})\big(2^jt\big)^{\frac{1}{2}}\int_0^\infty
\chi(s,2^jt)e^{\frac{2^jt}{4is}}e^{i2^{-j}tsH_{\alpha,B_0,\sigma}}f\mathrm{d}s\right\|_{L^\infty(X)}\\
&\qquad\quad\lesssim\big(2^jt\big)^{\frac{1}{2}}(2^{-j}t)^{-1}\int_0^\infty\chi(s,2^jt)
\mathrm{d}s\big\|\varphi(2^{-j}\sqrt{H_{\alpha,B_0,\sigma}})f\big\|_{L^1(X)}\\
&\qquad\quad\lesssim2^{2j}\big(1+2^jt\big)^{-\frac{1}{2}}\big\|\varphi(2^{-j}\sqrt{H_{\alpha,B_0,\sigma}})f\big\|_{L^1(X)}.
\end{split}
\end{equation}
Collecting the inequalities \eqref{<1} and \eqref{>1}, we obtain the desired estimate \eqref{est:mic-decay1}. To illustrate the special case \eqref{est:mic-decay2}, let us consider $0<t<T$ for any finite $T<+\infty$.
It is not difficult to see that for any finite $T>0$, there exists some integer $j_0\in\mathbb{Z}_+$ such that $2^{-j_0}T\leq\frac{\pi}{2B_0}$. On the one hand, for $j\leq j_0$, we have $2^jt\lesssim1$ and \eqref{est:mic-decay2} is covered by {\bf Case 1}. On the other, for $j\geq j_0$, if $2^jt\lesssim1$, the special case \eqref{est:mic-decay2} is still covered by {\bf Case 1}; if $2^jt\geq1$, the special case \eqref{est:mic-decay2} is covered by {\bf Case 2} since it always holds $2^{-j}t\leq\frac{\pi}{2B_0}$ for $j\geq j_0$ and $0<t\leq T$.
\end{proof}

\subsection{Strichartz estimate}\label{sec:str}

In this subsection, we present the Strichartz estimates \eqref{str:W} for the wave equation \eqref{eq:W} by using the decay estimate for the frequency-truncated half-wave propagator \eqref{est:mic-decay2}.
The proof is rather standard and relies on Keel-Tao's abstract Strichartz estimates theorem (see \cite{KT98}).
\begin{proposition}\label{prop:semi}
Let $(X,\mathcal{M},\mu)$ be a $\sigma$-finite measurable space and $U:I=[0,T]\rightarrow\mathcal{B}(L^2(X,\mathcal{M},\mu))$ be a weakly measurable map satisfying, for some constant $C>0$ perhaps depending on $T,\alpha\geq0,\sigma,h>0$,
\begin{equation}\label{md}
\begin{split}
\|U(t)\|_{L^2\rightarrow L^2}&\leq C,\quad t\in \mathbb{R},\\
\|U(t)U(s)^*f\|_{L^\infty}&\leq
Ch^{-\alpha}(h+|t-s|)^{-\sigma}\|f\|_{L^1},
\end{split}
\end{equation}
then, for any $q,p\in[1,\infty]$ satisfying
\begin{equation*}
\frac{1}{q}+\frac{\sigma}{p}\leq\frac{\sigma}{2},\quad q\geq2,\quad (q,p,\sigma)\neq(2,\infty,1),
\end{equation*}
there exists some constant $\tilde{C}>0$ only depending on $C,\sigma,q,r$ such that
\begin{equation*}
\left(\int_{I}\|U(t)f\|_{L^r}^q\mathrm{d}t\right)^{\frac{1}{q}}\leq\tilde{C}h^{-(\alpha+\sigma)(\frac{1}{2}-\frac{1}{p})+\frac{1}{q}}\|f\|_{L^2}.
\end{equation*}
\end{proposition}
\begin{proof}
The proof is analogous to the semiclassical Strichartz estimates for the Schr\"odinger equation as in \cite{KTZ07,Zwo12} (see also \cite{Zhang15} for a proof).
\end{proof}
For completeness, we give a quick sketch for the proof of the expected Strichartz estimates \eqref{str:W} in Corollary \ref{cor:S-W}.
\begin{proof}[Proof of \eqref{str:W}]
For each $k\in\mathbb{Z}$, we set
\begin{equation*}
u_k(t,\cdot)=\varphi_k(\sqrt{H_{\alpha,B_0,\sigma}})u(t,\cdot),
\end{equation*}
where $\varphi_k$ is a smoothing cut-off function given by \eqref{LP-dp} and $u(t,x)$ denotes the solution of the wave equation \eqref{eq:W}. It is easy to verify that $u_k(t,x)$ solves the Cauchy problem for the wave equation
\begin{equation*}
\begin{cases}
\partial_{tt}u_k+H_{\alpha,B_0,\sigma}u_k=0, &\\
u_k(0)=f_k,\quad \partial_tu_k(0)=g_k&
\end{cases}
\end{equation*}
with $f_k=\varphi_k(\sqrt{H_{\alpha,B_0,\sigma}})u_0$ and $g_k=\varphi_k(\sqrt{H_{\alpha,B_0,\sigma}})u_1$.
For any wave-admissible pair $(q,p)\in\Lambda_s^W$, it follows by the square-function estimate \eqref{square} and the Minkowski inequality that
\begin{equation}\label{LP}
\|u\|_{L_t^q(I;L_x^p(X))}\lesssim\left(\sum_{k\in\mathbb{Z}}\|u_k\|^2_{L_t^q(I;L_x^p(X))}\right)^{1/2},\quad I=[0,T].
\end{equation}
If we set $U(t)=e^{it\sqrt{H_{\alpha,B_0,\sigma}}}$ for the moment, then we have
\begin{equation}\label{sleq}
u_k(t,z)=\frac{U(t)+U(-t)}2f_k+\frac{U(t)-U(-t)}{2i\sqrt{H_{\alpha,B_0,\sigma}}}g_k.
\end{equation}
Setting $U_k(t):=\varphi_k(\sqrt{H_{\alpha,B_0,\sigma}})e^{it\sqrt{H_{\alpha,B_0,\sigma}}}$, we obtain by the spectral theorem
\begin{equation*}
\|U_k(t)f\|_{L^2(X)}\lesssim\|f\|_{L^2(X)}.
\end{equation*}
By the decay estimate \eqref{est:mic-decay2}, we obtain
\begin{equation*}
\begin{split}
\|U_k(t)U_k^*(s)f\|_{L^\infty(X)}&=\|U_k(t-s)f\|_{L^\infty(X)}\\
&\leq C_T2^{\frac{3}{2}k}\big(2^{-k}+|t-s|\big)^{-\frac{1}{2}}\|f\|_{L^1(X)},
\end{split}
\end{equation*}
which is exactly the inequality \eqref{md} specialized to $U(t)=U_k(t)$ and $\alpha=3/2,\sigma=1/2,h=2^{-k}$. Hence, it follows from Proposition \ref{prop:semi} that
\begin{equation*}
\|U_k(t)f\|_{L_t^q(I;L_x^p(X))}\leq C_T2^{[2(\frac{1}{2}-\frac{1}{p})-\frac{1}{q}]k}\|f\|_{L^2(X)},
\end{equation*}
which, together with \eqref{LP} and \eqref{sleq}, gives the desired Strichartz estimate \eqref{str:W} since $s=1-\frac{2}{p}-\frac{1}{q}$ in \eqref{gap}.
\end{proof}

\begin{center}

\end{center}

\end{document}